\documentclass{article}

\usepackage[in]{fullpage}  
\usepackage{natbib} 
\usepackage{mathpazo} 

\usepackage[utf8]{inputenc} %
\usepackage[T1]{fontenc}    %
\usepackage{hyperref}       %
\usepackage{url}            %
\usepackage{booktabs}       %
\usepackage{amsfonts}       %
\usepackage{nicefrac}       %
\usepackage{microtype}      %
\usepackage{lipsum}         %
\usepackage{todonotes}
\usepackage{amsmath}
\usepackage{amssymb}
\usepackage{amsthm}

\usepackage{tcolorbox}
\usepackage{wrapfig}

\usepackage{url}            %
\usepackage{booktabs}       %
\usepackage{multirow}    
\usepackage{amsfonts}       %
\usepackage{nicefrac}       %
\usepackage{microtype}      %
\usepackage{natbib}
\usepackage{enumerate}
\usepackage{hhline}
\usepackage{makecell}
\usepackage{pifont}

\usepackage{graphicx} %
\usepackage{subfigure}
\usepackage{caption}
\usepackage{amsmath}
\usepackage{amsthm}
\usepackage{amssymb}
\usepackage{tikz}
\usepackage{xcolor}
\usetikzlibrary{arrows}

\usepackage{mathrsfs}

\usepackage{algorithm}
\usepackage{algorithmic}
\usepackage{hyperref}
\usepackage{bm}

\newtheorem{thm}{Theorem}
\newtheorem{lem}{Lemma}

\newtheorem{prop}{Proposition}

\usepackage[capitalize,noabbrev]{cleveref}
\crefname{thm}{Theorem}{Theorems}
\crefname{lem}{Lemma}{Lemmas}
\crefname{cor}{Corollary}{Corollaries}
\crefname{prop}{Proposition}{Propositions}
\crefname{asmp}{Assumption}{Assumptions}
\crefname{defn}{Definition}{Definitions}
\crefname{oracle}{Oracle}{Oracles}
\crefname{fact}{Fact}{Facts}
\crefname{conj}{Conjecture}{Conjectures}
\crefname{rem}{Remark}{Remarks}
\crefname{example}{Example}{Examples}
\crefname{condition}{Condition}{Conditions}
\crefname{exercise}{Exercise}{Exercises}
\crefname{algorithm}{Algorithm}{Algorithms}
\crefname{table}{Table}{Tables}
\crefname{figure}{Figure}{Figures}
\crefname{section}{Section}{Sections}
\crefname{subsection}{Section}{Sections}
\crefname{appendix}{Appendix}{Appendices}
\crefname{message}{Message}{Messages}

\definecolor{red}{rgb}{1, 0, 0}

\definecolor{green}{rgb}{0, 1, 0}

\definecolor{blue}{rgb}{0, 0, 1}
\newcommand{\BLUE}[1]{{\color{blue} #1}}

\definecolor{orange}{rgb}{1, 0.4, 0.0}

\input{packages/math_commands}

\usepackage{soul}

\usepackage{framed} %
\colorlet{shadecolor}{orange!15}

\AtBeginDocument{\setlength\abovedisplayskip{5pt}}
\AtBeginDocument{\setlength\belowdisplayskip{5pt}}

\newcommand\nnfootnote[1]{%
  \begin{NoHyper}
  \renewcommand\thefootnote{}\footnote{#1}%
  \addtocounter{footnote}{-1}%
  \end{NoHyper}
}

\newcommand{\Ex}{\mathbb{E}}

\newcommand{\pproj}{p_{\operatorname{proj}}}

\hypersetup{
    colorlinks=true,
    citecolor=blue,
    linkcolor=blue,
}

\usepackage{listings} %

\definecolor{codegreen}{rgb}{0,0.6,0}
\definecolor{codegray}{rgb}{0.5,0.5,0.5}
\definecolor{codepurple}{rgb}{0.58,0,0.82}
\definecolor{codeblue}{rgb}{0,0,1}
\definecolor{backcolour}{rgb}{0.95,0.95,0.92}
\definecolor{key-color}{rgb}{0.8, 0.47, 0.196}

\usepackage{enumitem} %

\lstdefinestyle{mystyle}{
    backgroundcolor=\color{backcolour},   
    commentstyle=\color{codegreen},
    numberstyle=\tiny\color{codegray},
    stringstyle=\color{codepurple},
    basicstyle=\ttfamily\footnotesize,
    breakatwhitespace=false,         
    breaklines=true,                 
    captionpos=b,                    
    keepspaces=true,                 
    numbers=left,                    
    numbersep=5pt,                  
    showspaces=false,                
    showstringspaces=false,
    showtabs=false,                  
    tabsize=2,
    language=Python,
    emph={lm},
    emphstyle={\color{blue}},
    classoffset=1, %
    otherkeywords={sum},
    morekeywords={rm, mean},
    keywordstyle=\color{codegreen},
    classoffset=0,
}
\title{Finite Horizon Optimization: Framework and Applications $^*$}

\makeatletter
\def\@fnsymbol#1{\ensuremath{\ifcase#1\or *\or \dagger\or \ddagger\or
   \mathsection\or \sharp\or \Diamond\or \mathparagraph\or \|\or
   \or \ddagger\ddagger \else\@ctrerr\fi}}
\makeatother

\author{%
  Yushun Zhang$^{12}$, Dmitry Rybin$^{12}$, 
    Zhi-Quan Luo$^{12\dagger}$
   \\
   \\
   $^1$The Chinese University of Hong Kong, Shenzhen, China \\
  $^2$Shenzhen Research Institute of Big Data 
  \\
  \texttt{\{yushunzhang,dmitryrybin\}@link.cuhk.edu.cn, luozq@cuhk.edu.cn} \\
}

\date{}

\begin{document}

\maketitle
\nnfootnote{$*$: This is a preliminary report of an ongoing research.}
\nnfootnote{$\dagger$: Correspondence author.}
\begin{abstract}

 In modern engineering scenarios, there is often a strict upper bound on the number of algorithm iterations that can be performed within a given time limit. This raises the question of optimal algorithmic configuration for a fixed and finite iteration budget. In this work, we introduce the framework of {\it finite horizon optimization}, which focuses on optimizing the algorithm performance under a strict iteration budget $T$. We apply this framework to linear programming (LP) and propose Finite Horizon stepsize rule for the primal-dual method.
 The main challenge in the stepsize design is controlling the singular values of $T$ cumulative product of non-symmetric matrices, which appears to be a highly nonconvex problem, and there are very few helpful tools. Fortunately, in the special case of the primal-dual method, we find that the optimal stepsize design problem admits hidden convexity, and we propose a convex semidefinite programming (SDP) reformulation. This SDP only involves matrix constraints of size \( 4 \times 4 \) and can be solved efficiently in negligible time. Theoretical acceleration guarantee is also provided at the pre-fixed $T$-th iteration, but with no asymptotic guarantee. On more than 90 real-world LP instances,  Finite Horizon stepsize rule reaches an average 3.9$\times$ speed-up over the optimal constant stepsize, saving 75\% wall-clock time. Our numerical results reveal substantial room for improvement when we abandon asymptotic guarantees, and instead focus on the performance under finite horizon. We highlight that the benefits are not merely theoretical - they translate directly into computational speed-up on real-world problems.
\footnote{Our code is available at \url{https://github.com/zyushun/Finite-Horizon-Stepsize-Rule}.}

\end{abstract}

\section{Introduction}
\label{sec_intro}

In the field of optimization, iterative algorithms are often equipped with tunable hyperparameters. For example, Gradient Descent (GD) has a hyperparameter called stepsize, which determines the magnitude of the adjustment made in each iteration \citep{cauchy1847methode}. These hyperparameters significantly affect the algorithm's performance, influencing whether it diverges, converges, or how fast it converges.   Theoretical analysis is crucial for guiding the proper selection of hyperparameters. The mainstream theoretical analysis guides the selection of hyperparameters through the following procedure. First, we fix the class of functions and fix the rule for selecting a category of hyperparameters (e.g., constant stepsize); second, we study the decay behavior of an optimality residue versus the iteration number $T$ through asymptotic analysis, i.e., we aim to guarantee the good performance as $T \rightarrow \infty$.   
This methodology leads to many well-known results. For instance,  
the optimal constant stepsize for GD is $2/ (\mu + L) $ for $\mu$-strongly convex $L$-Lipschitz smooth functions \citep{nesterov2018lectures}. Such theoretical insights have significantly deepened our understanding of iterative algorithms and have guided engineers for decades.

\begin{figure}[t]
  \centering
  \subfigure[Traditional optimization theory v.s. modern applications such as auto-vehicle]{\includegraphics[width=0.65\textwidth]{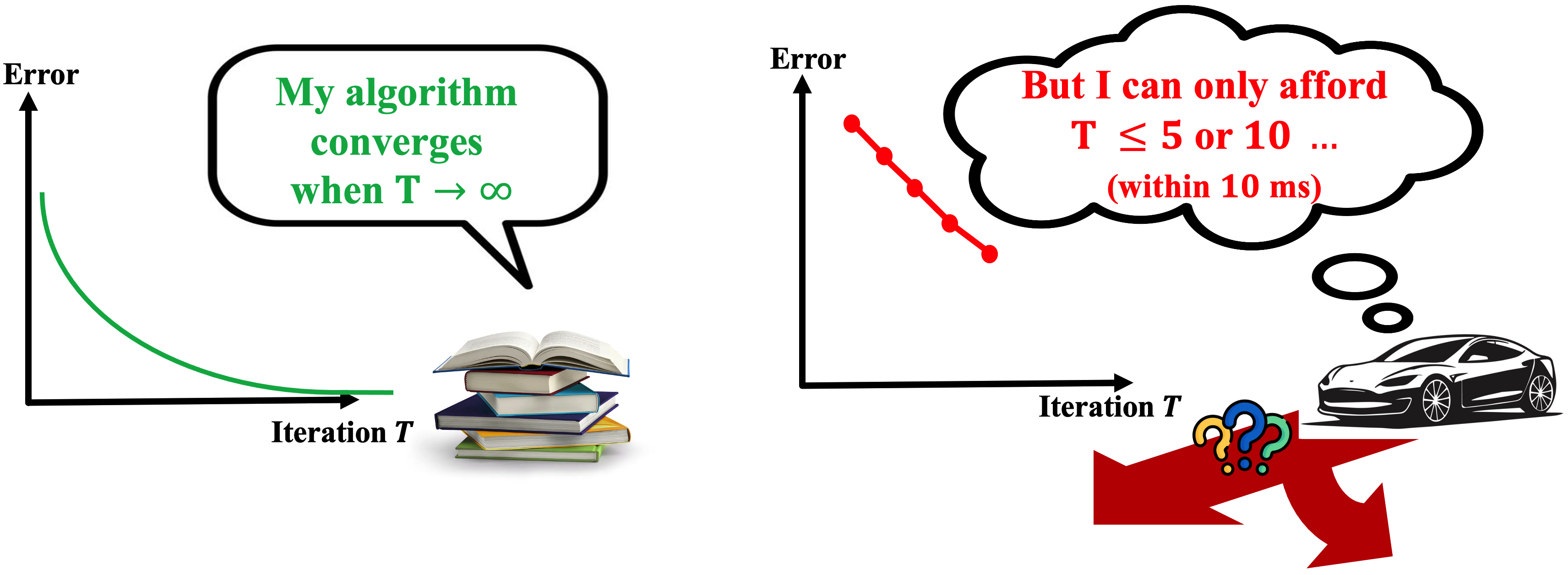}}
\subfigure[Finite Horizon stepsize rule v.s. constant stepsize]{\includegraphics[width=0.30\textwidth]{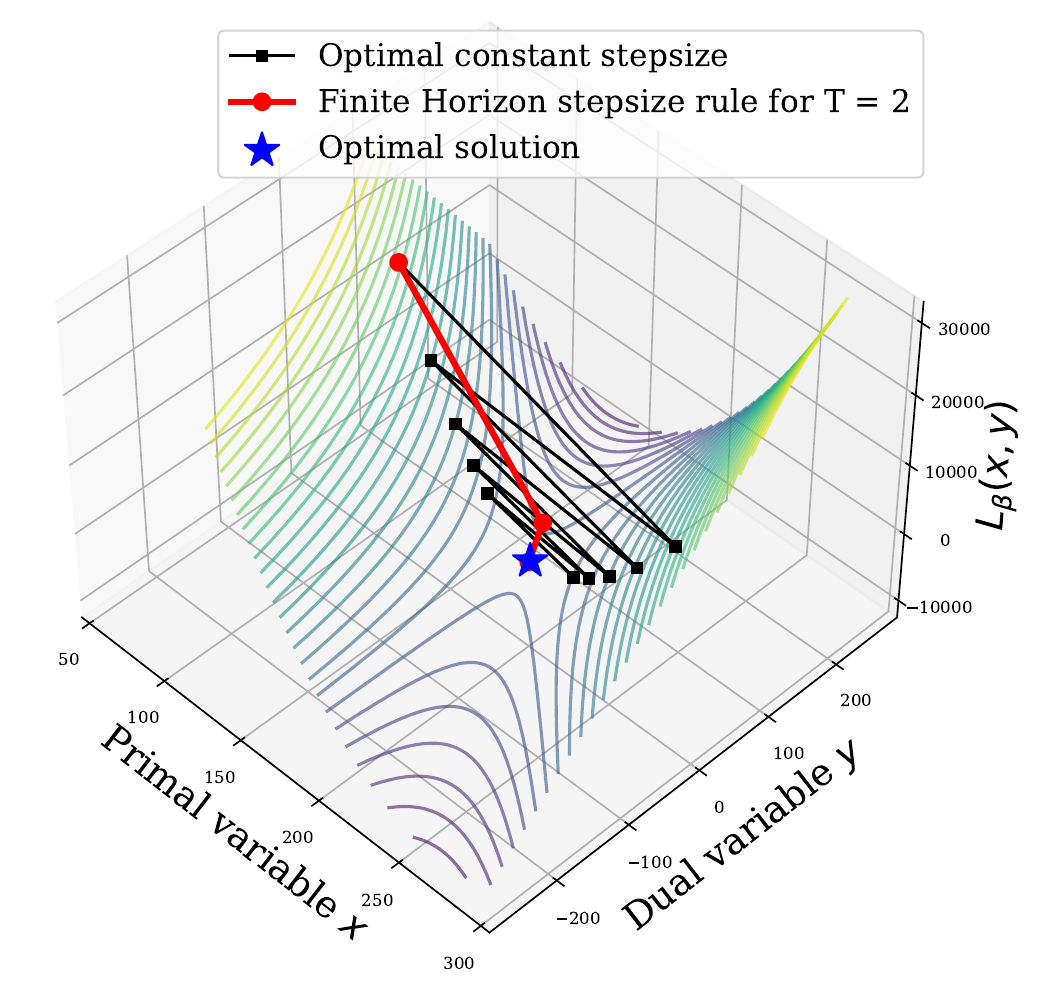}}
  \caption{ {\bf (a):} The growing mismatch between traditional optimization theory and modern applications such as auto-vehicle: classical theory focuses on $T\rightarrow \infty$, while many applications can only afford small $T$. {\bf (b):} The trajectories of Finite Horizon stepsize rule and the optimal constant stepsize on the augmented Lagrangian function of a simple linear programming (LP). The Finite Horizon stepsize rule for $T =2$ reaches the optimal solution (the saddle) in 2 steps, while the optimal constant stepsize cannot. }
\label{fig_intro}
\vspace{-0.3cm}
\end{figure}

Despite the long-standing and far-reaching impact of these theories,  modern real-world applications are revealing new challenges.  
In particular, there is often a strict upper bound on the number of algorithm iterations that can be performed within a given time limit. For instance,
in 5G or 6G wireless communication systems, iterative algorithms must operate at high throughput to ensure low-latency communications among users \citep{latva2019key}. Typically, the iteration budget for the channel decoder algorithms is restricted to fewer than {\bf 5 iterations}   (see Table 1 in \citep{sy2023optimization})
or  {\bf 15 iterations} (see Table IV in \citep{ferraz2021survey}). A comparable situation arises in modern power systems,  where the average running time for iterative algorithms is {\bf 0.062 seconds}
with {\bf 12 iterations} \citep{tang2017real}.
Similarly,  for autonomous vehicles, real-time decision-making problems are often formulated as linear or quadratic programming, and the solvers must return solutions within {\bf 10 milliseconds (ms)} or {\bf 20 ms} \citep{chen2012three,li2023real} to generate a safe and efficient path for the next 10 seconds.

The above applications reveal a growing gap between classical theory and practical scenarios:
classical theory pursues optimal performance guarantees when the iteration number $T$ approaches infinity, and the resulting designs may not be well-suited in scenarios with a small finite $T$ (see Figure~\ref{fig_intro} (a) as an illustration).
This raises the new question of optimal hyperparameter design for a fixed and finite iteration budget.  

In this work, we introduce the framework of {\it finite horizon optimization}, which focus on optimizing the algorithmic behaviors under a strict iteration budget.  We will review some existing approaches that fall within this finite horizon framework, and also introduce several practical problems that could be improved using this framework. In the main body of this work,  we will apply the finite horizon framework to linear programming  (LP) \citep{kantorovich1960mathematical,schrijver1998theory,dantzig2002linear,luenberger1984linear,boyd2004convex}.   We study LP because LP (and its variants) frequently emerge in engineering scenarios where solvers must operate under strict iteration limit (some examples are provided later in Section~\ref{sec_general_framework}). We will primarily work on a most simple yet fundamental algorithm to solve LP, namely, the primal-dual method. We aim to design a new stepsize rule for the primal-dual method,  tailored for the scenarios where the total iteration budget $T$ is fixed and finite.

Perhaps a bit surprisingly, we find that: if we abandon the pursuit of asymptotic performance and instead focus on how the algorithm performs within a finite number of iterations, the primal-dual method can be substantially accelerated by merely changing its stepsize design. 
We provide an illustration in Figure \ref{fig_intro} (b). This figure visualizes the augmented Lagrangian function of a simple LP, and the saddle point of this function is the optimal primal-dual solution to the LP. When the stepsize rule is particularly designed for a total iteration budget $T =2$, we find that the primal-dual method can reach the saddle point in 2 steps, while the conventional constant stepsize cannot solve the problem until $>200$ steps \footnote{We provide the detailed curves of optimality gap and the stepsize rules in Appendix \ref{appendix_more_results}.}.  Note that 
such acceleration is reached without changing the algorithm, e.g., storing extra building blocks like momentum \citep{nesterov2018lectures}. We summarize our contribution as follows.

\begin{itemize}[topsep=1pt,parsep=1pt,partopsep=1pt, leftmargin=*]
\item We introduce the framework of {\it finite horizon optimization}, which focuses on optimizing algorithmic performance under a strict iteration budget. We review existing approaches that fall within this framework and present several practical problems that could benefit from this methodology.
\item We introduce Finite Horizon stepsize rule for the primal-dual method for solving LP, especially for the scenarios where the total iteration  budget \( T \) is fixed and finite. The main challenge in designing this stepsize rule is controlling the singular values of $T$ cumulative product of {\it non-symmetric} matrices, which appears to be a highly nonconvex problem, and there is very few helpful tool. In the special case of the primal-dual method, we find that
the optimal stepsize design problem admits hidden convexity and we propose a convex semidefinite programming (SDP) reformulation. This reformulation involves only matrix constraints of size \( 4 \times 4 \) and can be solved efficiently by SDP solvers in negligible computation time. A theoretical acceleration guarantee is also provided at the pre-fixed $T$-th iteration, but with no asymptotic guarantee.

\item We numerically verify the effectiveness of Finite Horizon stepsize rule on real-world LP datasets.  On the \texttt{Netlib} LP benchmark with more than 90 real-world instances, our stepsize rule archives an average 3.9$\times$ speed-up over the optimal constant stepsize, saving 75\% wall-clock time to achieve the same level of precision. Our results reveal substantial room for improvement when we abandon asymptotic guarantees, and instead focus on the performance under finite horizon. We highlight that the benefits are shown not only in worst-case theory, but also manifest immediately in real-world problems.

\end{itemize}

\subsection{Notations}
\begin{itemize}[topsep=1pt,parsep=1pt,partopsep=1pt, leftmargin=*]
\item {\bf Matrix-related notations.} Given a matrix $X \in \mathbb{R}^{m\times n}$, we denote $X^\top$ as the transpose of $X$. Similarly, for $X \in \mathbb{C}^{m\times n}$, we denote $X^H$ as  hermitian of $X$. 
We write  the singular value decomposition (SVD) of $X$ as  $X = U \Sigma V^{\top} = \sum_{i = 1}^m \sigma_i u_i v_i^{\top}$, where $U \in \mathbb{R}^{m \times m}, V \in \mathbb{R}^{n \times n}$ are orthonormal matrices and $u_i \in \mathbb{R}^m, v_i \in \mathbb{R}^n$  are the $i$-th column of $U$ and $V$, respectively;   $\Sigma$ is the rectangular diagonal matrix with diagonal entries $\sigma_1 \geq \sigma_2 \cdots \geq \sigma_m \geq 0$. We also use $\sigma(X)$ to denote the set of singular values of $X$. 
We denote $\|X\|_{\operatorname{op}}$ as the spectral norm (a.k.a., operator norm) of $X$, i.e., the square root of the maximal eigenvalue of $X^{\top}X$, or equivalently, the maximal singular value $\sigma_1$ of $X$. When $X$ is a square matrix with $m= n$,   we use $\lambda(X)$ to denote the set of eigenvalues of $X$.  we denote $\rho(X):=\max \left\{\left|\lambda(X)\right|\right\}$ as the spectral radius of $X$. When $X$ is positive semi-definite (PSD), we use $\kappa = \frac{\lambda_{\max}(X)}{\lambda_{\min}^{+}(X)}$ to denote the condition number of $X$, where $\lambda_{\max}(X)$, $\lambda_{\min}^{+}(X)$ are the maximal eigenvalue and the smallest non-zero eigenvalue, respectively. We use $\operatorname{null}(X)$ to denote the null space of $X$.
 We use $I_{n\times n}$ and $0_{n \times n}$ to denote the identity matrix and the zero matrix of size $n \times n$.

\item {\bf Algorithm-related notations.} We say an algorithm has complexity $\tilde{\mathcal{O}}(C)$ or $\mathcal{O}(C \log (1 / \epsilon))$ if it takes (at most) $\mathcal{O}(C \log (1 / \epsilon))$ iterations to achieve error $\frac{ \text{dist}(z, \mathcal{Z}^*)}{\text{dist}(z^0, \mathcal{Z}^*)} \leq \epsilon$, where $\text{dist}(\cdot, \cdot)$ denotes the Euclidean distance, $z$ is the decision variable, $z^0$ is the initial point, $\mathcal{Z}^*$ is the solution set. Similarly, we say an algorithm has complexity $\tilde{\Omega}(C)$ or $\Omega(C \log (1 / \epsilon))$ if it takes (at least) $\Omega(C \log (1 / \epsilon))$ iterations to achieve error $\frac{ \text{dist}(z, \mathcal{Z}^*)}{\text{dist}(z^0, \mathcal{Z}^*)} \leq \epsilon$.

\item {\bf Other notations.} We use $[n]$ to denote the index set $\{1,2,\cdots ,n\}$.  For a vector $x \in \mathbb{R}^n$,  $\|x\|_2$ denotes the Euclidean norm of $x$, $x^{+}$ and $x^{-}$ denote the vector of positive parts and negative parts of $x$, respectively. That is: the components of $x^{+}$ and $x^{-}$ are $\left(x^{+}\right)_i=\max \left\{x_i, 0\right\}$ and $\left(x^{-}\right)_i=\max \left\{-x_i, 0\right\}$ for $i \in[n]$. 
\end{itemize}

\subsection{General Framework of Finite Horizon Optimization}
\label{sec_general_framework}

We now state the general framework for finite horizon optimization. Suppose we are given a function class $\mathcal{F}$ and an iterative algorithm $\mathcal{A}(x,\theta)$, where $x \in \mathcal{X}$ is the optimization variable and $\theta \in \Theta$ 
is the algorithmic configuration (e.g., hyperparameters, or building components). We define the performance measurement $\epsilon(x)$ to be a nonnegative function of $x$. We aim to adjust the algorithmic configuration $\theta$ to optimize the worst-case performance of the algorithm after $T$ iterations, where $T$ is pre-determined. The general formulation of finite horizon optimization is shown below.

\begin{equation}
    \label{eq_general_formulation}
    \begin{aligned}
        \min_{\theta \in \Theta} & \max_{f \in \mathcal{F}, x_0 \in \mathcal{X}} \epsilon(x_T) \\
        \text{s.t. } & x_T = \mathcal{A}(x_{T-1}, \theta), \\
                     & x_{T-1} = \mathcal{A}(x_{T-2}, \theta), \\
                     & \vdots \\
                     & x_1 = \mathcal{A}(x_{0}, \theta).
    \end{aligned}
\end{equation}

One may also easily generalize this formulation to the ``average case in $\mathcal{F}$" instead of ``worst case in $\mathcal{F}$", or  ``average performance within $T$ iterations" instead of ``the final performance at the $T$-th step".  We now provide some concrete examples that belong to \eqref{eq_general_formulation} ({\bf Example 1, 2}), as well as some examples that can be fit into \eqref{eq_general_formulation} in the future ({\bf Example 3, 4, 5}).

\paragraph{Example 1: stepsize design of GD for unconstrained minimization.} Formulation \eqref{eq_general_formulation} can be applied to select the optimal stepsize rules of GD after $T$ iterations, where $T$ is pre-determined. The stepsize selection problem can be modeled as follows.

\begin{equation}
    \label{eq_example_gd}
    \begin{aligned}
        \min_{\eta_0,\cdots, \eta_{T-1}} & \max_{f \in \mathcal{F}, x_0 \in \mathbb{R}^d} \epsilon(x_T) \\
        \text{s.t. } & x_T = x_{T-1} - \eta_{T-1} \nabla f(x_{T-1}),  \\
                     & x_{T-1} = x_{T-2} - \eta_{T-2} \nabla f(x_{T-2}),   \\
                     & \vdots \\
                     & x_1 = x_{0} - \eta_{0} \nabla f(x_{0}).  
    \end{aligned}
\end{equation}

Solving Problem \eqref{eq_example_gd} for a general $\mathcal{F}$ remains challenging. Fortunately, there have been some exciting breakthroughs for certain special classes of $\mathcal{F}$.
The pioneer work \citep{young1953richardson} solved \eqref{eq_example_gd} when $\mathcal{F}$ is the $\mu$-strongly-convex $L$-smooth {\it quadratic} functions and $\epsilon(x) := \|x - x^*\|_2^2 $ denotes the distance to optimal solution. Perhaps a bit surprisingly, the optimal $(\eta_0,\cdots, \eta_{T-1})$ in this case has the following closed-form solution.

\begin{equation}
\label{eq_young_stepsize}
    \eta_t=2\left((L+\mu)+(L-\mu) \cos \left(\left(t-\frac{1}{2}\right) \frac{\pi}{T}\right)\right)^{-1}, \quad t = 0,\cdots, T-1.
\end{equation}

Eq. \eqref{eq_young_stepsize} is related to the inverse of the roots of the $T$-th order Chebyshev polynomial. We will refer to  \eqref{eq_young_stepsize} as Young's stepsize. The derivation of Young's stepsize relies on a special connection between quadratic minimization and the minimax polynomial theory. The proof is presented in recent works \citep{altschuler2018greed,pedregosa2021residual,d2021acceleration}.   With Young's stepsize, the optimal value of $\epsilon(x_T)$ satisfies $\epsilon(x_T) \leq 2\left(1 - \frac{2}{ \sqrt{\kappa} + 1}\right)^{T}$, which can be translated into the complexity of $\tilde{\mathcal{O}}(\sqrt{\kappa})$. Note that the square-root dependency on $\kappa$ matches the lower bound of first-order methods \citep{nesterov2018lectures}, which achieves acceleration over the complexity lower bound of constant stepsize $\tilde{\Omega}(\kappa)$.
 It is worth mentioning that Young's stepsize is the first non-constant stepsize rule that achieves acceleration over the constant stepsize. It is also the first stepsize rule that helps vanilla GD achieve the optimal dependency on $\kappa$ without modifying the algorithm (e.g., introducing extra building blocks such as momentum \citep{polyak1964some,nesterov2018lectures}).

Recently, researchers have tried to design better stepsize rule of GD beyond the quadratic case. Now let $\mathcal{F}$ be the generic $L$-smooth convex functions and $\epsilon(x_t) :=  f\left(x_t\right)-\inf f $ and we still try to solve \eqref{eq_example_gd}. 
Unfortunately, the optimal $(\eta_0,\cdots, \eta_{T-1})$ does not admit a closed-form analytic expression in this case. Rather,   \citet{das2024branch} pointed out that \eqref{eq_example_gd}  can be reformulated as a nonconvex QCQP and proposed to solve it via branch-and-bound methods. They then discovered the optimal stepsize rule for $T \leq 25$.  Similar stepsize rule for $T \leq  127$  are also discovered by  ``brute force searching" \citep{grimmer2023accelerated}. Although these stepsize rules do not reach the lower bound complexity like Young's stepsize, they still achieve highly non-trivial improvement over constant stepsize.

We emphasize that the complexity of the aforementioned stepsize rules only holds at $T$-th step but not at any other iterations, so it does not have asymptotic guarantees like constant stepsize. 
These stepsize rules serve as examples that: GD can be substantially accelerated if we abandon the asymptotic guarantees and focus on the performance within $T$ iterations.

\paragraph{Example 2: Algorithm Unrolling (AU).} AU is a data-driven approach that employs a finite-depth neural network to learn the behavior of a specific algorithm after an infinite number of iterations (e.g., \citep{gregor2010learning,sun2016deep,sun2018learning,yang2018admm,adler2018learned,monga2021algorithm,chen2022learning,li2024pdhg}). We now interpolate AU under the framework of finite horizon optimization. We denote the $t$-th layer of  a neural network as $\varphi_t(x,W_t)$, where  $x$ is the input,  $W_t$ is the trainable weight, $t =0, \cdots ,T-1$. Given a neural network architecture with $T$ total layers, AU aims to train this neural network to solve the following problem. 

\begin{equation}
    \label{eq_example_unrolling}
    \begin{aligned}
        \min_{W_0,\dots, W_{T-1}} & \sum_{f\in \mathcal{F}} \sum_{x_0 \in \mathcal{X}} \|x_T - x^*(f,x_0) \|_2^2 \\
        \text{s.t. } & x_T = \varphi_{T-1}(x_{T-1}, W_{T-1}), \\
                     & x_{T-1} = \varphi_{T-2}(x_{T-2}, W_{T-2}), \\
                     & \vdots \\
                     & x_1 = \varphi_{0}(x_{0}, W_0), 
    \end{aligned}
\end{equation}

where $x^*(f,x_0) = \lim_{n \rightarrow \infty} \underbrace{\mathcal{A} (\mathcal{A}(\cdots \mathcal{A}(f,x_0)))}_{n \text{ iterations}}$ denotes the output of the target algorithm on $f$ and initialization $x_0$.
With an appropriately designed neural network architecture, a diverse collection of training instances from $\mathcal{F}$, and sufficient training of the weights $\{W_t\}_{t=0}^{T-1}$, the trained neural network can achieve performance comparable to $\mathcal{A}$ on test problem instances, all while significantly reducing computational costs. For instance,  a neural network with $T =4$ layers can perform on par with $\geq 1000$ iterations of Primal-Dual Hybrid Gradient algorithm (PDHG) on real-world LP benchmark, which helps reduce $\geq 45\%$ wall-clock running time \citep{li2024pdhg}. 

From the perspective of finite horizon optimization, 
the objective of  AU \eqref{eq_example_unrolling} is closely related to the objective in \eqref{eq_general_formulation}:
both \eqref{eq_general_formulation} and \eqref{eq_example_unrolling} aims to accelerate certain algorithm under finite $T$, up to the minor difference that  ``worst case in $\mathcal{F}$" in \eqref{eq_general_formulation} is changed to ``average case in $\mathcal{F}$" in \eqref{eq_example_unrolling}.  We believe AU serves as an example that traditional algorithms can be substantially improved under the finite horizon framework, e.g., such as through re-parametrization using finite-depth neural networks.

In the sequel, we provide some practical engineering problems where there is a strict upper bound on the total iteration number $T$. We believe there is substantial room for improvement if we re-design the relevant algorithms under the finite horizon framework \eqref{eq_general_formulation}.

\paragraph{Example 3: Weighted Sum Mean Square Error Minimization (WMMSE).}  WMMSE \citep{shi2011iteratively} is a popular optimization framework for efficient resource allocation in multiple-input multiple-output (MIMO) wireless communication systems. The goal of WMMSE is to maximize the weighted sum rate (WSR) of the system, which is a common objective in wireless networks to balance fairness and throughput among users. However, directly optimizing the WSR is often challenging due to its non-convex nature. WMMSE addresses this by introducing an additional weight matrix and lifting the problem into a higher dimensional space, where the problem has a more tractable form. The objective of WMMSE is formulated as follows.

\begin{equation}
    \label{eq_example_wmmse}
    \begin{aligned}
        \min_{\{W_i,V_i, U_i\}_{i = 1}^I} &\sum_{i=1}^I \alpha_i\left(\operatorname{Tr}\left(W_i E(U_i, V_i)\right)-\log \operatorname{det}\left(W_i\right)\right) \\
        \text{s.t. } & \sum_{i=1}^I \operatorname{Tr}\left(V_i V_i^H\right) \leq P_{\max }, 
    \end{aligned}
\end{equation}

where $W_i \in \mathbb{R}^{d_i \times d_i}$ is a positive-semidefinite weight matrix for receiver $i$; $V_i \in \mathbb{C}^{M \times d_i}$  and $U_i \in \mathbb{C}^{N_i \times d_i}$ is the beamformer and decoder for receiver $i$, respectively; $M$ and $N_i$ are the number transmit antennas of the base station and receive antennas of user $i$, respectively; $\alpha_i$ is the priority of user $i$ in the system; $E(U_i, V_i)$ is a quadratic function in  $U_i$ and $V_i$; $P_{\max }$ is the power budget of transmitter. More descriptions can be seen in \citep{shi2011iteratively}. Problem \eqref{eq_example_wmmse} is then solved iteratively by the WMMSE algorithm, which cyclically update $\{W_i,V_i, U_i\}_{i = 1}^I$  in a block-coordinate-descent fashion. The authors established the convergence guarantee of the WMMSE algorithm by leveraging the observation that the optimization problem in \eqref{eq_example_wmmse} is convex with respect to each individual variable when the other variables are held fixed.

\begin{figure}[t]
  \vspace{-1.5cm}
  \centering
  \subfigure[MMSE equalization in wireless communication system]{\includegraphics[width=0.25\textwidth]{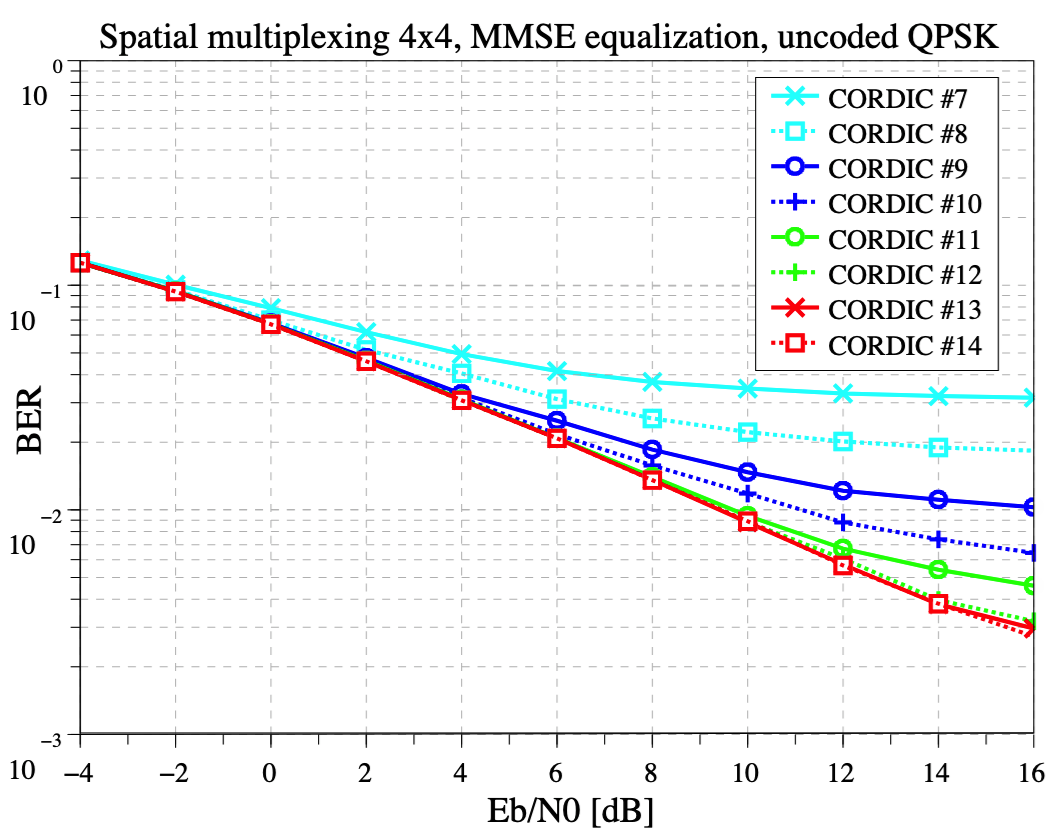}}
\subfigure[Real-time OPF problem in power system]{\includegraphics[width=0.38\textwidth]{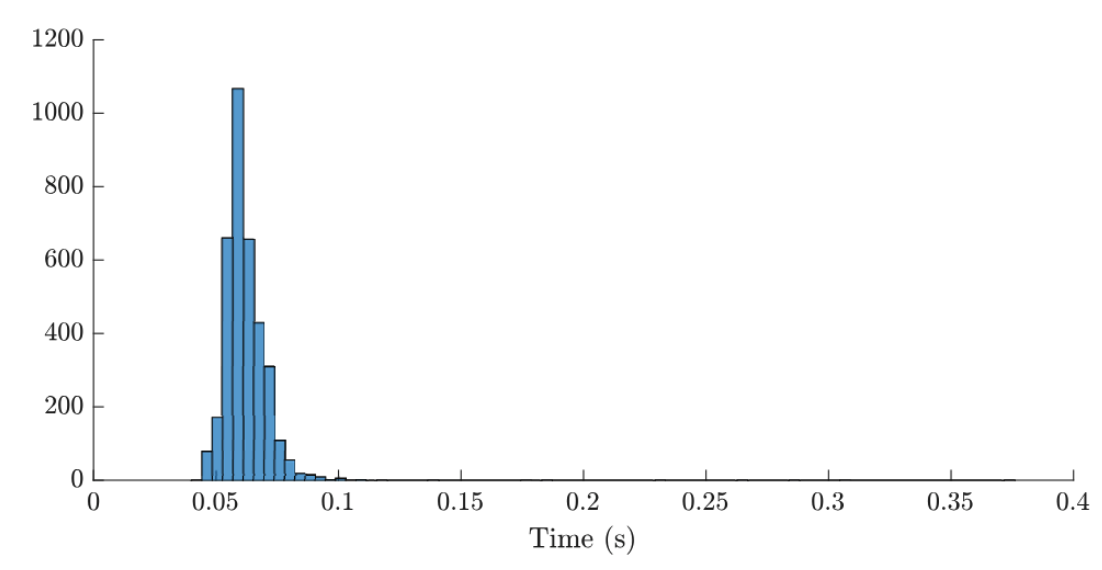}}
\subfigure[Path planning in auto-vehicles]{\includegraphics[width=0.25\textwidth]{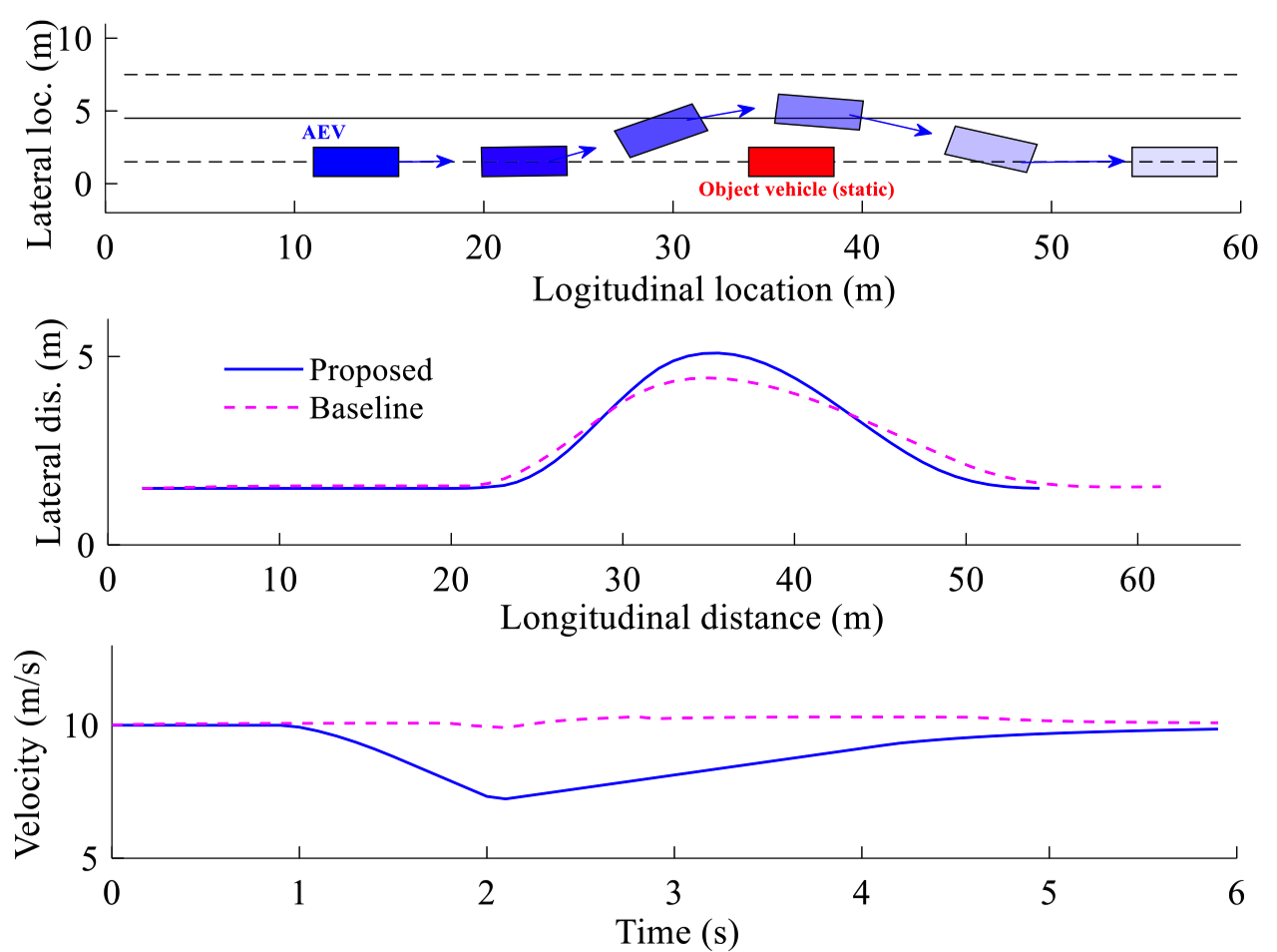}}
  \caption{Some practical applications with strict iteration budget for solvers. {\bf (a):} Figure 4 in \citep{boher2008fpga}. Performance of  MMSE equalization algorithms in MIMO  wireless communication systems after $\# x$ number of iterations. The maximal affordable iteration number is 14.  {\bf (b)} Figure 5 in \citep{tang2017real}. For real-time OPF problems in power system,  the average running time for iterative algorithms is 0.062 seconds with 12 iterations,
and the maximum is 0.376 seconds. {\bf (c)} Figure 6 in \citep{li2023real}. The path planner of autonomous vehicles needs to make decisions in 20 ms to avoid collision in the next 10 seconds.}
\label{fig_example}
\vspace{-0.3cm}
\end{figure}

Here, we highlight the relation between WMMSE and the finite horizon framework. To deploy the WMMSE framework in a real-time communication system, the WMMSE algorithm needs to find good solutions of beamformer $V_i$ and decoder $U_i$ within {\bf 4 iterations} (in Single-Input-Single-Output case) or {\bf 10 iterations} (in MIMO case)  \citep[Figure 1]{shi2011iteratively}. Otherwise, the users will suffer from high-latency communications.
Such requirement is ubiquitous in wireless communication systems. Another similar example is presented  in Figure \ref{fig_example} (a), where the iterative algorithms in the MIMO system are only allowed to operate {\bf within 14 iterations.}
This calls for the need to re-design the WMMSE algorithm under the finite horizon framework \eqref{eq_general_formulation} and boost performance within strict constraints of total iteration number $T$. Nevertheless, WMMSE is not the main focus of this script, and the relevant progress will be reported elsewhere.

\paragraph{Example 4: Optimal power flow (OPF) problem.}  OPF is a fundamental problem in power systems
to determine the most economical and secure way to distribute electrical power while satisfying network constraints \citep{chatzivasileiadis2018optimal}. DC-OPF problem is a  linear approximations of the actual nonlinear OPF problem
\footnote{The name of DC-OPF is a bit misleading, as this method does not assume the use of direct current (as opposed to alternating current). Rather, it is merely a linearization of AC-OPF problem, which is another name for the original  OPF problem \citep{mones2021gentle}.},   
and it is widely used for real-time operational planning and control due to its computational efficiency and ability to handle large-scale networks \citep{frank2012optimal}.
DC-OPF minimizes the total generation cost subject to the limits of generator operation, the power balance equation, and
the transmission line capacity constraints \citep{christie2000transmission}.  DC-OPF admits an LP formulation as follows \citep{pan2021deepopf}.

\begin{equation}
    \label{eq_example_opf}
    \begin{aligned}
        \min_{\{P_{G i}\}_{i =1}^{ N_{\text {gen}}},\{\theta_i\}_{i =1}^{ N_{\text {bus}}}} &\sum_{i=1}^{ N_{\text{gen }}} c_i P_{G i} \\
        \text{s.t. } & P_{G i}^{\min } \leq P_{G i} \leq P_{G i}^{\max }, i=1,2, \cdots, N_{\text {bus }}, \\
                     & \sum_{j = 1}^{N_{\text {bus}}}B_{i,j}\theta_j =  P_{G i} -  P_{D i},  i = 1, 2, \cdots, N_{\text {bus }}, \\
                     & \frac{1}{x_{i j}}\left(\theta_i-\theta_j\right) \leq P_{i j}^{\max }, i, j=1,2, \cdots, N_{\text {bus }}, 
    \end{aligned}
\end{equation}

where $N_{\text {gen }}$ denotes the number of generators and $N_{\text {bus }}$ is the number of buses.  $P_{G i}$ (and $P_{D i}$) denotes the power output (and consumption) of the generator in the $i$-th bus. $P_{G i}^{\min }$ and $P_{G i}^{\max }$ are the output limits of generators in the $i$-th bus, respectively. $\theta_i$ is the phase angle at the $i$-th bus.  $B_{i,j}$ and $x_{i,j}$ are constant terms from the power balance equations. More descriptions can be seen in \citep{pan2021deepopf}.

We highlight the relation between OPF problems and the finite horizon framework. In practice,  LP solvers must return their solutions to DC-OPF {\bf with far less than 1 second}; otherwise, it will be too slow to adjust the system's operating points in response to {\it real-time} changes in renewable power injection \citep{zhang2021convex,babaeinejadsarookolaee2019power,tang2017real}.   Unfortunately, This makes most standard solvers incompetent in these scenarios \citep{tang2017real}.
Currently, the average algorithm running time for real-time OPF is  {\bf 0.062 second with 12 iterations}  \citep{tang2017real} (see Figure \ref{fig_example} below).  \citet{tang2017real} also pointed out that  ``{\it most algorithms in the literature must wait until the iteration has converged because
the intermediate iterates typically do not satisfy power flow
equations and are not implementable.}" They also pointed out that  the traditional algorithms
will be ``{\it inadequate for future power grids}", and we need ``{\it real-time OPF algorithms}"
 that can respond quickly to network changes and maintain (sub)optimality."

To sum up, future power system requires algorithms to operate within a strict upper bound on the total iteration number $T$, and this raises the question of optimal algorithm design under framework \eqref{eq_general_formulation}.

\paragraph{Example 5: Path planning in autonomous vehicles.} Path planning is a critical component of autonomous vehicles, responsible for determining the optimal route of the vehicle within a certain time window. 
A qualified path planner must generate trajectories that allow the vehicle to reach destinations, avoid obstacles, and minimize fuel consumption, all while taking into account for environmental constraints such as buildings, pedestrians, and other moving vehicles. The path planner must make decisions every few seconds to ensure a safe and efficient path in real-time \citep{culligan2006online}.  Typically, the search for a good planner can be formulated as various optimization problems including  LP \citep{chasparis2005linear,chen2012three,kiessling2022feasible}, mix-integer LP (MILP) \citep{culligan2006online,toupet2006real}, convex QCQP \citep{subosits2019racetrack}, and others \citep{dolgov2008practical}.  For example, the following LP often arises in time-optimal motion planning \citep{kiessling2022feasible}.

\begin{equation}
\label{eq_example_auto_vehicles}
\begin{aligned} \min _{w \in \mathbb{R}^{N_w} } & c^{\top} w \\ 
\text { s.t. } & C w+G ^{\top} P_y(w-\hat{w})=0, \\ 
& A w+b \leq 0, \\ 
& \left\|P_y(w-\hat{w})\right\|_{\infty} \leq \Delta,
\end{aligned}
\end{equation}

where the optimization variable $w$ is the concatenation of state, control, and slack variables. More introduction for $w$ and the physical meanings of constant terms including $c, C, G, P_y, \hat{w}, A, b, \Delta$ can be seen in \citep{kiessling2022feasible}. Here, we highlight the relation between the autonomous path planner and the finite horizon framework. To generate a safe and efficient path in real-time, the sub-problems of LP \eqref{eq_example_auto_vehicles} must be solved in {\bf within 10 ms} \citep[Table VI]{chen2012three}. Similarly, the  MILP formulation in \citep{culligan2006online} and \citep[Table 1]{yu2024rigid} need to be solved in {\bf within 100 ms} to get the plan for the next 10 seconds; the convex QCQP in \citep{subosits2019racetrack} needs to be solved in {\bf within 20 ms} for the plan in the next 10 seconds (Figure \ref{fig_example} (c)); the  MIQP in \citep[Table II]{qian2016optimal} needs to be solved in {\bf within 228 ms}. All these strict constraints on operation time post an upper bound on the total iteration number $T$.

\paragraph{Goal.} Motivated by the examples above, we aim to design new stepsize rules for LP algorithms under the finite horizon framework \eqref{eq_general_formulation}. We study LP because LP (and its variants) frequently arise in engineering scenarios with strict iteration budget constraints for solvers, such as {\bf Example 4} and {\bf 5} above.  As an initial attempt, we will work on a most simple yet fundamental algorithm: the primal-dual method. We will show that the primal-dual method can be accelerated by about $3.9\times$ if the stepsize rules are re-designed using formulation \eqref{eq_general_formulation}.  The acceleration comes at a price of sacrificing asymptotic guarantees, but  as argued above, this trade-off is inconsequential for the a wide range of applications.

It would be intriguing to apply framework \eqref{eq_general_formulation} to more advanced methods such as the  momentum or the preconditioned variants of primal-dual method \citep{applegate2021practical}, where we believe substantial potential for improvement remains. Given the technical complexity involved in the analysis of the primal-dual method, we will concentrate on this algorithm for now and defer the extensions to future research.

\paragraph{Challenges.} Our stepsize design of the primal-dual method is inspired by the seminal work \citep{young1953richardson} in {\bf Example 1}. \citet{young1953richardson} focused on GD for quadratic functions, and we provide a (rather non-trivial) extension of their idea to the primal-dual method. There are at least two-fold challenges. 

\begin{itemize}[topsep=1pt,parsep=1pt,partopsep=1pt, leftmargin=*]
    \item {\bf First: non-symmetric update matrix.} The primal-dual method is driven by a  {\it non-symmetric} update matrix (i.e., $M$ is non-symmetric in \eqref{eq_lp_no_inequality}), whereas  GD for quadratic is driven by a symmetric update matrix (i.e., Hessian).  We highlight that for non-symmetric matrices,  their singular values are usually largely different from their eigenvalues.  This posts the following difficulties. First, most existing stepsize designs of GD including \citep{young1953richardson} focus on controlling Hessian eigenvalues. Since we aim to control singular values, these eigenvalue-based approaches cannot be directly applied.  Second, in general, there are very few tools to analyze the singular values of non-symmetric matrices.  Similar difficulties are also reported in  \citep{kittaneh2006spectral,sun2020efficiency,sun2021worst}.

\item {\bf Second: non-asymptotic analysis.} Since we focus on the {\it non-asymptotic} regime with finite $T$, we cannot use any asymptotic properties of the singular value spectrum. For instance, we cannot apply the asymptotic relation between singular values and eigenvalues \citep{saad2003iterative} or any similar theoretical results. To our knowledge, there are very few helpful tools in this non-asymptotic regime.
\end{itemize}

We present our stepsize rule in {\bf Algorithm \ref{algo_stepsize_rule}}.  We address the above two challenges by uncovering the hidden convexity inherent in \eqref{eq_general_formulation} for any $T$, a property uniquely brought up by the primal-dual method. We then propose a convex SDP reformulation of \eqref{eq_general_formulation}  and numerically search for the desired stepsize. The SDP only involves matrices with size $4\times 4$ and can be solved efficiently by modern solvers. While having an (exact) convex reformulation may not be as convenient as having a closed-form solution as in \eqref{eq_young_stepsize} 
 \citep{young1953richardson}, it is still better news compared to the (exact) nonconvex reformulation for GD stepsize \citep{das2024branch}. We introduce the detailed derivation of {\bf Algorithm \ref{algo_stepsize_rule}} in the next section.

 \begin{algorithm}[H]
  \caption{Finite Horizon stepsize rule}
  \label{algo_stepsize_rule}
  \begin{algorithmic}[1]
  \STATE{Input the total iteration budget $T$.}
  \STATE{Input the smallest \& largest singular value of $A$ (the constraint matrix in LP), denoted as $\mu, L$.}
  \STATE{Get $n_{\text{sample}}$ evenly spaced numbers over $[\mu,L]$, denoted as $\{\sigma_i\}_{i = 1}^{n_{\text{sample}}}$.}
  \STATE{Formulate the SDP \eqref{eq_sdp} using $\{\sigma_i\}_{i = 1}^{n_{\text{sample}}}$ and $T$. Solve it using the Interior Point Methods, return the optimal solution $\{a_t\}_{t= 1}^T$ }
  \STATE{Find the roots of the polynomial $p(x) = 1+ a_1x + a_2 x^2 \cdots + a_T x^T$ using the algorithm in  \citep{horn2012matrix}, denoted as $\{r_t\}_{t=1}^T$. }
  \STATE{Calculate $\eta_t = 1 / r_t$, for $t = [T]$.}
  \STATE{Return $\{\eta_t\}_{t=1}^T$ as the stepsize rule. }
  \end{algorithmic}
  \end{algorithm}

\section{Main Methods}
\label{sec_main_methods}

\subsection{LP Formulation and Algorithms}
\label{sec_formulation}

We consider the standard form of LP as follows.

\begin{equation}
\label{eq_lp}
\begin{aligned}
\min _{x \in \mathbb{R}^n} & \quad c^{\top} x \\
\text { s.t. } & A x=b \\
& x \geq 0,
\end{aligned}
\end{equation}
where $c \in \mathbb{R}^n,  b \in \mathbb{R}^m$ and we assume $A \in \mathbb{R}^{m\times n}$ is  full rank  and $m \leq n$. The singular values of $A$ satisfy   $\sigma(A)\subseteq [\mu,L]$, where $\mu, L >0$. We consider the primal-dual form of \eqref{eq_lp}: 
\begin{equation}
\label{eq_primal_dual_LP}
\min_{x\geq 0} \max_{y \in \mathbb{R}^m} \mathcal{L}_{\beta}(x, y):= c^{\top}x + y^{\top}(A x-b)+\frac{\beta}{2}\|A x-b\|_2^2,
\end{equation}
where $\mathcal{L}_{\beta}(x, y)$ is the augmented Lagrangian function. To solve the primal-dual formulation of LP \eqref{eq_primal_dual_LP}, a most simple first-order method (FOM) is  the primal-dual method, or equivalently,  the projected gradient descent-ascent method (GDA):

\begin{equation}
\label{eq_primal_dual}
\left\{\begin{array}{l}
x^{t+1} =\left(x^t - \eta_t \left(c + A^{\top} y^t  + \beta A^{\top}(Ax^t -b) \right) \right)^+   \\
y^{t+1} = y^t  + \eta_t \left(A x^t - b\right),
\end{array}\right.
\end{equation}

where $(\cdot)^+$ operator is the projection onto the non-negative orthant and is cheap to implement. We denote  $\eta_t$ as the stepsize at $t$-th step for the primal and dual update.
In this work, we will use the same stepsize $\eta_t$ for both primal and dual updates. We ask the following question: 

\begin{snugshade}
\begin{center}
Given a finite iteration budget $T$,  how to find the optimal stepsize rule $\{\eta_t\}_{t=1}^T$ of \eqref{eq_primal_dual} for these $T$ iterations? 
\end{center}
\end{snugshade}

\subsection{Finite Horizon Stepsize Rule}
We now design Finite Horizon stepsize rule for the primal-dual method \eqref{eq_primal_dual}. We first rewrite the primal update as follows.

\begin{eqnarray}
    x^{t+1} - x^* &=& \left(x^t - \eta_t \left(c + A^{\top} y^t  + \beta A^{\top}(Ax^t -b) \right) \right)^+  - x^* \nonumber \\ 
     &\overset{(i)}{=}& \left((x^t - x^*) - \eta_t \left(c + A^{\top} y^t  + \beta A^{\top}A(x^t -x^*) \right) \right)^+,  \nonumber \\
    &\overset{(ii)}{=}&\left((x^t - x^*) - \eta_t \left(  A^{\top}( y^t - y^*)  + \beta A^{\top}A(x^t -x^*) \right) \right)^+,   \nonumber\\
    &\overset{}{=}& \left( (I - \eta_t \beta A^{\top}A)(x^t - x^*) - \eta_t \left(  A^{\top}( y^t - y^*)   \right) \right)^+,   \label{eq_primal_update}
\end{eqnarray}

where $x^*$ and $y^*$ are optimal primal and dual solutions, respectively; $(i)$ is due to: $\left(x^*\right)^+ = x^*, Ax^* = b$. $(ii)$ is due to the KKT condition $A^T y^* = -c$. Similarly, we rewrite the dual update as follows.

\begin{equation}
\label{eq_dual_update}
    y^{t+1} - y^* = y^t  -y^*  + \eta_t A(x-x^*). 
\end{equation}
 
Define $z = [(x)^\top, (y)^\top  ]^\top \in \mathbb{R}^{(n +m)}$, the update equations \eqref{eq_primal_update} and \eqref{eq_dual_update} can be rewritten as

\begin{equation}
\label{eq_primal_dual_matrix}
 z^{t+1} - z^*= P_t \left[\begin{array}{cc}
(I_{n\times n} - \eta_t \beta A^\top A) &   - \eta_t A^\top\\
\eta_t A & I_{m\times m} \\
\end{array}\right] \left(z^t - z^*\right),
\end{equation}

where $P_t$ is a diagonal matrix defined as follows:

\begin{equation}
    P_t = \left[\begin{array}{cccc}
p_{t,1} & \cdots & 0 & 0   \\
  & \ddots &   &    \\
 0 & \cdots & p_{t,n} & 0   \\
0& \cdots & 0 & I_{m \times m}
\end{array}\right] \in \mathbb{R}^{(n+m) \times (n+m)}, p_{t,i}  = \left\{\begin{array}{cc}
0 & \text{if $x_{t,i} <0$ } \\
1 & \text{otherwise.} 
\end{array}\right .
\end{equation}

Define 

\begin{equation}
\label{eq_M}
M = \left[\begin{array}{cc}
 \beta A^\top A & A^\top \\
-A & 0_{m\times m}
\end{array}\right],
\end{equation}

then the update equation \eqref{eq_primal_dual_matrix} can be further rewritten as 

\begin{eqnarray}
    z^{t+1} - z^* &=& P_t(I - \eta_t M)(z^t - z^*)\\
    &=& P_t(I - \eta_tM) P_{t-1}(I - \eta_{t-1}M) \cdots P_{1}(I - \eta_{1}M)(z^0-z^*).
\end{eqnarray}

 Given a class of LP instances with $\sigma(A)\subseteq [\mu,L]$, we aim to find the optimal stepsize rule  $\{\eta_t\}_{t=1}^T$ for  $T$ total iterations, where $T>0$ is a finite integer. In particular, we aim to minimize the "worst-case" distance to the optimal solution at iteration $T$. Here, the "worst-case" is in the sense of LP instance $A$ and initial point $z^0$. 
This gives rise to the following optimization problem.

 \begin{equation}
 \label{eq_minmax_with_projection}
      \min_{\BLUE{\eta_1, \cdots, \eta_T}}\max_{A, s.t., \sigma(A)\subseteq [\mu,L], z \not\in \mathcal{Z}^*}\frac{\|P_T(I - \BLUE{\eta_T} M) P_{T-1}(I - \BLUE{\eta_{T-1}}M) \cdots P_{1}(I - \BLUE{\eta_{1}}M) (z - z^*) \|_{2}}{\|z - z^*\|_2} ,
 \end{equation}

 where $\mathcal{Z}^* = \{z; z \text{ is the optimal primal-dual solution of \eqref{eq_lp}}, z\in \operatorname{null}(M)\}$, and $z^*$ is a closest point  in $\mathcal{Z}^*$ to $z$.
To solve \eqref{eq_minmax_with_projection}, one immediate obstacle lies in how to handle the projection matrices $P_t$. There are at least two challenges. 
{\bf First}, the realization of $P_t$ depends on the specific algorithmic trajectory, which is almost unpredictable in the worst-case-based formulation \eqref{eq_minmax_with_projection}. {\bf Second},  $P_t$ does not commute with $M$, i.e., $P_t (I - \eta_t M)$ does not share the same singular vectors with $(I - \eta_t M)$. In the worst case, the singular vectors of the update system changes at every step, making it hard to analyze.

  In this work, we view the projection $P_t$ as random projection matrices, i.e., $p_{t,i} \overset{\iid}{\sim} \text{Bernoulli}(1- \pproj)$, and we consider the expected output of the random projection at iteration $T$:

\begin{eqnarray}
    \Ex \left[z^T - z^*\right]&=&\mathbb{E}_{P_1, \cdots, P_T}\left[  P_T(I - \eta_TM) {P_{T-1}}(I - \eta_{T-1}M) \cdots {P_{1}}(I - \eta_{1}M)(z^0-z^*)  \right] \\
     & = & P(I - \eta_TM) P(I - \eta_{T-1}M) \cdots P(I - \eta_{1}M)(z^0-z^*),
\end{eqnarray}

where

\begin{equation}
\label{eq_random_projection}
    P = \left[\begin{array}{cccc}
1- \pproj & \cdots & 0 & 0   \\
  & \ddots &  &    \\
 0 & \cdots & 1- \pproj& 0   \\
0& \cdots & 0 & I_{m \times m}
\end{array}\right] \in \mathbb{R}^{(n+m) \times (n+m)}.
\end{equation}

We now rewrite the expected output. We first  define   $\Tilde{I} \in \mathbb{R}^{(n+m)\times (n+m)}$ and $\Tilde{M} \in \mathbb{R}^{(n+m)\times (n+m)}$ as follows.

\begin{equation}
\label{eq_I_tilde}
    \Tilde{I} = \left[\begin{array}{cccc}
I_{n\times n} & 0_{n\times m} \\
0_{m\times n} & 0_{m\times m}
\end{array}\right] \in \mathbb{R}^{(n+m) \times (n+m)}, \Tilde{M} = \Tilde{I}  M.
\end{equation}

Then we have

{\small
\begin{eqnarray}
&&\mathbb{E}_{P_1, \cdots, P_T}\left[  P_T(I - \eta_TM)  \cdots P_{1}(I - \eta_{1}M) \right] \nonumber \\
   & = &   P(I - \eta_TM) \cdots P(I - \eta_{1}M) \label{eq_expected_product}
   \\ 
    & =& \left(I - \eta_T M - \pproj (\Tilde{I} - \eta_T \Tilde{M}) \right) \cdots \left(I - \eta_1  M - \pproj (\Tilde{I} - \eta_1 \Tilde{M}) \right)  \nonumber\\
    & = & (I - \eta_TM)\cdots (I - \eta_{1}M) 
    \nonumber\\
    && + \pproj \sum_{t=1}^T (I - \eta_T M ) \cdots (I - \eta_{t+1} M )(\Tilde{I} - \eta_t \Tilde{M})(I - \eta_{t-1} M )\cdots (I - \eta_{1} M ) \nonumber \\
     &&  + \pproj^2 \sum_{1\leq i,j \leq T}  (I - \eta_T M ) \cdots (I - \eta_{j+1} M )(\Tilde{I} - \eta_j \Tilde{M})(I - \eta_{j-1} M )\cdots (I - \eta_{i+1} M )(\Tilde{I} - \eta_i \Tilde{M})(I - \eta_{i-1} M ) \cdots (I - \eta_{1} M ) 
 \nonumber\\
  && + \cdots \nonumber \\
  && + \pproj^T (\Tilde{I} - \eta_T \Tilde{M}) \cdots  (\Tilde{I} - \eta_1 \Tilde{M}) \nonumber \\
    & := & \underbrace{(I - \eta_TM)\cdots (I - \eta_{1}M)}_{(a)} + \underbrace{  \pproj \operatorname{Residue} }_{(b)}.
    \label{eq_expected_product_3}
\end{eqnarray}
}

In the sequel, we focus on optimizing the leading term $(a)$. We focus on $(a)$ because: we numerically find that $\pproj$ is small when initializing $z^0$ in the positive orthant, and the resulting $\pproj$ is usually $\leq 0.1$.
Further, when  $\pproj$  is small, we find that the operator norm of \eqref{eq_expected_product} is close to the operator norm of $(a)$. As such, term $(a)$  usually plays a more crucial role in the overall convergence rate.  We relegate the numerical evidence in Section \ref{sec_experiments}. 
Now we focus on term $(a)$ and solve the following problem \eqref{eq_minmax_no_projection}.

\begin{equation}
\label{eq_minmax_no_projection}
      \min_{\BLUE{\eta_1, \cdots, \eta_T}}\max_{A, s.t., \sigma(A)\subseteq [\mu,L], z \not\in \mathcal{Z}^*}\frac{\|(I - \BLUE{\eta_T} M) (I - \BLUE{\eta_{T-1}}M) \cdots (I - \BLUE{\eta_{1}}M) (z - z^*) \|_{2}}{\|z - z^*\|_2}.
\end{equation}

\paragraph{Challenges for solving \eqref{eq_minmax_no_projection}.} We note that \eqref{eq_minmax_no_projection} is non-trivial to solve. There are at least three challenges. 
\begin{itemize}[topsep=1pt,parsep=1pt,partopsep=1pt, leftmargin=*]
    \item {\bf C1: Nonconvexity.} The optimization variable $\BLUE{\eta_1, \cdots, \eta_T}$ appears in \eqref{eq_minmax_no_projection} as a product form. This makes the objective function nonconvex. 
    \item {\bf C2: Non-symmetric $M$.} Eq. \eqref{eq_minmax_no_projection} involves controlling the  singular values of $(I - \eta_T M)\cdots (I - \eta_1 M)$, which is the product of non-symmetric matrices. It is worth mentioning that: for non-symmetric matrices,  their singular values are usually largely different from their eigenvalues, and there is no clear relation between these two.  
    This posts the following difficulties. First, most existing stepsize design of GD (e.g., \citep{young1953richardson}) focus on controlling eigenvalues of symmetric matrices (Hessian). Since we aim to control singular values, which are {\it not} equal to eigenvalues anymore in the non-symmetric case, these eigenvalue-based approaches cannot be directly applied.  Second, in general, there are very few tools to analyze non-symmetric matrices (actually, both singular values and eigenvalues of non-symmetric matrices are difficult topics with very limited results).  Similar difficulties in analyzing 
    non-symmetric matrices  are also reported in  \citep{kittaneh2006spectral,sun2020efficiency,sun2021worst}.

    Here, we post a simple example to illustrate that: for non-symmetric matrices, the gap between singular values and eigenvalues can be arbitrarily large. Consider $M=\left(\begin{array}{ll}0 & \alpha \\ 0 & 0\end{array}\right)$, its eigenvalues are 0 with multiplicity 2, while its maximal singular values equals $|\alpha|$, which could be unbounded as $\alpha \rightarrow \infty$. This huge gap makes all eigenvalue-based methods inapplicable.
  
\item {\bf C3: Non-asymptotic analysis.} One possible solution to {\bf C2} is to transform singular values into eigenvalues and then resort to the eigenvalue analysis in \citep{young1953richardson}. Unfortunately,  for non-symmetric matrices, most existing relations between singular values and eigenvalues only apply in the asymptotic sense.  For instance, \citep[Theorem 1.12]{saad2003iterative} states that  $\lim _{T \rightarrow \infty}\sigma_1(M^T)^{1 / T}=|\lambda_1(M)|$. In indeed, we have $\lim _{T \rightarrow \infty}\sigma_1(M^T)^{1 / T} = 0$ for the $2\times 2$ example above.  However,  in the non-asymptotic regime, the gap between singular values and eigenvalues can still be arbitrarily large. 
In the non-asymptotic regime, one potential approach is to use the symmetrization trick, i.e., to analyze the eigenvalues of $M^\top M$. However, this is still not easy since the eigenbasis of $M$ is {\it not} orthonormal, and thus it is unclear whether we can diagonalize $M^\top M$ to extract its eigenvalues.
\end{itemize}

Perhaps a bit surprisingly, we find \eqref{eq_minmax_no_projection} admits hidden convexity for any $T$ under proper changes of variables. We will show that \eqref{eq_minmax_no_projection} has an exact convex SDP reformulation and can be solved efficiently in polynomial time. 

\paragraph{Proposed solution.} We will use the symmetrization trick to analyze the eigenvalues of $M^\top M$, but in a rather non-conventional way. 
Our derivation 
contains three steps.  {\bf Step 1:} Instead of diagonalizing $M^\top M$, we do the next-best thing: we block diagonalize $M$ into multiple $2\times 2$ blocks using orthonormal matrix. {\bf Step 2:} we rewrite the objective into a polynomial of these $2\times 2$ matrix blocks. {\bf Step 3:} we control the eigenvalues in each block using linear matrix inequalities.  

{\bf Step 1.} Based on the singular value decomposition (SVD) of $A = U \Sigma V^{\top} = \sum_{i = 1}^m \sigma_i u_i v_i^{\top}$, we rewrite $M$ as follows.

\begin{equation}
    M = \left[\begin{array}{cc}
 \beta V \Sigma^\top \Sigma V^{\top} & V \Sigma^\top  U^\top  \\
- U\Sigma V^{\top}  &  0_{m \times m}  
\end{array}\right]_{(m+n)\times (m+n) } = \left[\begin{array}{cc}
 \sum_{i = 1}^m \beta\sigma_i^2  v_i v_i^\top &  \sum_{i = 1}^m \sigma_i  v_i u_i^\top \\
-  \sum_{i = 1}^m \sigma_i  u_i v_i^\top &  0_{m \times m}  
\end{array}\right]_{(m+n)\times (m+n) }.
\end{equation}

Notice that for any $i = [m]$, we have

\begin{eqnarray}
\left[\begin{array}{cc}
\beta \sigma_i^2  v_i v_i^\top &  \sigma_i  v_i u_i^\top \\
-  \sigma_i  u_i v_i^\top &  0_{m \times m}  
\end{array}\right]_{(m+n)\times (m+n) }  &=& \left[\begin{array}{cc}
  v_i  &  0_{n\times 1} \\
0_{m\times 1} &  u_i 
\end{array}\right]_{(m+n)\times 2} \left[\begin{array}{cc}
 \beta \sigma_i^2 &  \sigma_i\\
 - \sigma_i &  0
\end{array}\right]_{2\times 2}  \left[\begin{array}{cc}
  v_i^\top  &  0_{1\times m } \\
0_{1\times n} &  u_i^\top 
\end{array}\right]_{2\times (m+n)} 
\end{eqnarray}

Define 
\begin{equation}
  Q_i =   \left[\begin{array}{cc}
  v_i  &  0_{n\times 1} \\
0_{m\times 1} &  u_i 
\end{array}\right] \in \mathbb{R}^{(m+n)\times 2} , \quad
B(\sigma_i) = \left[\begin{array}{cc}
  \beta \sigma_i^2 &  \sigma_i\\
 - \sigma_i &  0
\end{array}\right]\in \mathbb{R}^{2\times 2}.  
\end{equation}

We now block-diagonalize  $M$ using an orthonormal matrix.

{\small
\begin{eqnarray}
\label{eq_block_diagonalize_M}
    M &=& \left[\begin{array}{cc}
 \sum_{i = 1}^m \sigma_i^2 \beta v_i v_i^\top &  \sum_{i = 1}^m \sigma_i  v_i u_i^\top \\
-  \sum_{i = 1}^m \sigma_i  u_i v_i^\top &  0_{m \times m}  
\end{array}\right] = \sum_{i= 1}^m \left[\begin{array}{cc}
  v_i  &  0_{n\times 1} \\
0_{m\times 1} &  u_i 
\end{array}\right] \left[\begin{array}{cc}
 \beta \sigma_i^2 &  \sigma_i\\
 - \sigma_i &  0
\end{array}\right]\left[\begin{array}{cc}
  v_i^\top  &  0_{1\times m } \\
0_{1\times n} &  u_i^\top 
\end{array}\right]   \nonumber \\
&=&  
\left[\begin{array}{ccc}
Q_1,   & \cdots, & Q_m 
\end{array}\right]_{(n+m)\times 2m}
\left[\begin{array}{ccc}
  B(\sigma_1) &   \cdots & 0 \\
   & \ddots   &    \\
 0 &  \cdots  & B(\sigma_m) \\
\end{array}\right]_{2m \times 2m}
\left[\begin{array}{c}
Q_1^\top \\
\cdots \\
Q_m^\top \\
\end{array}\right]_{2m\times (n+m)}. \nonumber\\
&=& \!\!\!\!\!\!\left[\begin{array}{ccccc}
Q_1,  & \cdots, & Q_m, & \cdots, & Q_{m+n} 
\end{array}\right]_{(n+m)\times (n+m)}
\left[\begin{array}{cccc}
  B(\sigma_1) &   \cdots & 0  & 0\\
   &  \ddots  &  &  \\
 0 &  \cdots & B(\sigma_m) & 0 \\
  0 &  \cdots & 0 & 0_{(n-m)\times (n-m)} \\
\end{array}\right]_{(n+m) \times (n+m)}
\left[\begin{array}{c}
Q_1^\top \\
\cdots \\
Q_m^\top \\
\cdots \\
Q_{m+n}^\top \\
\end{array}\right]_{(n+m)\times (n+m)} \nonumber \\
&:=& Q \Lambda Q^\top,\nonumber
\end{eqnarray}
}
where $Q_{m+1},\cdots, Q_{m+n}$ form the orthonormal basis of $\operatorname{null}(M)$ and they make $Q$  an orthonormal matrix. This property will be important for the analysis in the sequel.

{\bf Step 2.}  Now we rewrite \eqref{eq_minmax_no_projection} based on the block diagonalization of $M$.

\begin{eqnarray}
\label{eq_minmax_no_projection_block_diagonal}
   \eqref{eq_minmax_no_projection}& = &       \min_{\BLUE{\eta_1, \cdots, \eta_T}}\max_{A, s.t., \sigma(A)\subseteq [\mu,L], z \not\in \mathcal{Z}^*}\frac{\|(I - \BLUE{\eta_T} Q \Lambda Q^\top) \cdots (I - \BLUE{\eta_{1}}Q \Lambda Q^\top) (z - z^*) \|_{2}}{\|z - z^*\|_2} \nonumber \\
   & \overset{}{ = }&  \min_{\BLUE{\eta_1, \cdots, \eta_T}}\max_{A, s.t., \sigma(A)\subseteq [\mu,L], z \not\in \mathcal{Z}^*}\frac{\|(I - \BLUE{\eta_T}  \Lambda) \cdots (I - \BLUE{\eta_{1}} \Lambda) Q^\top(z - z^*) \|_{2}}{\|Q^\top (z - z^*)\|_2}  \nonumber\\
    & \overset{(a)}{ = }&\min_{\BLUE{\eta_1, \cdots, \eta_T}}\max_{\sigma\subseteq [\mu,L]}\left\|\left(I_{2\times 2} - \BLUE{\eta_T}B(\sigma) \right) \cdots \left(I_{2\times 2}  - \BLUE{\eta_{1}}B (\sigma)\right) \right\|_{\operatorname{op}} \nonumber\\
    & \overset{(b)}{ = }&\min_{\BLUE{a_1, \cdots, a_T}}\max_{\sigma\subseteq [\mu,L]}\left\| I_{2\times 2}  + \BLUE{a_1} B(\sigma) + \BLUE{a_2} B(\sigma)^2 \cdots +\BLUE{a_T} B(\sigma)^T \right\|_{\operatorname{op}}   \nonumber
\end{eqnarray}

where $(a)$ is because we are maximizing over $z\not\in \mathcal{\mathcal{Z}^*}$ so we only need to consider the non-zero eigenvlaues of $M$, which all reside in $B(\sigma)$. In 
$(b)$, we  rearrange the product term  as a $T$-th order polynomial of $B$, and we optimize $(\eta_1, \cdots \eta_T)$ is equivalent to optimizing the coefficient of the polynomial $(a_1, \cdots, a_T )$.

{\bf Step 3.}  Note that \eqref{eq_minmax_no_projection_block_diagonal} is equivalent to the formulation below.

\begin{equation}
 \label{eq_minmax_polynomial}
 \begin{array}{rl}
\min _{a_1, \ldots, a_T, s} & s \\
\text { s.t. } & \left( I_{2\times 2}  + \BLUE{a_1} B(\sigma) +\cdots +\BLUE{a_T} B(\sigma)^T \right)^\top  \left( I_{2\times 2}  + \BLUE{a_1} B(\sigma) + \cdots  +\BLUE{a_T} B(\sigma)^T \right) \preceq s^2 I_{2\times 2}, \nonumber\\
& \sigma\subseteq [\mu,L]. \\
\end{array} 
\end{equation}

Note that \eqref{eq_minmax_polynomial} can be rewritten as follows using linear matrix inequalities:

\begin{equation}
 \begin{array}{rl}
\min _{\BLUE{a_1, \ldots, a_T}, s} & s \\
\text { s.t. } & 
\left[\begin{array}{cc}
s I_{2\times 2} &\left( I_{2\times 2}  + \BLUE{a_1} B(\sigma) + \BLUE{a_2} B(\sigma)^2 \cdots \BLUE{a_T} B(\sigma)^T \right)^\top \\
I_{2\times 2}  + \BLUE{a_1} B(\sigma)  + \BLUE{a_2} B(\sigma)^2 \cdots \BLUE{a_T} B(\sigma)^T & s I_{2\times 2}  \\
\end{array}\right] \succeq 0, \\
& \sigma\subseteq [\mu,L].
\end{array}
\end{equation}

We further re-arrange it as follows:

\begin{equation}
\label{eq_sdp_inf_constraint}
 \begin{array}{rl}
\min _{\BLUE{a_1, \ldots, a_T}, s} & s \\
\text { s.t. } & s I_{4\times 4} + 
\left[\begin{array}{cc}
0 & I_{2\times 2} \\
I_{2\times 2} & 0 \\
\end{array}\right]_{4\times 4}  + \BLUE{a_1}  \left[\begin{array}{cc}
0 & B(\sigma)^\top \\
B(\sigma ) & 0 \\
\end{array}\right]_{4\times 4}  + 
\cdots  + \BLUE{a_T} \left[\begin{array}{cc}
0 & ( B(\sigma)^T)^\top \\
B(\sigma)^T & 0 \\
\end{array}\right]_{4\times 4} \succeq 0, \\
& \sigma\subseteq [\mu,L].
\end{array}
\end{equation}

Note that \eqref{eq_sdp_inf_constraint} is an SDP, but admits an infinite number of matrix-inequality constraints. For the ease of computation, we further discretize the interval $\sigma\subseteq [\mu,L]$ by
 returning $n_{\text{sample}}$ evenly spaced $\sigma_i$ over $[\mu,L]$, for $i= [n_{\text{sample}}]$.

\begin{equation}
\label{eq_sdp}
 \begin{array}{rl}
\min _{\BLUE{a_1, \ldots, a_T}, s} & s \\
\text { s.t. } & s I_{4\times 4} + 
\left[\begin{array}{cc}
0 & I_{2\times 2} \\
I_{2\times 2} & 0 \\
\end{array}\right]_{4\times 4}  + \BLUE{a_1}  \left[\begin{array}{cc}
0 & B(\sigma_i)^\top \\
B(\sigma_i) & 0 \\
\end{array}\right]_{4\times 4}  + 
\cdots  + \BLUE{a_T} \left[\begin{array}{cc}
0 & (B(\sigma_i)^T)^\top \\
B(\sigma_i)^T & 0 \\
\end{array}\right]_{4\times 4} \succeq 0, \\
& i = [n_{\text{sample}}].
\end{array}
\end{equation}

Now Problem \eqref{eq_sdp} has a finite number of constraints and only involves small-scaled matrices with size $4 \times 4$. This problem can be efficiently solved in polynomial time using the Interior Point Methods (IPMs).

After solving \eqref{eq_sdp}, we factorize the polynomial $p(x) = 1+ a_1x + a_2 x^2 \cdots + a_T x^T = (1- \eta_1 x) (1-\eta_2x) \cdots (1-\eta_Tx)$ and use $\{\eta_t\}_{t= 1}^T$ as the final stepsize rule. Note that the polynomial factorization process can be efficiently done via two steps: First,  find the roots of $p(x)$, denoted as $\{r_t\}_{t= 1}^T$. This can be efficiently done by computing the eigenvalues of the companion matrix \citep{horn2012matrix}. Second, take the inverse over these roots and get $\eta_t = 1 / r_t$, for $t = [T]$.

The complete algorithm for finding $\{\eta_t\}_{t= 1}^T$ is presented in {\bf Algorithm \ref{algo_stepsize_rule}}, and we call it Finite Horizon stepsize rule (for the primal-dual method).

\section{Theoretical Guarantee for Finite Horizon Stepsize Rule}

We now prove that Finite Horizon stepsize rule can reach accelerated convergence over the constant stepsize. Consider the case where is no non-negativity constraint $x \geq 0$ in \eqref{eq_lp}, that is

\begin{equation}
    \label{eq_lp_no_inequality}
\begin{aligned}
\min _{x \in \mathbb{R}^n} & \quad c^{\top} x \\
\text { s.t. } & A x=b.
\end{aligned}
\end{equation}

To ensure the boundedness of this  LP, we will assume $c \in \operatorname{range}(A^\top)$.
In this case,  the objective function will be constant within the feasible set, i.e., $c^\top x = c^\top A^\top (A A^\top)^{-1}b$ for all feasible $x$, and the problem is equivalent to solving the linear system $A x = b$. In this case, the primal-dual method admits the following update rule:
\begin{equation}
    \label{eq_primal_dual_no_projection}
z^{t+1} - z^*  = (I - \eta_t M) (z^t - z^*),
\end{equation}
where all the variables follow the same definition in Section \ref{sec_formulation}. 
We now prove that:  At the $T$-th (pre-determined) iteration,  the complexity of our stepsize rule is at most  $\mathcal{O}(\sqrt{\kappa}\log(\frac{1}{\epsilon}))$, which achieves acceleration over the optimal constant stepsize $\Omega(\kappa\log(\frac{1}{\epsilon}))$, where $\kappa$ is the condition number of the update matrix. We first recall the classical lower bound for the asymptotic constant step size.

\begin{prop} ({\bf Lower bound for constant stepsize})
\label{prop_constant_lower_bound}
    Consider problem \eqref{eq_lp_no_inequality} with $A \in \mathbb{R}^{m\times n}$ full row rank and $m \leq n$. Assume the  singular values of $A$ lie in $[\mu, L]$, where $\mu, L >0$. Suppose $b \in \operatorname{range}(A)$ and  $c \in \operatorname{range}(A^\top)$ so that \eqref{eq_lp_no_inequality} is feasible and bounded.  Consider the primal-dual method \eqref{eq_primal_dual_no_projection} with constant stepsize. 
    Denote $z^\top = (x^\top, y^\top) \in \mathbb{R}^{n+m}$ as the concatenation of the primal variable $x \in  \mathbb{R}^{n}$ and dual variable $y \in  \mathbb{R}^{m}$. Assume the augmented Lagrangian coefficient $\beta$ is large such that $\beta \geq \frac{2}{\mu}$. Define the  asymptotic optimal constant stepsize as:

    \begin{equation}
    \label{eq_def_optimal_constant_stepsize}
    \eta^*  = \arg min_\eta\max_{A, s.t., \sigma(A)\subseteq [\mu,L], z \not\in \mathcal{Z}^*} \limsup_{T\rightarrow \infty} \frac{\| (I -\eta M)^T(z^0 - z^*)\|_2}{\|z^0 - z^*\|_2}.
    \end{equation}

    Then  we have (i): $\eta^* = (\beta L^2 + \sqrt{\beta^2 L^4 - 4 L^2}  + \beta \mu^2  -  \sqrt{ \beta^2 \mu^4 - 4 \mu^2} )^{-1}$. (ii): For any $A, s.t., \sigma(A)\subseteq [\mu,L]$, there exists a real-valued initialization $z^0$ such that the following lower bound hold  at any $T>0$:
    
    \begin{equation}
    \label{eq_constant_step_rate}
        \text{dist}(z^T, \mathcal{Z}^*) \geq \left(1 -  \frac{2}{\frac{\beta L^2 + \sqrt{\beta^2 L^4 - 4 L^2}}{\beta \mu^2 -  \sqrt{ \beta^2 \mu^4 - 4 \mu^2}} + 1}\right)^T \text{dist}(z^0, \mathcal{Z}^*),
     \end{equation}
    where  $\text{dist}(z, \mathcal{Z}^*)$ denotes the Euclidean distance of $z$ to the set of optimal primal-dual solutions. Eq. \eqref{eq_constant_step_rate} can also be rewritten as follows: 
\begin{equation}
\label{eq_constant_step_rate_kappa}
        \text{dist}(z^T, \mathcal{Z}^*) \geq \left(1  - \frac{2}{\kappa +1}\right)^T \text{dist}(z^0, \mathcal{Z}^*),
    \end{equation}
    where $\kappa$ as the condition number of matrix $M$ (as defined in \eqref{eq_M}).
\end{prop}

The proof can be seen in Section \ref{sec_proof_prop_constant_lower_bound}.  This result provides a complexity lower bound $\Omega\left(\kappa \log(1/\epsilon) \right)$ for the constant stepsize. Now we show that Finite Horizon stepsize rule achieves a faster convergence rate at the pre-determined $T$-th step.

\begin{thm}
\label{theorem1} ({\bf Upper bound for Finite Horizon stepsize rule})
    Consider problem \eqref{eq_lp_no_inequality} with $A \in \mathbb{R}^{m\times n}$ full row rank and $m \leq n$. Assume the  singular values of $A$ lie in $[\mu, L]$, where $\mu, L >0$. Suppose $b \in \operatorname{range}(A)$ and  $c \in \operatorname{range}(A^\top)$ so that \eqref{eq_lp_no_inequality} is feasible and bounded. Consider the primal-dual method with Finite Horizon stepsize rule $(\eta_1, ... \eta_T)$ as in Algorithm \ref{algo_stepsize_rule}.
    Denote $z^\top = (x^\top, y^\top) \in \mathbb{R}^{n+m}$ as the concatenation of the primal variable $x \in  \mathbb{R}^{n}$ and dual variable $y \in  \mathbb{R}^{m}$. Assume the augmented Lagrangian coefficient $\beta$ is large such that $\beta \geq \frac{2}{\mu}$, then for any $z^0$ and any pre-determined $T$, Algorithm \ref{algo_stepsize_rule} achieves the following rate at the $T$-th step:  
    \begin{equation}
    \label{eq_finite_horizon_rate}
        \text{dist}(z^T, \mathcal{Z}^*) \leq \sqrt{2+ 4\gamma}\left(1 -  \frac{2}{\sqrt{\frac{\beta L^2 + \sqrt{\beta^2 L^4 - 4 L^2}}{\beta \mu^2 -  \sqrt{ \beta^2 \mu^4 - 4 \mu^2}} } + 1}\right)^T \text{dist}(z^0, \mathcal{Z}^*),
    \end{equation}
    where $\gamma = 4 / (\beta^{2}\mu^2 - 4)$ and $\text{dist}(z, \mathcal{Z}^*)$ denotes the Euclidean distance of $z$ to the set of optimal primal-dual solutions. Eq. \eqref{eq_finite_horizon_rate} can also be rewritten as follows: 
\begin{equation}
\label{eq_finite_horizon_rate_kappa}
        \text{dist}(z^T, \mathcal{Z}^*) \leq \sqrt{2+ 4\gamma}\left(1  - \frac{2}{\sqrt{\kappa} +1}\right)^T \text{dist}(z^0, \mathcal{Z}^*),
    \end{equation}
    where $\kappa$ as the condition number of matrix $M$ (as defined in \eqref{eq_M}).
\end{thm}

The proof of Theorem \ref{theorem1} is shown in Section
\ref{sec_proof_thm1}. Theorem \ref{theorem1} provides a complexity upper bound $\mathcal{O}(\sqrt{\kappa}\log(\frac{1}{\epsilon}))$ for Finite Horizon stepsize rule at $T$-th step, and the upper bound achieves acceleration over the lower bound $\Omega(\kappa \log(\frac{1}{\epsilon}) )$ for the constant step. It is important to note that the acceleration is reached  without changing the algorithm, e.g.,
storing extra building blocks like momentum, and the advantage merely comes from the stepsize choice. Meanwhile, the acceleration in Theorem \ref{theorem1} is only
guaranteed at the $T$-th iteration, which is fixed and finite. This is different from the classical theory that focuses on asymptotic guarantees. Our analysis is inspired by the pioneer work of GD stepsize \citep{young1953richardson}, yet we encounter new challenges raised by the non-symmetric nature of $M$. We provide a detailed explanation as follows.

\section{Proof of Main Results}
\subsection{Proof of Proposition \ref{prop_constant_lower_bound}}
\label{sec_proof_prop_constant_lower_bound}

The proof of Proposition \ref{prop_constant_lower_bound}
is just two steps. First, by Thereom 1.1 in \citep{young2014iterative}, we have $\eqref{eq_def_optimal_constant_stepsize} = \max_{A, s.t., \sigma(A)\subseteq [\mu,L], z \not\in \mathcal{Z}^*}  \rho^T$ where $\rho := \max_{\lambda >0, \lambda\in\lambda( I - \eta M)} |\lambda|.$ Then, it is easy to see that  $\eta^* = \frac{1}{\lambda_1 + \lambda_{2m}}$, where $\lambda_1, \lambda_{2m}$ are the largest and smallest non-zero eigenvalues of $M$, respectively.
Second, when initializing at $z^0 - z^* = x_{2m}$, where $x_{2m}$ is the eigenvector of $M$ associated with  $\lambda_{2m}$, we have the following equality and this concludes the proof. 
$$ 
 \text{dist}(z^T, \mathcal{Z}^*) = \left( 1 - \eta^* \lambda_{2m} \right)^T   \text{dist}(z^0, \mathcal{Z}^*) = \left( 1 - \frac{2}{\kappa + 1}\right)^T   \text{dist}(z^0, \mathcal{Z}^*).
$$

\subsection{Proof of Theorem \ref{theorem1}}
\label{sec_proof_thm1}

We now prove Theorem \ref{theorem1}. To simplify notations, we define \( \Gamma \) as follows. To prove Theorem \ref{theorem1}, we need to bound the singular values of $\Gamma$.
\begin{equation}
    \label{eq_gamma}
    \Gamma =(I - \eta_T M)(I - \eta_{T-1} M) \cdots (I - \eta_1 M).
\end{equation}

There are at least two difficulties to prove Theorem \ref{theorem1} (or to bound the singular of $\Gamma$). First, $M$ is a non-symmetric matrix, and there are not many helpful tools to bound the singular values of {\it the product of a series of non-symmetric matrices}. One potential way is to calculate the eigenvalues of $\Gamma^\top \Gamma$, but it is still not easy since $M^\top$ and $M$ do not commute and cannot be simultaneously diagonalized. Second, even if we explicitly calculate the singular values of $\Gamma$, they would have complicated dependency on $(\eta_1,\cdots, \eta_T)$, and it is not clear how to relate these explicit expressions to the condition number $\kappa$ of $M$.

The proof of Theorem \ref{theorem1} is summarized in the following four steps. 

\begin{itemize}[topsep=1pt,parsep=1pt,partopsep=1pt, leftmargin=*]
    \item {\bf Step 1.} We explicitly calculate the eigenvalues and eigenvectors of $M$.  This serves as the foundation of the whole proof.
    \item {\bf Step 2.}  We calculate the eigenvalues of $\Gamma^\top\Gamma$. The conventional way is to diagonalize $\Gamma^\top \Gamma$ using the eigenvectors of $M$. Although this method is useful when $M = M^\top$, we find that the diagonalization process is rather difficult when $M$ is non-symmetric, primarily because $M^\top$ and $M$ cannot be simultaneously diagonalized.   Here, we propose a rather unconventional way: First,  we  {\it block-diagonalize} $\Gamma^\top \Gamma$ into $2\times 2$ blocks. Then, we calculate the eigenvalues of each block via characteristic polynomials.
    \item {\bf Step 3.} We bound the eigenvalues of $\Gamma^\top \Gamma$ by the polynomial of $\lambda(M)$.
    \item {\bf Step 4.} Using the optimality condition of Finite Horizon rule, we bound the polynomial of $\lambda(M)$ by the condition number $\kappa$ of $M$. This step is done using the properties of the Chebyshev polynomial.
\end{itemize}

\paragraph{Step 1.}  We prove the following lemma to obtain the eigenvalues and eigenvectors of $M$. 
The proof of Lemma \ref{lemma_M} can be seen in Section \ref{sec_proof_lemma_M}.

\begin{lem}
\label{lemma_M}
    Consider the update matrix of the primal-dual method $M$ as in \eqref{eq_M}. Suppose the singular values of $A$ are $0<\mu\leq\sigma_m\leq \cdots \leq \sigma_1 \leq L$ and suppose $\beta$ is large such that $\beta \geq \frac{2}{\mu}$, then $M$ is diagonalizable with $2m$ non-zero eigenvalues and $(n-m)$ zero eigenvlaues. In particular, all $2m$ nonzero eigenvalues and the corresponding eigenvectors of $M$ are shown as follows.
    \begin{equation}
        \label{eq_eigen_M}
        \lambda_{i,1},  \lambda_{i,2} =  \frac{ \beta \sigma_{i}^2 \pm \sqrt{ \beta^2 \sigma_{i}^4 - 4 \sigma_{i}^2 } }{2}, x_{i,1} = \begin{bmatrix} v_{i} \\ -\frac{1}{\lambda_{i,1}} A v_i \end{bmatrix},  x_{i,2} = \begin{bmatrix} v_i \\ -\frac{1}{\lambda_{i,2}} A v_i \end{bmatrix}, i = 1, \cdots, m,
    \end{equation}
    where $v_i$ are the eigenvectors of $A^\top A$ associated with the nonzero eigenvalue $\sigma_i^2$, $i = 1,\cdots, m$.
    For the $(n-m)$  zero eigenvalues of $M$, the corresponding 
    eigenvectors are 
    \begin{equation}
        \label{eq_null_M}
        x_i =  \begin{bmatrix} v_i \\ 0 \end{bmatrix}, i = 2m+1, \cdots, n+m, 
    \end{equation}
    where $v_{2m+1}, \cdots, v_{n+m}$ form the basis of $\operatorname{null}(A^\top A)$.    
\end{lem}

\paragraph{Step 2.} Now we explicitly calculate the eigenvalues of $\Gamma^\top \Gamma$ via block-diagonalization.   

\begin{lem}
\label{lemma_block_diag_gamma}
    Consider the matrix $M$  and $\Gamma$ as defined in \eqref{eq_M} and \eqref{eq_gamma}, respectively. Denote the eigen-pairs of $M$ using the same notations as in Lemma \ref{lemma_M}. Define $X:= [x_{1,1}, \cdots, x_{m,2}, x_{2m+1},\cdots, x_{n+m}]  \in \mathbb{R}^{(n+m\times (n+m)}$ as the matrix of eigenvectors of $M$, then we have 
\[X^{-1} \Gamma^{\top} \Gamma X = 
 \left[\begin{array}{cccc}
  B_1^2  & \cdots & 0  & 0\\
   &   \ddots &  &  \\
 0 &  \cdots  & B_m^2 &0  \\
  0 &  \cdots & 0 & I_{(n-m) \times (n-m)} \\
\end{array}\right]_{(n+m) \times (n+m)},
\]
  where   $B_i :=  
\frac{1}{\lambda_{i, 1} - \lambda_{i, 2}} \begin{bmatrix}
  \delta_{i,1} (\lambda_{i, 1} + \lambda_{i, 2}) & -2 \delta_{i,1} \lambda_{i, 2}  \\
   2\delta_{i,2} \lambda_{i, 1}  & -\delta_{i,2} (\lambda_{i, 1} + \lambda_{i, 2}) 
\end{bmatrix} \in \mathbb{R}^{2\times 2}$,  $\delta_{i,k} = \prod_{t=1}^T (1 - \eta_t \lambda_{i,k})$, $i = [m]$, $k = [2]$.

Further, the eigenvalues of $B_i^2$ are

\begin{equation}
\label{eq_lambda_bi}
\left\{\begin{array}{l}
\lambda_1(B_i^2) = \frac{1}{2}\left(\delta_{i,1}^{2} + \delta_{i,2}^{2} + \varepsilon_i - \sqrt{(\delta_{i,1}^{2} - \delta_{i,2}^{2})^{2} + 2\varepsilon_i(\delta_{i,1}^{2} + \delta_{i,2}^{2}) + \varepsilon_i^{2}}\right) \\
\lambda_2(B_i^2) = \frac{1}{2}\left(\delta_{i,1}^{2} + \delta_{i,2}^{2} + \varepsilon_i + \sqrt{(\delta_{i,1}^{2} - \delta_{i,2}^{2})^{2} + 2\varepsilon_i(\delta_{i,1}^{2} + \delta_{i,2}^{2}) + \varepsilon_i^{2}}\right),
\end{array}\right.
\end{equation}
where
\[
\varepsilon_i = \frac{4\lambda_{i, 1}\lambda_{i, 2}}{(\lambda_{i, 1} - \lambda_{i, 2})^{2}}(\delta_{i,1} - \delta_{i,2})^{2} = \frac{4\sigma_{i}^2}{\beta^{2}\sigma_{i}^{4} - 4\sigma_{i}^2}(\delta_{i,1} - \delta_{i,2})^{2} = \frac{4}{\beta^{2}\sigma_{i}^2 - 4}(\delta_{i,1} - \delta_{i,2})^{2}.
\]
\end{lem}

The proof of Lemma \ref{lemma_block_diag_gamma} is shown in Section \ref{sec_proof_lemma_block_diag_gamma}.  Note that the eigenvalues of $B_i^2$ in \eqref{eq_lambda_bi} are also the eigenvalues of $\Gamma^\top \Gamma$ because of the block-diagonalization. 

\paragraph{Step 3.} We now bound $\lambda(\Gamma^\top \Gamma)$  using the polynomial of $\lambda(M)$.
Since $\varepsilon_i \geq 0$, we have $r_{i,1} \leq r_{i,2}$. We notice a bound
\[
\varepsilon_i \leq \gamma (\delta_{i,1} - \delta_{i,2})^{2} \leq 2\gamma (\delta_{i,1}^{2} + \delta_{i,2}^{2}),
\]
for some constant $\gamma = 4 / (\beta^{2}\mu^2 - 4)$  not depending on the iteration number $T$. The bound implies that
\[
r_{i,2} \leq \delta_{i,1}^{2} + \delta_{i,2}^{2} + 2\gamma(\delta_{i,1}^{2} + \delta_{i,2}^{2}) \leq (2 + 4\gamma)\max(\delta_{i,1}^{2}, \delta_{i,2}^{2}).
\]
Therefore, we have
\begin{equation}
\label{eq_gamma_upper_bound}
\max_{\lambda \in \lambda(\Gamma^\top\Gamma), \lambda \neq 1} |\lambda|^{\frac{1}{2}} \leq \sqrt{2 + 4\gamma}\max_{i\in [2m]}(|\delta_{i,1}|, |\delta_{i,2}|) \leq \sqrt{2 + 4\gamma} \max_{\lambda \in [\lambda_{2m}, \lambda_1]} \left|(1-\eta_1 \lambda) \cdot (1 - \eta_{2}\lambda) \dots (1 - \eta_{T}\lambda)\right|.
\end{equation}

\paragraph{Step 4.} Finally, we bound \eqref{eq_gamma_upper_bound} using the properties of Chebyshev polynomial.
Before that, we first write down the convergence rate using the optimality condition of Finite Horizon stepsize rule in {\bf Algorithm \ref{algo_stepsize_rule}}.
\begin{eqnarray}
\text{dist}(z^T, \mathcal{Z}^*) &\leq &\max_{A, s.t., \sigma(A)\subseteq [\mu,L], z \not\in \mathcal{Z}^*}\frac{\|(I -\eta_T M) \cdots ( I - \eta_1 M)(z - z^*) \|_{2}}{\|z - z^*\|_2} \text{dist}(z^0, \mathcal{Z}^*) \nonumber \\
& \overset{\text{Algorithm \ref{algo_stepsize_rule}}}{ = } & \min_{\eta_1,\cdots,\eta_T} \max_{A, s.t., \sigma(A)\subseteq [\mu,L], z \not\in \mathcal{Z}^*}\frac{\|(I -\eta_T M) \cdots ( I - \eta_1 M)(z - z^*) \|_{2}}{\|z - z^*\|_2} \text{dist}(z^0, \mathcal{Z}^*) \nonumber  \\
& \overset{\eqref{eq_gamma}}{ = } & \min_{\eta_1,\cdots,\eta_T} \max_{A, s.t., \sigma(A)\subseteq [\mu,L], z \not\in \mathcal{Z}^*}\frac{\|\Gamma(z - z^*) \|_{2}}{\|z - z^*\|_2} \text{dist}(z^0, \mathcal{Z}^*) \nonumber  \\
& \overset{}{ = } & \min_{\eta_1,\cdots,\eta_T} \max_{A, s.t., \sigma(A) \in [\mu, L], \lambda \in \lambda(\Gamma^\top\Gamma), \lambda \neq 1}  |\lambda|^{\frac{1}{2}} \text{dist}(z^0, \mathcal{Z}^*) \nonumber  \\
& \overset{\eqref{eq_gamma_upper_bound}}{\leq} & \sqrt{2 + 4\gamma} \underbrace{\min_{\eta_1,\cdots,\eta_T} \max_{\lambda \in [\lambda_{2m}, \lambda_1]}   \left|(1-\eta_1 \lambda) \cdot (1 - \eta_{2}\lambda) \dots (1 - \eta_{T}\lambda)\right|}_{(a)} \text{dist}(z^0, \mathcal{Z}^*). \nonumber 
\end{eqnarray}

. We denote $f_T(\lambda) = (1-\eta_1 \lambda) \cdot (1 - \eta_{2}\lambda) \dots (1 - \eta_{T}\lambda)$. We observe that $ f_T(\lambda ) $ is a $T$-th order polynomial with constraint $f_T(0) = 1$. We now bound  $(a)$ using the help of the Chebyshev polynomial \citep{chebyshev1853theorie}. We first recall a classical proposition of the Chebyshev polynomial from the literature \citep{markoff1916polynome,pedregosa2021residual}, and then use it to bound $(a)$.

\begin{prop}
\label{prop_chebyshev_minimax} ({\bf Minimax polynomial}) 
    Define the linear scaling mapping $\sigma(\lambda)=\frac{\lambda_{1}+\lambda_{2m}}{\lambda_1-\lambda_{2m}}-\frac{2}{\lambda_1-\lambda_{2m}} \lambda,$ which maps the positive interval $[\lambda_{2m}, \lambda_1]$ into $[-1,1]$. Denote the $t$-th order Chebyshev polynomial  as $C_t(x)=\cos (t * \arccos x)$, then we have
    \begin{equation}
        \label{eq_chebyshev_minimax}
        \frac{C_t(\sigma(\lambda))}{C_t(\sigma(0))}=\underset{f_t \in \mathcal{F}_t}{\arg \min } \max _{\lambda \in[\lambda_{2m}, \lambda_1]}\left|f_t(\lambda)\right|, \quad s.t., f_t(0)=1,
    \end{equation}

where $\mathcal{F}_t$ denotes the set of $t$-th order polynomials.   
\end{prop}

The proof of Proposition \ref{prop_chebyshev_minimax} is relegated in Section \ref{sec_proof_chebyshev}.
Now we bound $(a)$ using $\frac{C_t(\sigma(\lambda))}{C_t(\sigma(0))}$.

\begin{eqnarray}
\text{dist}(z^T, \mathcal{Z}^*) 
& \overset{\eqref{eq_gamma_upper_bound}}{\leq} & \sqrt{2 + 4\gamma} \underbrace{\min_{\eta_1,\cdots,\eta_T} \max_{\lambda \in [\lambda_{2m}, \lambda_1]}   \left|(1-\eta_1 \lambda) \cdot (1 - \eta_{2}\lambda) \dots (1 - \eta_{T}\lambda)\right|}_{(a)} \text{dist}(z^0, \mathcal{Z}^*). \nonumber  \\
& \overset{\eqref{eq_chebyshev_minimax}}{=} & \sqrt{2 + 4\gamma} \frac{\sup _{\lambda \in[\lambda_{2m}, \lambda_1]}\left|C_T(\sigma(\lambda))\right|}{\left|C_T(\sigma(0))\right|} \text{dist}(z^0, \mathcal{Z}^*) \\
& \overset{(a)}{=} &   \frac{\sqrt{2 + 4\gamma} }{C_T\left(\frac{\kappa + 1}{\kappa -1}\right)} \text{dist}(z^0, \mathcal{Z}^*)\\
&\overset{(b)}{\leq}& \frac{2 \sqrt{2 + 4\gamma} }{  \left(\frac{\kappa + 1}{\kappa - 1}+\sqrt{\left(\frac{\kappa+1}{\kappa - 1}\right)^2-1}\right)^t   } \text{dist}(z^0, \mathcal{Z}^*)\\
&\overset{(c)}{\leq}&\sqrt{2 + 4\gamma} \left(\frac{\sqrt{\kappa} - 1}{\sqrt{\kappa} + 1}\right)^T  \text{dist}(z^0, \mathcal{Z}^*)\\
&=&\sqrt{2 + 4\gamma} \left(1- \frac{2}{1+\sqrt{\kappa}}\right)^T  \text{dist}(z^0, \mathcal{Z}^*),
\end{eqnarray}

where $(a)$: $\sup _{x \in[-1,1]}\left|C_T(x)\right|=1$, $\sigma(0)=\frac{\lambda_1+\lambda_{2m}}{\lambda_1-\lambda_{2m}}$ and we denote $\kappa = \frac{\lambda_1}{\lambda_{2m}}$; $(b)$:  the property of Chebyshev polynomial \citep{pedregosa2021residual}: $C_T(x) \geq \frac{\left(x+\sqrt{x^2-1}\right)^T}{2}, \forall x \notin(-1,1)$; $(c)$: some basic re-arrangement.

\subsection{Proof of Lemma \ref{lemma_M}}
\label{sec_proof_lemma_M}

We solve the following equation to get the eigenpairs of $M$:
\[
M \begin{bmatrix} v \\ u \end{bmatrix} \overset{\eqref{eq_M}}{=}  \left[\begin{array}{cc}
 \beta A^\top A & A^\top \\
-A & 0_{m\times m}
\end{array}\right]  \begin{bmatrix} v \\ u \end{bmatrix} 
 =       \lambda \begin{bmatrix} v \\ u \end{bmatrix},
\]
where \( v \in \mathbb{R}^n \), \( u \in \mathbb{R}^m \), and \( \lambda \) is an eigenvalue of \( M \).
This yields:
\begin{align}
\beta A^\top A v + A^\top u &= \lambda v, \label{eq1} \\
- A v &= \lambda u. \label{eq2}
\end{align}

We first consider the case with $\lambda \neq 0$. 
From \eqref{eq2}, we have:
\begin{equation}
\label{eq_lemma1_y}
    u = -\frac{1}{\lambda} A v.
\end{equation}

Substituting \( u \) back into equation \eqref{eq1}, we have:
\[
\beta A^\top A v - \frac{1}{\lambda} A^\top A v = \lambda v,
\]
or equivalently, we have:
\[
\left( \beta - \frac{1}{\lambda} \right) A^\top A v = \lambda v.
\]
Since \( A^\top A v = \sigma_i^2 v \) for eigenvalues \( \sigma_i^2 \) of \( A^\top A \), we get:
\[
\left( \beta - \frac{1}{\lambda} \right) \sigma_i^2 v = \lambda v.
\]
This leads to the scalar equation:
\[
\left( \beta - \frac{1}{\lambda} \right) \sigma_i^2 = \lambda.
\]
Multiplying both sides by \( \lambda \):
\[
\left( \beta \lambda - 1 \right) \sigma_i^2 = \lambda^2.
\]
Solve this quadratic equation in \( \lambda \) and we will get:
\[
\lambda_{i,1}, \lambda_{i,2} = \frac{ \beta \sigma_i^2 \pm \sqrt{ \beta^2 \sigma_i^4 - 4 \sigma_i^2 } }{2}, i = 1, \cdots ,m.
\]
To get the corresponding eigenvectors,  we use \eqref{eq_lemma1_y} and conclude that 
    \[ x_{i,1} = \begin{bmatrix} v_i \\ -\frac{1}{\lambda_{i,1}} A v_i \end{bmatrix},  x_{i,2} = \begin{bmatrix} v_i \\ -\frac{1}{\lambda_{i,2}} A v_i \end{bmatrix}, i = 1, \cdots, m. \]

We now calculate the eigenvectors corresponding to 0 eigenvalues. 
In this case, we have $A u = 0 $ from \eqref{eq2}. Plugging it into \eqref{eq1} and we have
\[A^\top u = 0.\]
Since $A \in \mathbb{R}^{m\times n}$ is full row rank with $m <n$, we have $u = 0$. Therefore, the eigenvectors for the zero eigenvalues are 
 \[x_i =  \begin{bmatrix} v_i \\ 0 \end{bmatrix}, i = 2m+1, \cdots, n+m, \]
    where $v_{2m+1}, \cdots, v_{n+m}$ form the basis of $\operatorname{null}(A)$, or equivalently, the basis of $\operatorname{null}(A^\top A)$.   Finally, we note that all these eigenvectors $x_{i,1}, x_{i,2}, i = 1,\cdots, m$  and $x_{2m+1},\cdots, x_{n+m}$ are linearly independent, so $M$ is diagonalizable.

\subsection{Proof of Lemma \ref{lemma_block_diag_gamma}}
\label{sec_proof_lemma_block_diag_gamma}

We define the matrix \( \Omega \) as:
\[
\Omega = \begin{bmatrix}
I_{n\times n} & 0 \\
0 & -I_{m\times m}
\end{bmatrix}.
\]
Note that \( \Omega \) is involutive (\( \Omega^2 = I \)) and symmetric (\( \Omega^\top = \Omega \)), and it is easy to see that  $M^\top = \Omega M \Omega$. 
Now we investigate how \( \Omega \) affects the eigenvectors of \( M \).
Let \(x_{i,1}, x_{i,2}\) be eigenvectors of \( M \) corresponding to eigenvalues \( \lambda_{i,1},\lambda_{i,2} \) as defined in \eqref{eq_eigen_M}.
Applying \( \Omega \) to \( x_{i,1}, x_{i,2}\) and we get:
\[
\Omega x_{i,k} = \begin{bmatrix} v_i \\ \frac{1}{\lambda_{i, k}}Av_{i} \end{bmatrix}, k = 1, 2.
\]
Let us express \( \Omega x_{i,k} \) in terms of \( x_{i,1} \) and \( x_{i,2} \). We can write:
\[
\Omega x_{i,1}  = \begin{bmatrix} v_i \\ \frac{1}{\lambda_{i, 1}} Av_{i} \end{bmatrix} = c_{i, 1} \begin{bmatrix} v_i \\ -\frac{1}{\lambda_{i, 1}} Av_{i} \end{bmatrix} + c_{i, 2} \begin{bmatrix} v_i \\ -\frac{1}{\lambda_{i, 2}} Av_{i} \end{bmatrix} = c_{i,1} x_{i,1} + c_{i,2} x_{i,2}
\]
\[
\Omega x_{i,2}  =  \begin{bmatrix} v_i \\ \frac{1}{\lambda_{i, 2} } Av_{i} \end{bmatrix} = d_{i, 1} \begin{bmatrix} v_i \\ -\frac{1}{\lambda_{i, 1}} Av_{i} \end{bmatrix} + d_{i, 2} \begin{bmatrix} v_i \\ -\frac{1}{\lambda_{i, 2}} Av_{i} \end{bmatrix}  = d_{i,1} x_{i,1} + d_{i,2} x_{i,2}
\]
We determine the coefficients \( c_{i, 1}, c_{i, 2} \) and \( d_{i, 1}, d_{i, 2} \) by solving these equations. The results are shown below.

\begin{equation}
\label{eq_coefficient_c_d}
     \begin{bmatrix}
  c_{i,1}& c_{i,2} \\
   d_{i,1} & d_{i,2}
\end{bmatrix} = 
\frac{1}{\lambda_{i, 1} - \lambda_{i, 2}} \begin{bmatrix}
  \lambda_{i, 1} + \lambda_{i, 2} & -2\lambda_{i, 2} \\
   2\lambda_{i, 1} & -(\lambda_{i, 1} + \lambda_{i, 2})
\end{bmatrix}
\end{equation}

Since \(x_{i,k}\) is the eigenvectors of \( M \) associated with eigenvalue \(\lambda_{i,k} \),  we have 
\[
\Gamma x_{i,k} = \prod_{t=1}^T (1 - \eta_t \lambda_{i,k}) x_{i,k}.
\]
Denote $\delta_{i,k} = \prod_{t=1}^T (1 - \eta_t \lambda_{i,k})$. 
Using the coefficients in \eqref{eq_coefficient_c_d}, 
we can express \( \Omega \Gamma x_{i,k} \) in terms of \( x_{i,k}, k = 1,2 \).

\[
\Omega \Gamma x_{i,1}  = \delta_{i,1} \begin{bmatrix} v_i \\ \frac{1}{\lambda_{i, 1}} Av_{i} \end{bmatrix} = \delta_{i,1} c_{i, 1}  \begin{bmatrix} v_i \\ -\frac{1}{\lambda_{i, 1}} Av_{i} \end{bmatrix} + \delta_{i,1} c_{i, 2} \begin{bmatrix} v_i \\ -\frac{1}{\lambda_{i, 2}} Av_{i} \end{bmatrix} = \delta_{i,1} c_{i,1} x_{i,1} +  \delta_{i,1} c_{i,2} x_{i,2}.
\]
\[
\Omega \Gamma x_{i,2}  = \delta_{i,2}  \begin{bmatrix} v_i \\ \frac{1}{\lambda_{i, 2} } Av_{i} \end{bmatrix} = \delta_{i,2} d_{i, 1} \begin{bmatrix} v_i \\ -\frac{1}{\lambda_{i, 1}} Av_{i} \end{bmatrix} + \delta_{i,2} d_{i, 2} \begin{bmatrix} v_i \\ -\frac{1}{\lambda_{i, 2}} Av_{i} \end{bmatrix}  = \delta_{i,2} d_{i,1} x_{i,1} +  \delta_{i,2} d_{i,2} x_{i,2}.
\]

Define 
$B_i :=  
\frac{1}{\lambda_{i, 1} - \lambda_{i, 2}} \begin{bmatrix}
  \delta_{i,1} (\lambda_{i, 1} + \lambda_{i, 2}) & -2 \delta_{i,1} \lambda_{i, 2}  \\
   2\delta_{i,2} \lambda_{i, 1}  & -\delta_{i,2} (\lambda_{i, 1} + \lambda_{i, 2}) 
\end{bmatrix} \in \mathbb{R}^{2\times 2},$ $i = [m]$, and $X:= [x_1, \cdots, x_{n+m}]  \in \mathbb{R}^{(n+m)\times (n+m)} $. We have 

\[
\Omega \Gamma X  = X \left[\begin{array}{cccc}
  B_1  & \cdots & 0  & 0\\
   &   \ddots &  &  \\
 0 &  \cdots  & B_m &0  \\
  0 &  \cdots & 0 & I_{(n-m) \times (n-m)} \\
\end{array}\right]_{(n+m) \times (n+m)}. 
\]

Therefore, we have 

\[X^{-1} \Gamma^{\top} \Gamma X = 
X^{-1}\Omega \Gamma  \Omega \Gamma X  = \left[\begin{array}{cccc}
  B_1^2  & \cdots & 0  & 0\\
   &   \ddots &  &  \\
 0 &  \cdots  & B_m^2 &0  \\
  0 &  \cdots & 0 & I_{(n-m) \times (n-m)} \\
\end{array}\right]_{(n+m) \times (n+m)}. 
\]

To bound the largest non-one eigenvalue of $\Gamma^\top \Gamma$, we need to calculate the eigenvalues of $B_i^2$ and find the largest one. 
 We calculate the roots of the characteristic polynomial of $B_i^2$ and we get:

\[
\left\{\begin{array}{l}
\lambda_1(B_i^2) = \frac{1}{2}\left(\delta_{i,1}^{2} + \delta_{i,2}^{2} + \varepsilon_i - \sqrt{(\delta_{i,1}^{2} - \delta_{i,2}^{2})^{2} + 2\varepsilon_i(\delta_{i,1}^{2} + \delta_{i,2}^{2}) + \varepsilon_i^{2}}\right) \\
\lambda_2(B_i^2) = \frac{1}{2}\left(\delta_{i,1}^{2} + \delta_{i,2}^{2} + \varepsilon_i + \sqrt{(\delta_{i,1}^{2} - \delta_{i,2}^{2})^{2} + 2\varepsilon_i(\delta_{i,1}^{2} + \delta_{i,2}^{2}) + \varepsilon_i^{2}}\right),
\end{array}\right.
\]

where
\[
\varepsilon_i = \frac{4\lambda_{i, 1}\lambda_{i, 2}}{(\lambda_{i, 1} - \lambda_{i, 2})^{2}}(\delta_{i,1} - \delta_{i,2})^{2} = \frac{4\sigma_{i}^2}{\beta^{2}\sigma_{i}^{4} - 4\sigma_{i}^2}(\delta_{i,1} - \delta_{i,2})^{2} = \frac{4}{\beta^{2}\sigma_{i}^2 - 4}(\delta_{i,1} - \delta_{i,2})^{2}.
\]

\subsection{Proof of Proposition \ref{prop_chebyshev_minimax}}
\label{sec_proof_chebyshev}

Proposition  \ref{prop_chebyshev_minimax} is classical results of Chebyshev polynomial \citep{markoff1916polynome,pedregosa2021residual}. For completeness of the whole proof,  we restate the proof of this classical result under our notations.  To prove Proposition  \ref{prop_chebyshev_minimax}, we need to first prove the following Proposition \ref{prop_chebyshev_equi}.

\begin{prop}({\bf Equioscillation}) 
\label{prop_chebyshev_equi}
Consider  the $t$-th order Chebyshev polynomial $C_t(x)=\cos (t * \arccos x)$, then for $x \in[-1,1]$, there exists $x_0<x_1<\ldots<x_t \in[-1,1]$, s.t. $C_t\left(x_k\right)=(-1)^k, k = 0, \cdots n$. In other words, $T_t(x)$ has $(t+1)$ extreme points in $x \in[-1,1]$.
\end{prop}

The proof of Proposition \ref{prop_chebyshev_equi} is just one line: for any $k=0,1, \cdots t$, define $x_k=\cos \left(\frac{k \pi}{t}\right)$, then we have the following equation and conclude the proof. 

$$C_t\left(x_k\right)=\cos \left(t * \arccos \left(\cos \left(\frac{k \pi}{t}\right)\right)\right)=\cos (k \pi)=(-1)^k.$$

  We now prove Proposition \ref{prop_chebyshev_minimax} by contradiction. We will show that if the problem \eqref{eq_chebyshev_minimax} has a better solution, then that solution must at least be a $(t+1)$-th degree polynomial.

  Since $\sigma(\lambda)$ is a linear translation from $[\lambda_{2m}, \lambda_1]$ into $[-1,1]$, $\frac{C_t(\sigma(\lambda))}{C_t(\sigma(0))}$ has $(t+1)$ extreme points on $[\lambda_{2m}, \lambda_1]$.
  By the Equioscillation property of Chebyshev polynomial (Proposition \ref{prop_chebyshev_equi}), the image at these extreme points have the same absolute values and they are alternately positive and negative.
 Now suppose there exists $t$-th order polynomial $R_t(\lambda)$, s.t.: $R_t(0)=1$ but $R_t(\lambda)$ has smaller maximum absolute value.
Consider $Q(\lambda) \stackrel{\text { def }}{=}\frac{C_t(\sigma(\lambda))}{C_t(\sigma(0))}-R_t(\lambda),$ $ Q(\lambda)$ still alternately $>0$ or $<0$ at all $(t+1)$ extreme points, so $Q$ must have $t$ zeros in $[\lambda_{2m}, \lambda_1]$.
However, $Q(0)=P_t(0)-R_t(0)=0$, so $Q$ has $(t+1)$ zeros. This forms a Contradiction. Therefore,$\frac{C_t(\sigma(\lambda))}{C_t(\sigma(0))}$ is the optimal solution to the constrained minimax problem \eqref{eq_chebyshev_minimax}.

\section{Experiments}
\label{sec_experiments}

Now we verify the effectiveness of Finite Horizon stepsize rule ({\bf Algorithm \ref{algo_stepsize_rule}}) of the primal-dual method on various LP instances. All experiments are conducted on Intel(R) Xeon(R) Gold 6226 CPU @ 2.70GHz with 48 cores. 

{\bf Setups.} For all experiments, we choose $\beta = 4 / \mu$, where $\mu$ is the smallest non-zero singular value of $A$. This choice of $\beta$ satisfies the condition in Theorem \ref{theorem1}.
 To solve the SDP sub-problem in {\bf Algorithm \ref{algo_stepsize_rule}}, we use \texttt{SCS} solver in \texttt{cvxpy} library \footnote{\url{https://www.cvxgrp.org/scs/}} with the configuration $n_{\operatorname{sample}} = 200$, \texttt{alpha} =1.5, and \texttt{sdp\_iter} = 100 (the number of iterations of SDP solver to solve \eqref{eq_sdp}). We set this configuration as the default choice and will not change it unless mentioned otherwise.  For all experiments, we use Gaussian random initialization $\mathcal{N}(5,1)$ to initialize the primal and dual variables. We choose the mean value of Gaussian to be $5$ to keep the primal variable away from the origin and thus reduce the possibility of the projection.  We define the optimality gap as the relative KKT residue \citep{xiong2023computational},
, i.e., the sum of primal feasibility, dual feasibility, and primal-dual gap, in the relative sense: $\epsilon(x, y):=\frac{\left\|A x^{+}-b\right\|}{1+\|b\|}+\frac{\left\|\left(c-A^{\top} y\right)^{-}\right\|}{1+\|c\|}+\frac{\left|c^{\top} x^{+}-b^{\top} y\right|}{1+\left|c^{\top} x^{+}\right|+\left|b^{\top} y\right|} .$

\subsection{Case Study on the Toy Example: Intuition and Initial Results}
\label{sec_toy_example}

We first focus on a simple LP to gain intuition and boost understanding of Finite Horizon stepsize rule.
. We consider $c = A = 1, b = 200$, i.e.:

\begin{equation}
\label{eq_toy_lp}
\begin{aligned}
\min _{x \in \mathbb{R}^n} & \quad  x \\
\text { s.t. } &  x= 200 \\
& x \geq 0.
\end{aligned}
\end{equation}

For this example, we have $M = \left[\begin{array}{cc}
\beta & 1 \\
-1 & 0 \\
\end{array}\right] = \left[\begin{array}{cc}
4 & 1 \\
-1 & 0 \\
\end{array}\right]  $ under the choice of $\beta =\frac{4}{\mu}  = 4$. The optimal primal-dual solution is $(x^*, y^*) = (200, 1)$. 
 We consider two versions of optimal constant stepsize rules as the baseline methods.
  
\paragraph{Baseline method 1: $T$-step optimal constant stepsize.}   We consider the optimal constant stepsize  within $T$ total iterations, i.e., 

 \begin{equation}
 \label{eq_toy_example_minimax_T}
      \min_{\BLUE{\eta}}\max_{A, s.t., \sigma(A)\subseteq [\mu,L]}\|(I - \BLUE{\eta} M)^T  \|_{\operatorname{op}}. 
 \end{equation}
 
Since  $A = 1 \in \mathbb{R}$ in this toy example, the interval of $\sigma(A)$ degenerates into a single point $\mu = L = 1$ and the above minimax problem 
\eqref{eq_toy_example_minimax_T} reduces to the following problem.

 \begin{equation}
 \label{eq_toy_example_min_T_2}
\min_{\BLUE{\eta}}\|(I - \BLUE{\eta} M)^T  \|_{\operatorname{op}}. 
 \end{equation}

Problem \eqref{eq_toy_example_min_T_2}
 is a 1-dimensional  (nonconvex) optimization problem. It can be solved to near-global optimal solution using Stimulated Annealing methods \citep{tsallis1988possible, xiang2013generalized}. In this work, we use the algorithm described in \citep{xiang2013generalized} and use the implementation in \texttt{scipy} library \footnote{\url{https://docs.scipy.org/doc/scipy/reference/generated/scipy.optimize.dual_annealing.html}}.

Note that when $T = 1$,  \eqref{eq_toy_example_minimax_T} reduces to Finite Horizon stepsize rule for $T = 1$. We will see that the annealing methods to produce the same results as the SDP solver for \eqref{eq_sdp} in this case.

\paragraph{Baseline method 2: $\infty$-step optimal constant stepsize.} We consider the asymtopic optimal constant stepsize when $T \rightarrow \infty$.

 \begin{eqnarray}
 \label{eq_toy_example_minimax_inf}
     && \min_{\BLUE{\eta}}\max_{A, s.t., \sigma(A)\subseteq [\mu,L]} \lim_{T \rightarrow \infty}\|(I - \BLUE{\eta} M)^T\|_{\operatorname{op}}^{\frac{1}{T}}.  \nonumber\\
   &\overset{(a)}{=} & \min_{\BLUE{\eta}}\max_{A, s.t., \sigma(A)\subseteq [\mu,L]} \rho\left(I - \BLUE{\eta} M \right)  \nonumber\\
   & \overset{(b)}{=} & \frac{2}{\lambda_1 + \lambda_{2m}},
 \end{eqnarray}

where $\lambda_1$ and $\lambda_{2m} $ are the largest and smallest non-zero eigenvalues of $M$. $(a)$ is due to the property of spectral radius $\rho (A) =  \lim_{T \rightarrow \infty}\|A^T\|_{\operatorname{op}}^{\frac{1}{T}}$ \citep[Theorem 1.12]{saad2003iterative}.  $(b)$ relies on two properties of $M$: When the augmented Lagrangian coefficient satisfies $\beta \geq \frac{2}{\mu}$, we have: (i) $M$ is diagonalizable.  (ii) all eigenvalues of $M$ are real-valued. Further,    
 $\lambda_1$ and $\lambda_{2m}$ can both be calculated using the extreme singular values of  $A$, i.e., $\mu$ and $L$. We present their relation as below:

\begin{equation}
\label{eq_lambda_relation_sigma}
    \left\{\begin{array}{l}
\lambda_1 = \frac{\beta L^2}{2}  +  \sqrt{\frac{\beta^2 L^4 }{4} - L^2} \overset{(*)}{ = }  \frac{2 L^2}{\mu}  +  \sqrt{\frac{4 L^4 }{ \mu^2} - L^2}  \\
\lambda_{2m} =  \frac{\beta L^2}{2}  -  \sqrt{\frac{\beta^2 L^4 }{4} - L^2} \overset{(*)}{ = }  \frac{2 L^2}{\mu}  - \sqrt{\frac{4 L^4 }{ \mu^2} - L^2} ,
\end{array}\right.
\end{equation}

 $(*)$ holds under our default choice of $\beta  = \frac{4}{\mu}$. Detailed derivation and proof can be seen in the proof of Lemma \ref{lemma_M}.
  In the experiments, we refer to Baseline method 1 as ``optimal constant stepsize for $T=t$" and refer Baseline method 2 as  ``optimal constant stepsize".

  \begin{figure}[t]
    \centering
  \subfigure[Optimality gap]{\includegraphics[width=0.30\textwidth]{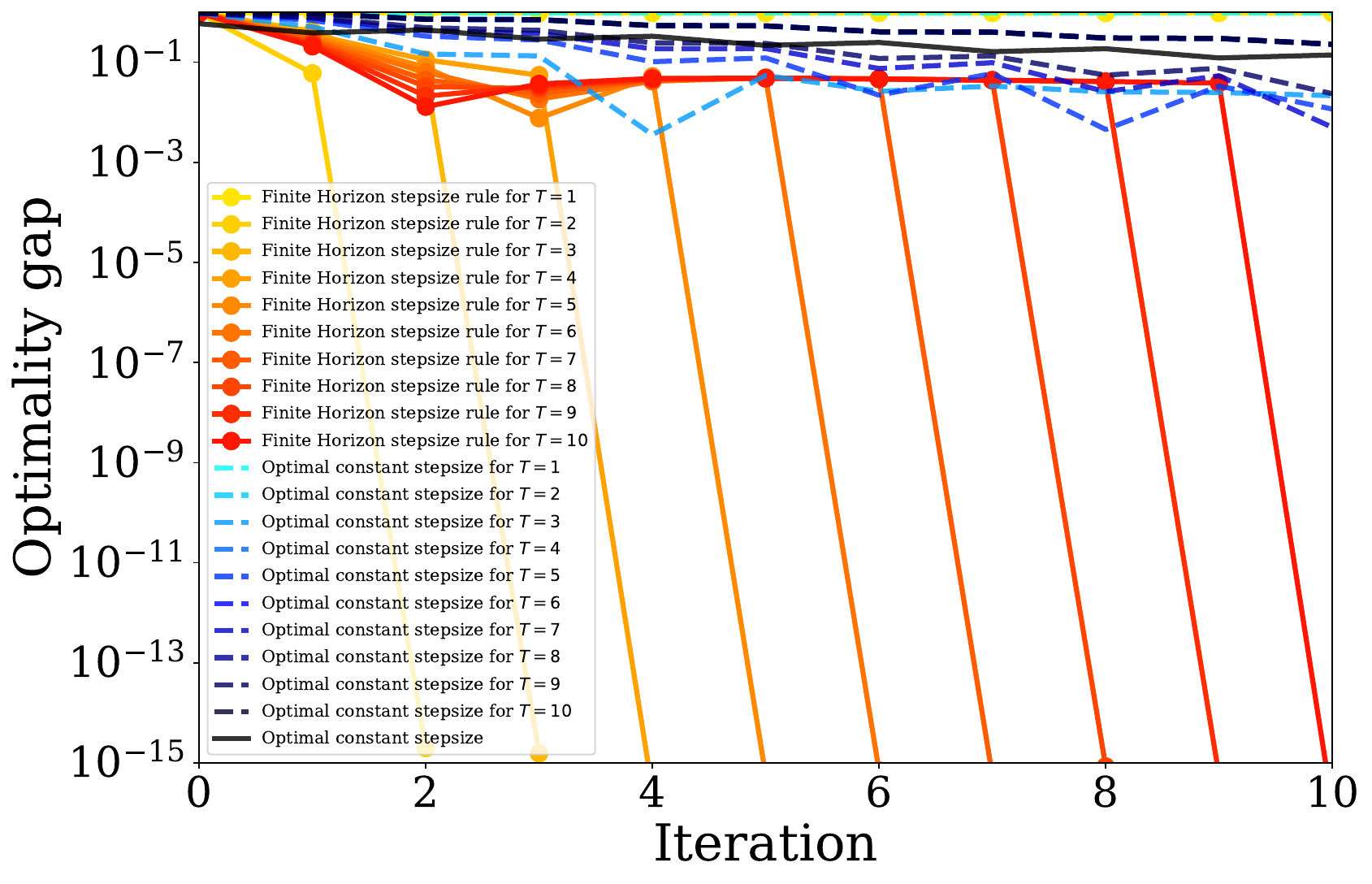}}
  \subfigure[Stepsize rules]{\includegraphics[width=0.30\textwidth]{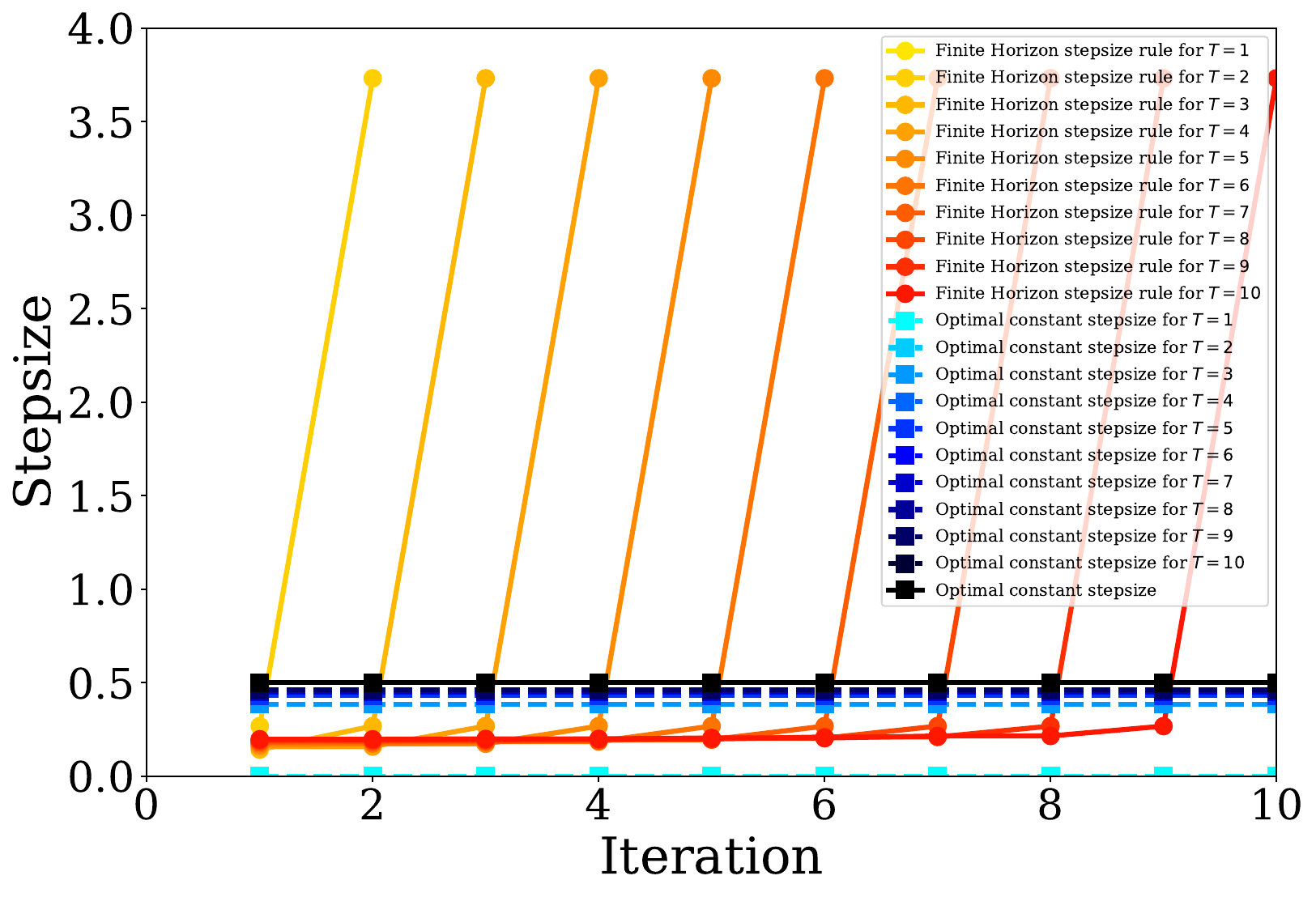}}
  \subfigure[Trajectories]{\includegraphics[width=0.30\textwidth]{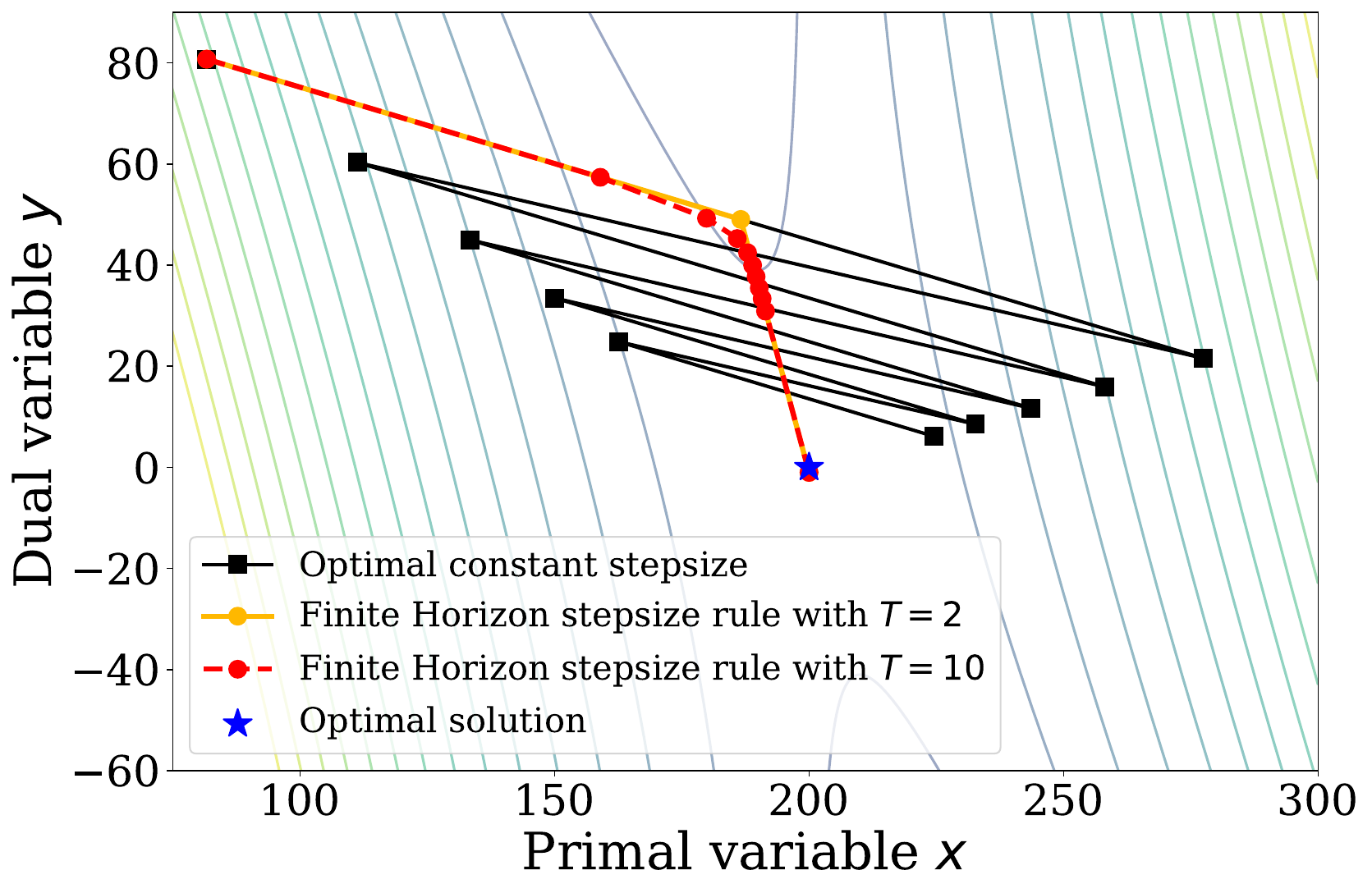}}
    \caption{Results on the toy LP instance \eqref{eq_toy_lp}. For Finite Horizon stepsize rule, we find that (a): the optimality gap will drop sharply at the $T$-th step; (b) there will be a sudden surge at the $T$-th step; (c): the iterates by Finite Horizon stepsize rules will take smaller steps in the sharp region, and then take large steps in the flat region.}
  \label{fig_toy_example}
  \vspace{-0.3cm}
  \end{figure}

\paragraph{Remark 1:} if $M$ is symmetric, then  Baseline method 1 and 2  are the same and both equal to $\eta^*  =  \frac{2}{\lambda_1 + \lambda_{n +m}}$. However, for non-symmetric $M$, these two baseline methods are different (see Figure \ref{fig_toy_example}). This is because Baseline method 1 controls singular values, while Baseline method 2 controls eigenvalues. For non-symmetric matrices, thse two are usually different (see discussion in Section \ref{sec_main_methods}).

\paragraph{Remark 2: Same information queried from LP class.} We emphasize that both methods query the same amount of problem information to determine the hyperparameters: they both query the largest \& smallest singular values of $A$, i.e., $\mu$ and $L$. However, they use this information in a different way. For the optimal constant stepsize, it uses  $\mu$ and $L$ to calculate stepsize using the formulas in \eqref{eq_lambda_relation_sigma}. For Finite Horizon stepsize rule, it uses  $\mu$ and $L$ to solve the SDP sub-problem \eqref{eq_sdp}.    In the following, we will show that Finite Horizon stepsize rule consistently outperforms the optimal constant stepsize.

\paragraph{Results on the toy example \eqref{eq_toy_lp}.}  In Figure \ref{fig_toy_example}, we run the primal-dual method with Finite Horizon stepsize rule for $T \leq 10$ and compare them to the corresponding optimal constant stepsize counterparts.
We find that: for Finite Horizon stepsize rule with the total iteration budget $T$,  the first $(T-1)$-th stepsizes will be smaller than the optimal constant step size, and then there will be a sudden surge at the $T$-th step.
Accordingly, the optimality gap decays similarly to that of the constant step size for the first $(T-1)$-th steps, but it will drop sharply at the $T$-th step and reach the optimal primal-dual solution.  This aligns with our theoretical predictions: Finite Horizon stepsize rule is designed to reach a fast convergence rate {\it at} $T$-th step, and there is no acceleration guarantee for steps {\it before} $T$.
In Figure \ref{fig_toy_example} (c), we see that the iterates by Finite Horizon stepsize rules will take smaller steps in the sharp region, and then take large steps in the flat region.

\begin{figure}[t]
  \centering
   \subfigure[$\pproj$ = 0.4]{\includegraphics[width=0.30\textwidth]{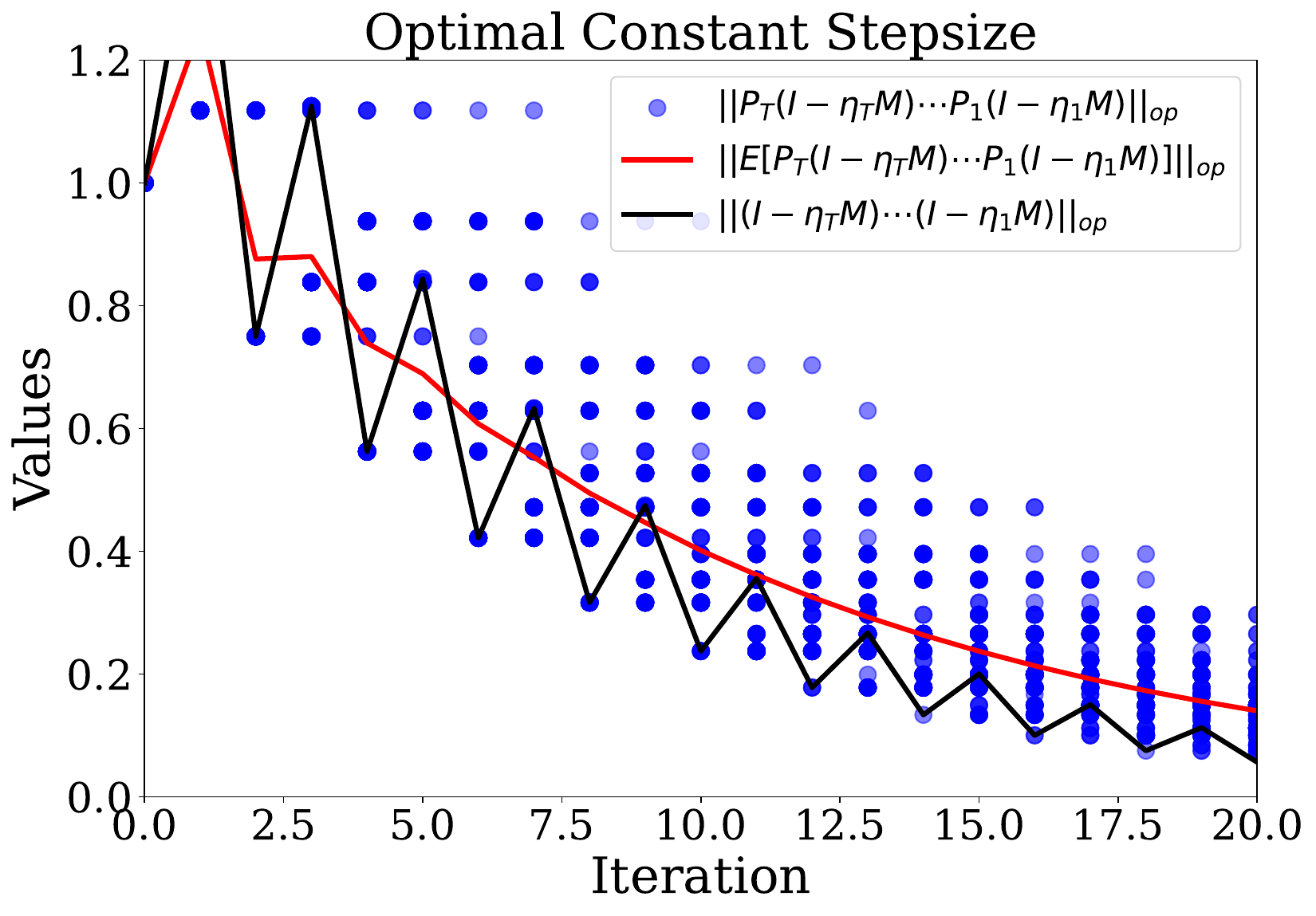}}
   \subfigure[$\pproj$ = 0.2]{\includegraphics[width=0.30\textwidth]{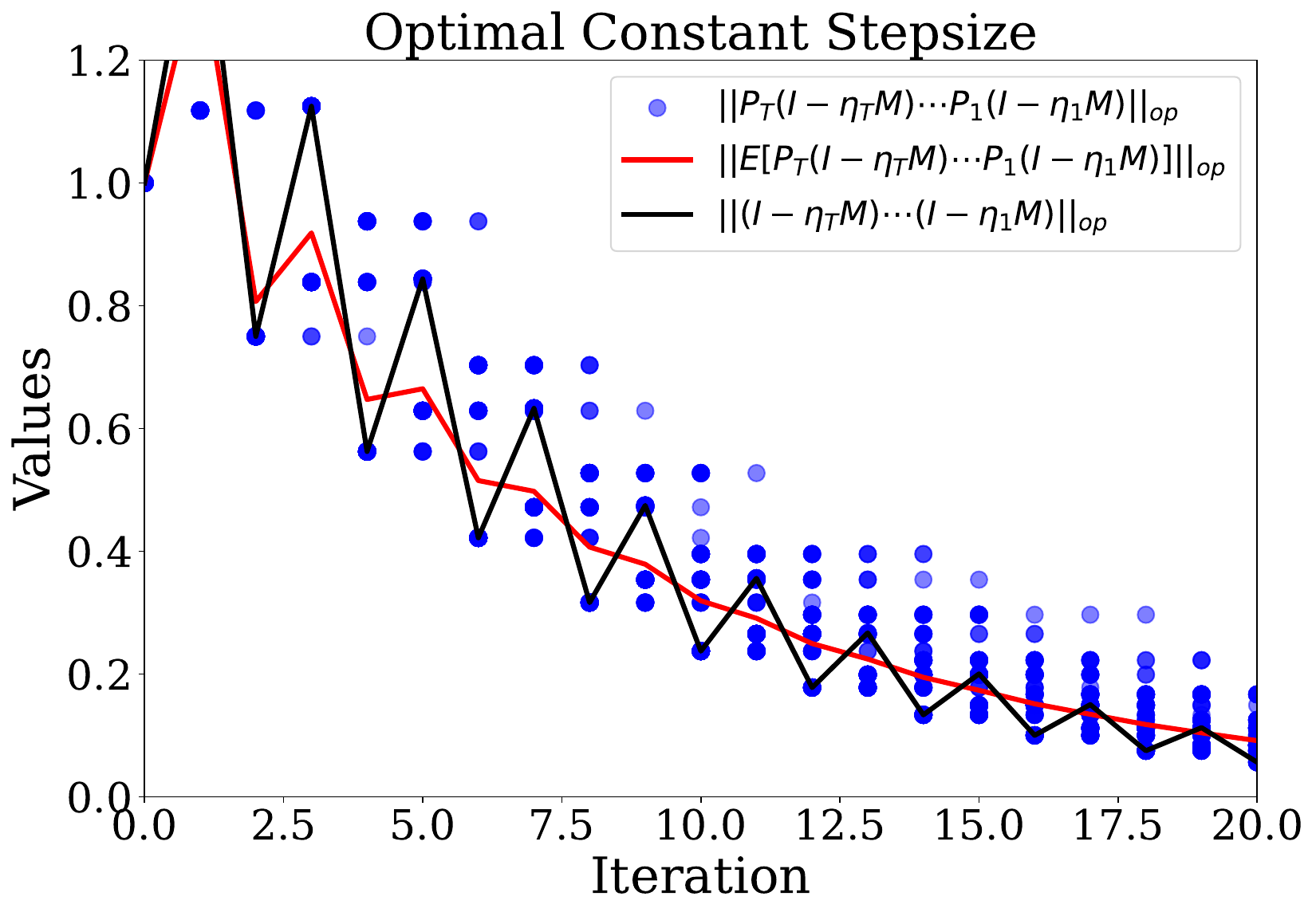}}
  \subfigure[$\pproj$ = 0.05]{\includegraphics[width=0.30\textwidth]{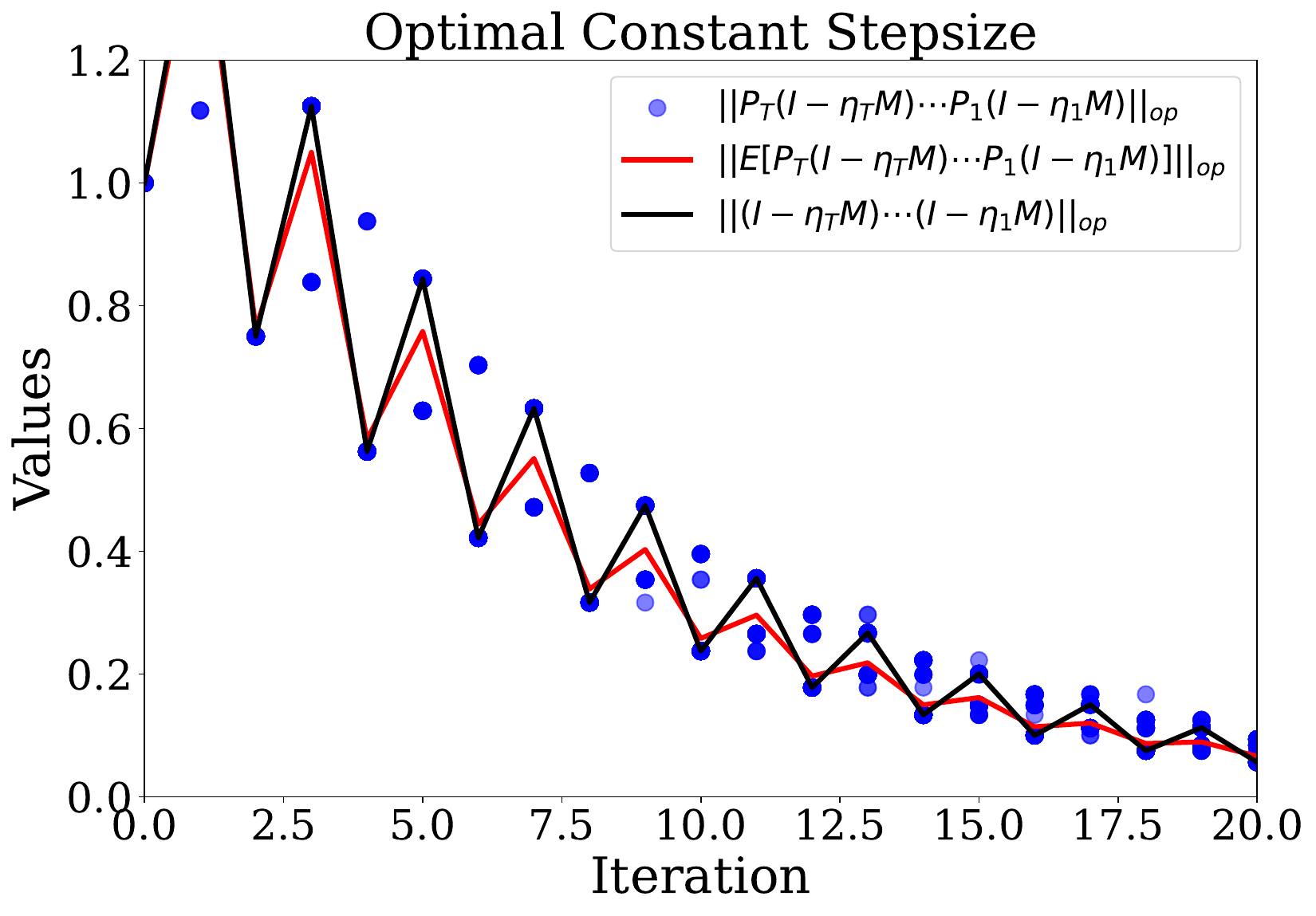}}
  \subfigure[$\pproj$ = 0.4]{\includegraphics[width=0.30\textwidth]{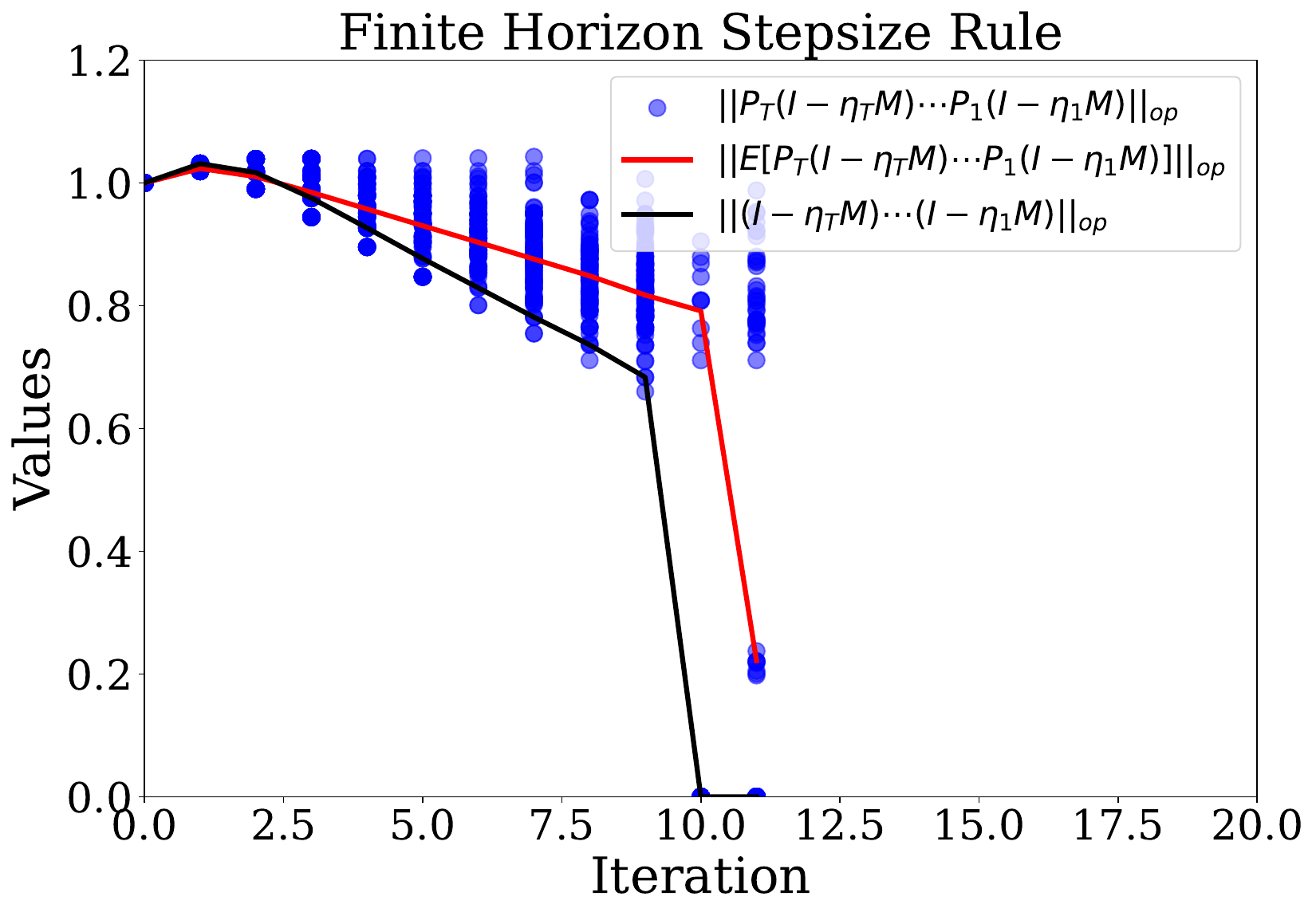}}
   \subfigure[$\pproj$ = 0.2]{\includegraphics[width=0.30\textwidth]{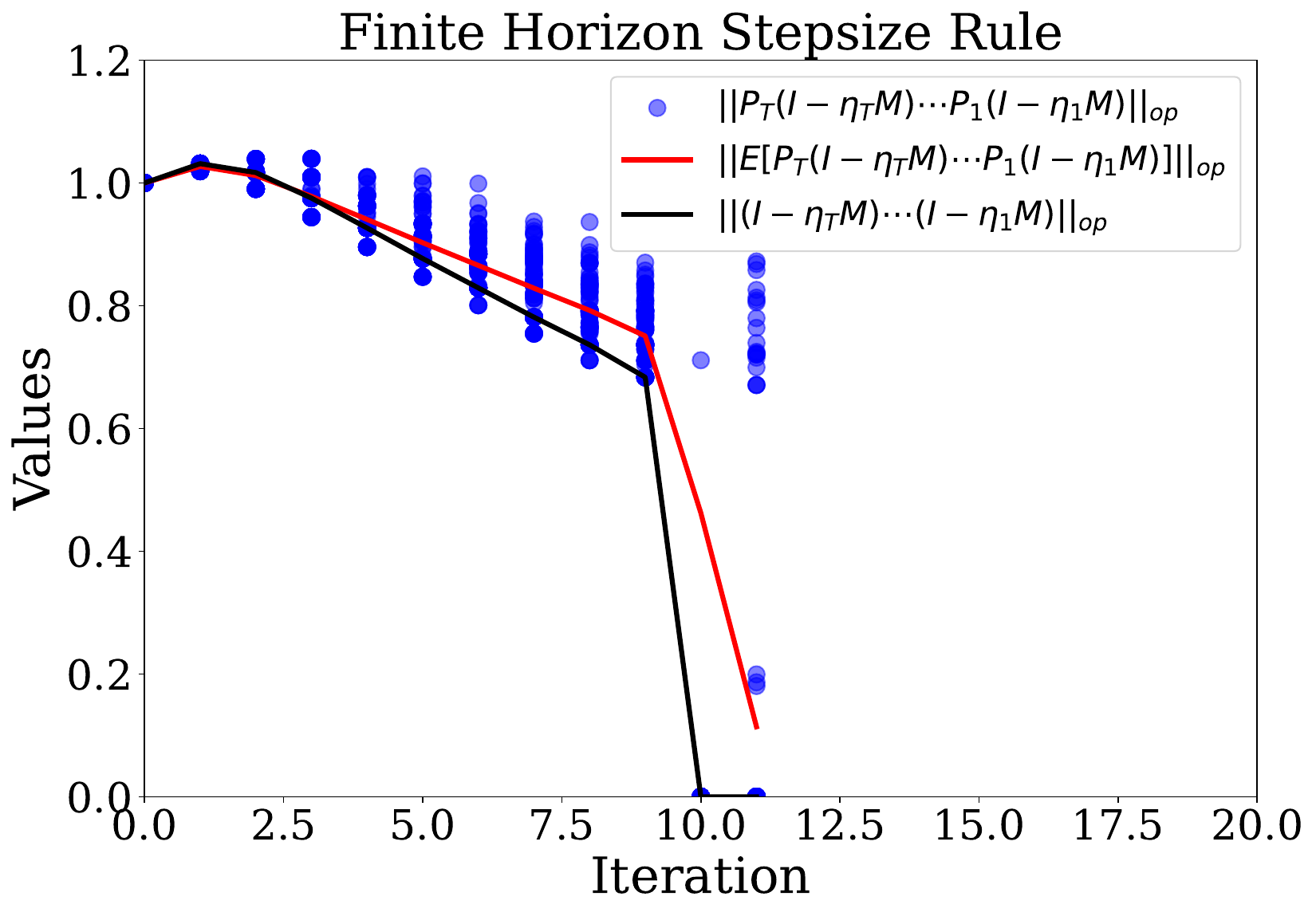}}
  \subfigure[$\pproj$ = 0.05]{\includegraphics[width=0.30\textwidth]{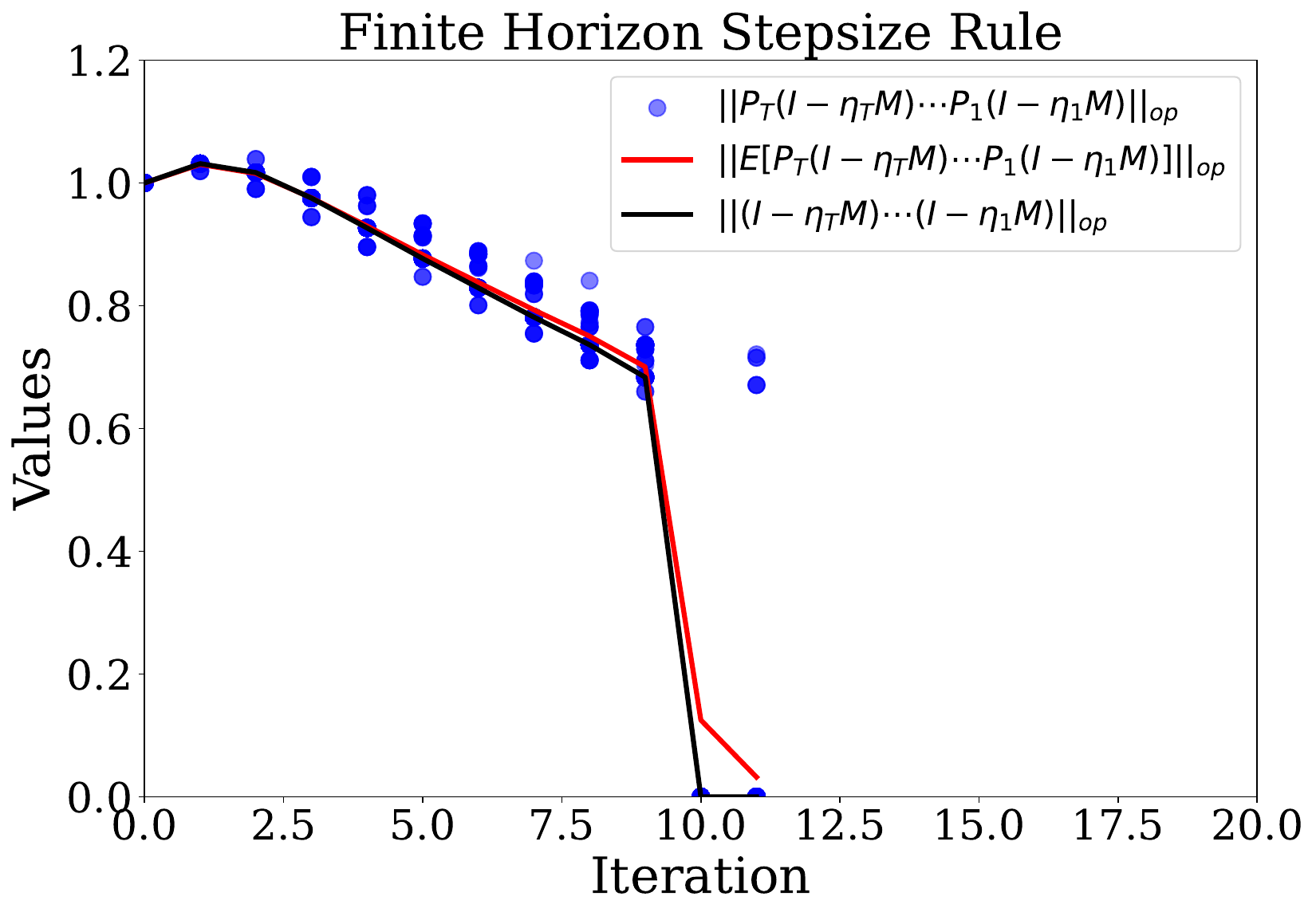}}
  \caption{ The evolution spectral norm of three types of error matrices along iterations.  We use the total iteration budget $T =10$ for Finite Horizon stepsize rule. We find that: the expected error matrix under projection (red curve) is close to that of the error matrix
without projection (black curve), and their similarity grows when the projection probability
$\pproj$ shrinks to 0. Further, Finite Horizon stepsize rule can help the black curve reach almost zero error at the final step ($T = 10$), while the optimal constant stepsize cannot.  }
\label{fig_spectral_norm}
\vspace{-0.3cm}
\end{figure}

\paragraph{Spectral norm investigation.} We further investigate the changes of error matrices along iterations. For Finite Horizon stepsize rule,  we set the total iteration budget to be $T = 10$. For each iteration $t$, we compare  the spectral norm of  three types of error 
matrices: 
\begin{itemize}[topsep=1pt,parsep=1pt,partopsep=1pt, leftmargin=*]
    \item (i) The error matrix with projection: $P_t(I - \eta_t M)  \cdots P_{1}(I - \eta_{1} M)$, where $P_1, \cdots P_t$ are random projection matrices as defined in \eqref{eq_random_projection}. For each iteration $t$, we randomly generate 100 samples of $P_t$.
    \item  (ii) The expected error matrix with projection: $\Ex \left[P_t(I - \BLUE{\eta_t} M)  \cdots P_{1}(I - \eta_{1}M)\right]$. 
    \item (iii) The error matrix without projection: $(I - \eta_t M) \cdots (I - \eta_{1}M)$.
\end{itemize}

The results are shown in Figure \ref{fig_spectral_norm}. We find that the spectral norm of the expected error matrix (red curve) is close to that of the error matrix without projection (black curve). Further, their similarity grows when the projection probability $\pproj$ shrinks to 0. This phenomenon indicates that: if the  $\pproj$ is small, we can effectively control the expected error rate (red curve) by manipulating the spectral norm of the error matrix without projection $\|(I - \eta_t M) \cdots (I - \eta_{1}M)\|_{\operatorname{op}}$ (black curve). As such, the SDP formulation of Finite Horizon stepsize rule \eqref{eq_sdp}, which focuses on optimizing the decay rate of the black curve, can provide good choices on stepsize rules. Indeed, the black curve reaches 0 error at $T = 10$ (see Figure \ref{fig_spectral_norm} (d,e,f)), while optimal constant stepsize cannot make it. In the sequel, we will verify that $\pproj$ is indeed small on real-world LP instances, and Finite Horizon stepsize rule is effective.

\subsection{Results on Real-World LP Instances}
\label{sec_netlib}

Now we verify the effectiveness of Finite Horizon stepsize rule on real-world LP instances. 

\paragraph{Setups.} We consider \texttt{Netlib} LP benchmark \citep{gay1985electronic} more than 90 real-world LP instances \footnote{\url{https://github.com/scipy/scipy/tree/main/benchmarks/benchmarks/linprog_benchmark_files}}. The largest instance has $m = 6084$ constraints and $n= 12243$ variables. We set the stopping criteria as follows: for some easy instances, we set the targeted error as $\epsilon(x,y) =$ 1e-4; for other hard instances, we set the targeted error as $\epsilon(x,y) =$ 1e-1  with an algorithm termination condition after 1e5 iterations. If the constant stepsize method hits the termination condition after 1e5 iterations, we will stop Finite Horizon method when it reaches the final precision of the constant stepsize method.   For Finite Horizon stepsize rule, we set the total iteration budget $T \leq 10$.  For those with $T = 10$, we adjust \texttt{sdp\_iter} (the number of iterations of SDP solver to solve \eqref{eq_sdp}) from 100 to 20 and find that it is sufficient to bring good performance. We did not use larger $T$ since we find that further increasing  $T$ will slow down the SDP solver for solving the sub-problem \eqref{eq_sdp}, and the SDP solver will encounter numerical instability due to the high-degree matrix polynomial constraints. If the algorithm does not reach the desired precision after $T$ steps, we will apply Finite Horizon stepsize strategy in a cyclical manner until the stopping criteria are met.
 We find that {\bf ``cyclically repeating the stepsize rule" is an effective trick} to extend Finite Horizon stepsize rule with small $T$ to the large-$T$ scenarios. 

We compare Finite Horizon stepsize rule to the optimal constant stepsize (Baseline 2), i.e.,  $\eta^* = 2/ (\lambda_1 + \lambda_{2m})$, where $\lambda_1 $ and $\lambda_{2m} $ are expressed in \eqref{eq_lambda_relation_sigma}. Here, we did not consider $T$-step optimal constant stepsize (Baseline 1) since it is unclear how to solve the minimax problem \eqref{eq_toy_example_minimax_T}. 

As mentioned above, both Finite Horizon stepsize rule and the optimal constant stepsize query the same problem information to determine the hyperparameters: they both query the largest \& smallest singular values of $A$, i.e., $\mu$ and $L$, but this information is used differently and results in different stepsize rules.

\begin{figure}[t]
  \vspace{-1.5cm}
    \centering
    \subfigure[Optimality gap curves]{\includegraphics[width=0.30\textwidth]{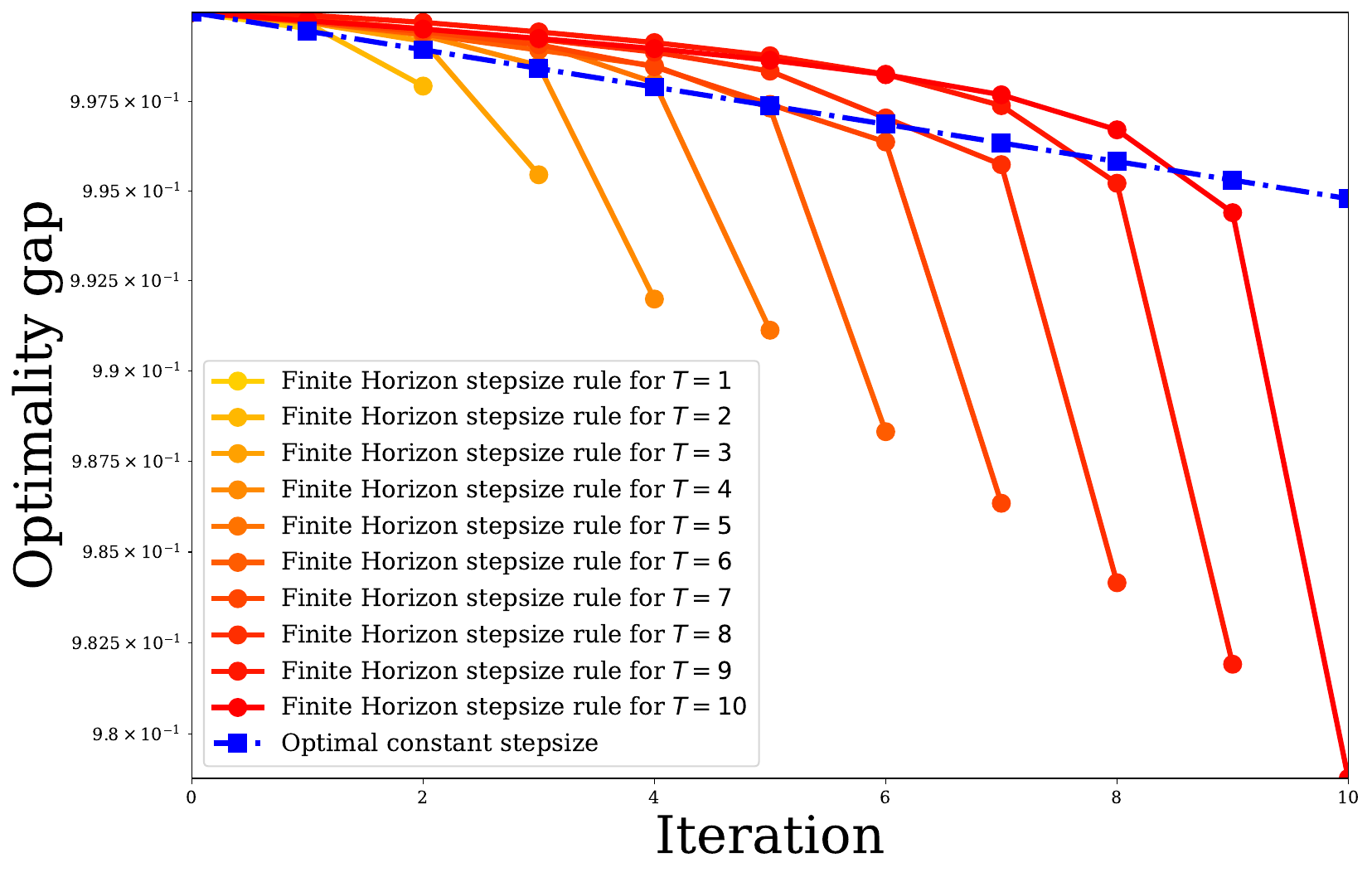}}
    \subfigure[Stepsize rules]{\includegraphics[width=0.28\textwidth]{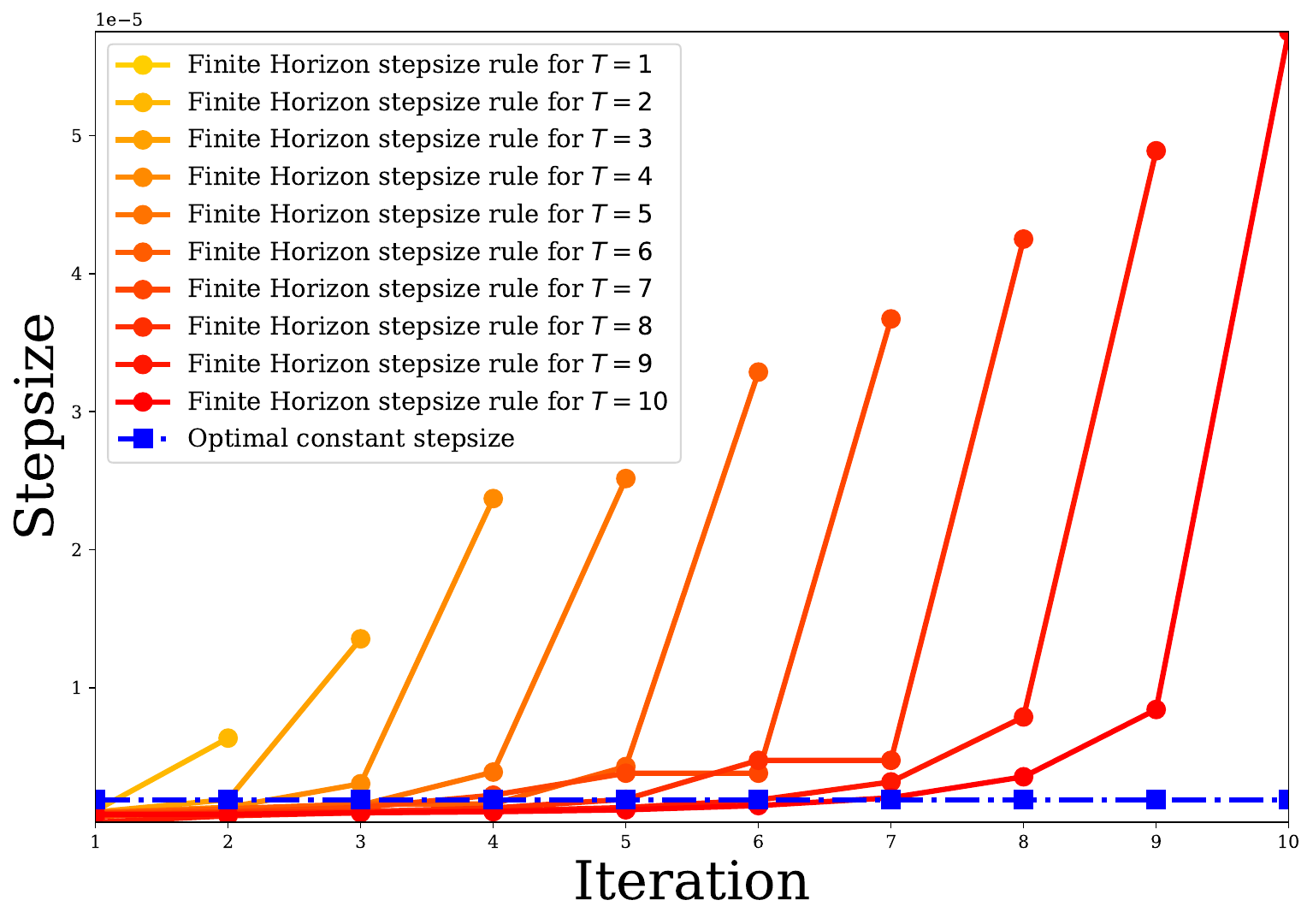}}
    \subfigure[Optimality gap curves]{\includegraphics[width=0.30\textwidth]{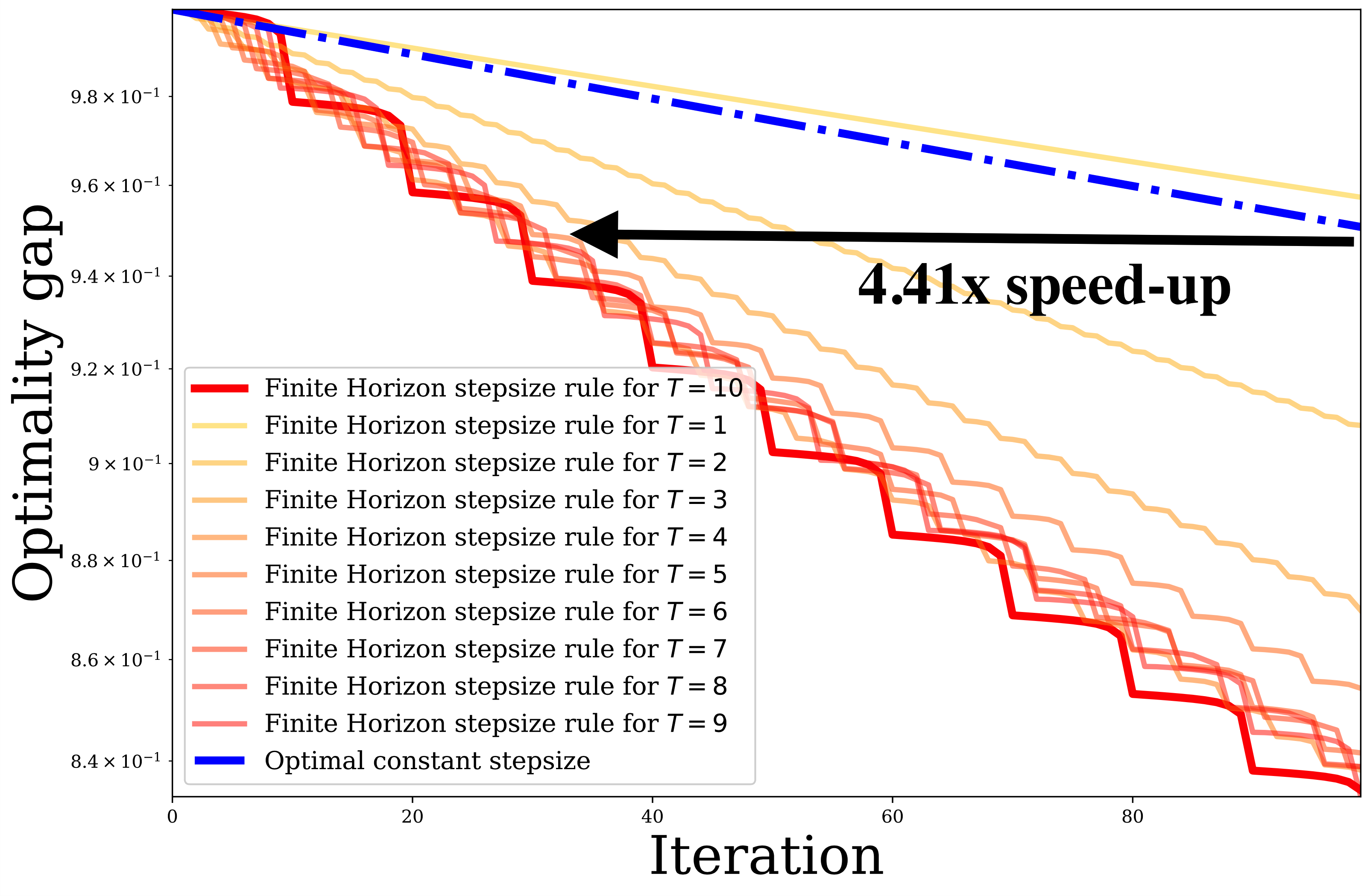}}
    \vspace{-0.2cm}
    \caption{ An example of Finite Horizon stepsize rule on real-world LP instance ( \texttt{Netlib-AGG}).  We find: (a, b) for Finite Horizon stepsize rule, the optimality gap will decrease significantly at the $T$-th step. Meanwhile, the stepsize will increase sharply; (c) Finite Horizon stepsize rule for $T>1$ consistently outperforms optimal constant stepsize. }
  \label{fig_netlib_agg}
  \vspace{-0.3cm}
\end{figure}

\paragraph{Case Study on \texttt{Netlib-AGG} Instance.} We first provide a case study on \texttt{Netlib-AGG} with $m =488 $ $n= 615$. We conduct Finite Horizon stepsize rule for $T = 1, 2, \cdots, 10$ and compare them to the optimal constant stepsize. The results are shown in Figure \ref{fig_netlib_agg}. Similarly to the observations on the toy example (Section \ref{sec_toy_example}), the stepsize will increase sharply at the $T$-th step and will result in a significant decrease in the optimality gap. 

We also observe some different phenomena from the previous toy example. In the previous toy example, the final stepsize was the same for different $T$ (see Figure \ref{fig_netlib_agg} (b)), whereas here, larger $T$ results in larger final stepsize. Accordingly, the optimality gap will decrease to a smaller value. The reason for this difference is that in the previous toy example, the optimal value is already reached at the $T$-th step, whereas it has not been reached here.

\begin{wrapfigure}{r}{0.35\textwidth}
  \includegraphics[width=0.33\textwidth]{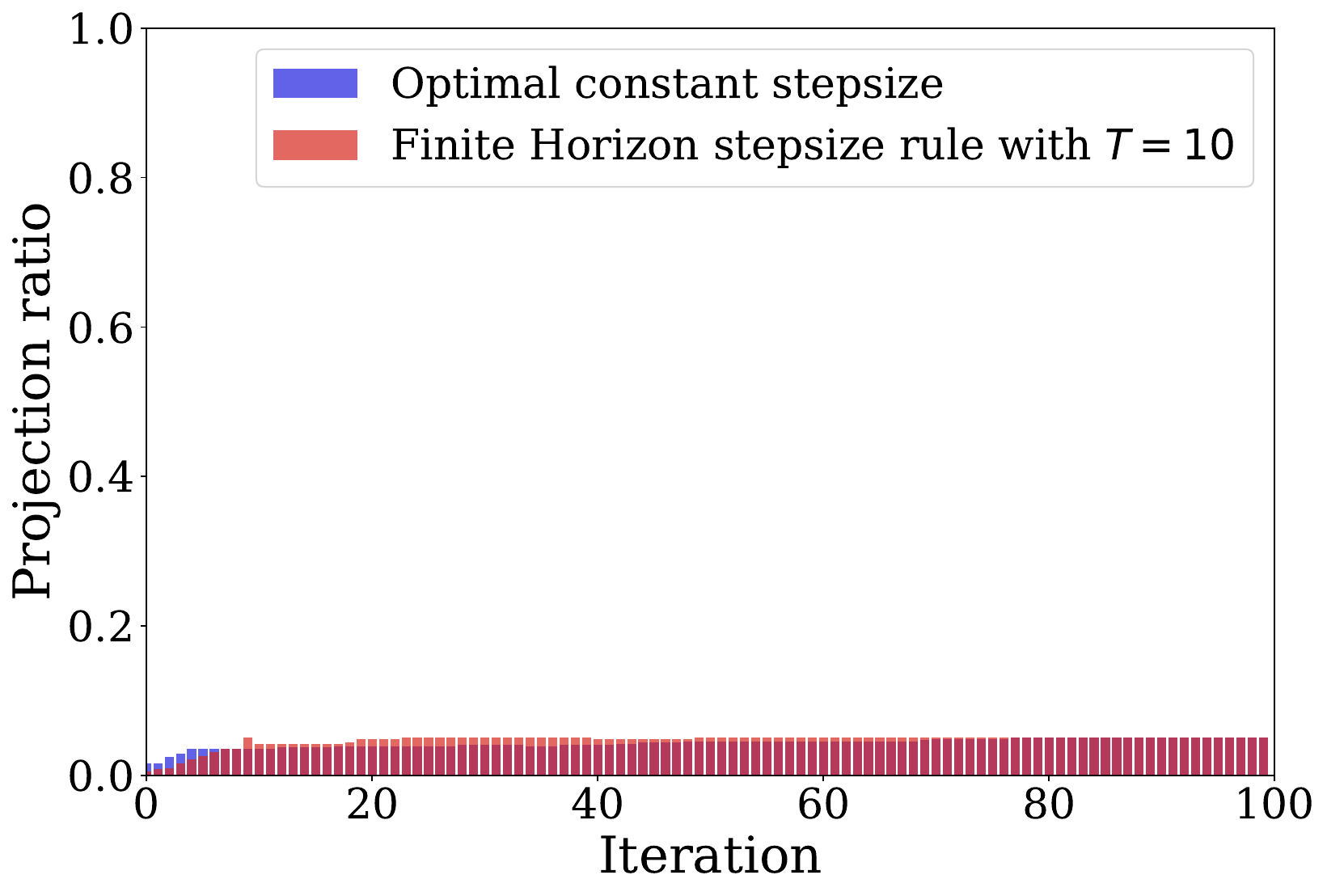}
      \caption{ The projection ratio (i.e., the proportion of entries in $x$ that are forced to 0) along iterations is consistently small ($<5\%$).  }
    \label{fig_projection}
    \vspace{-0.3cm}
  \end{wrapfigure}
In Figure \ref{fig_netlib_agg} (c), we further run the algorithm for 100 iterations by cyclically repeating the stepsize rule.  We find that Finite Horizon stepsize rule for $T > 1$ consistently outperforms the optimal constant stepsize and $T= 10$ performs the best. We observe {\bf 4.41$\times$} speed-up in this case.  For the rest of the instances in  \texttt{Netlib} benchmark, we will choose  $T=10$ and omit the others.

Figure~\ref{fig_projection} shows the projection ratio (i.e., the proportion of entries in $x$ that are forced to 0) along iterations on \texttt{Netlib-AGG}. We find that the ratio is consistently small ($<5\%$) along the trajectories. More specifically,  the average projection ratio is 4.3\% and 4.7\% for the constant and Finite Horizon stepsize rule, respectively. This supports the assumption on small $p_{\text{proj}}$ in Section~\ref{sec_main_methods}.

\begin{figure}[h]
  \centering
  \subfigure[Finite Horizon Stepsize rules]{\includegraphics[width=0.30\textwidth]{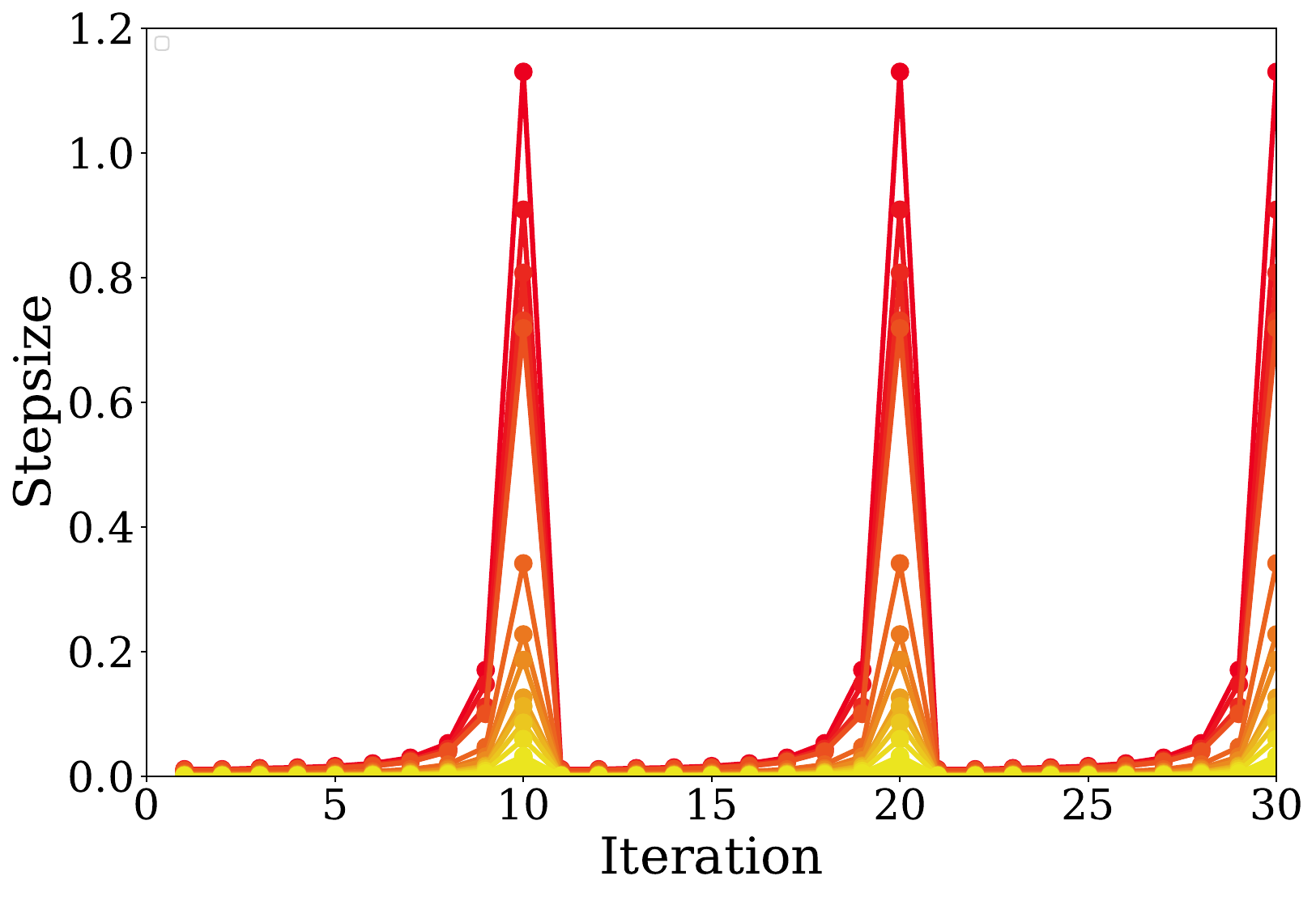}}
  \subfigure[Average iteration complexity]{\includegraphics[width=0.30\textwidth]{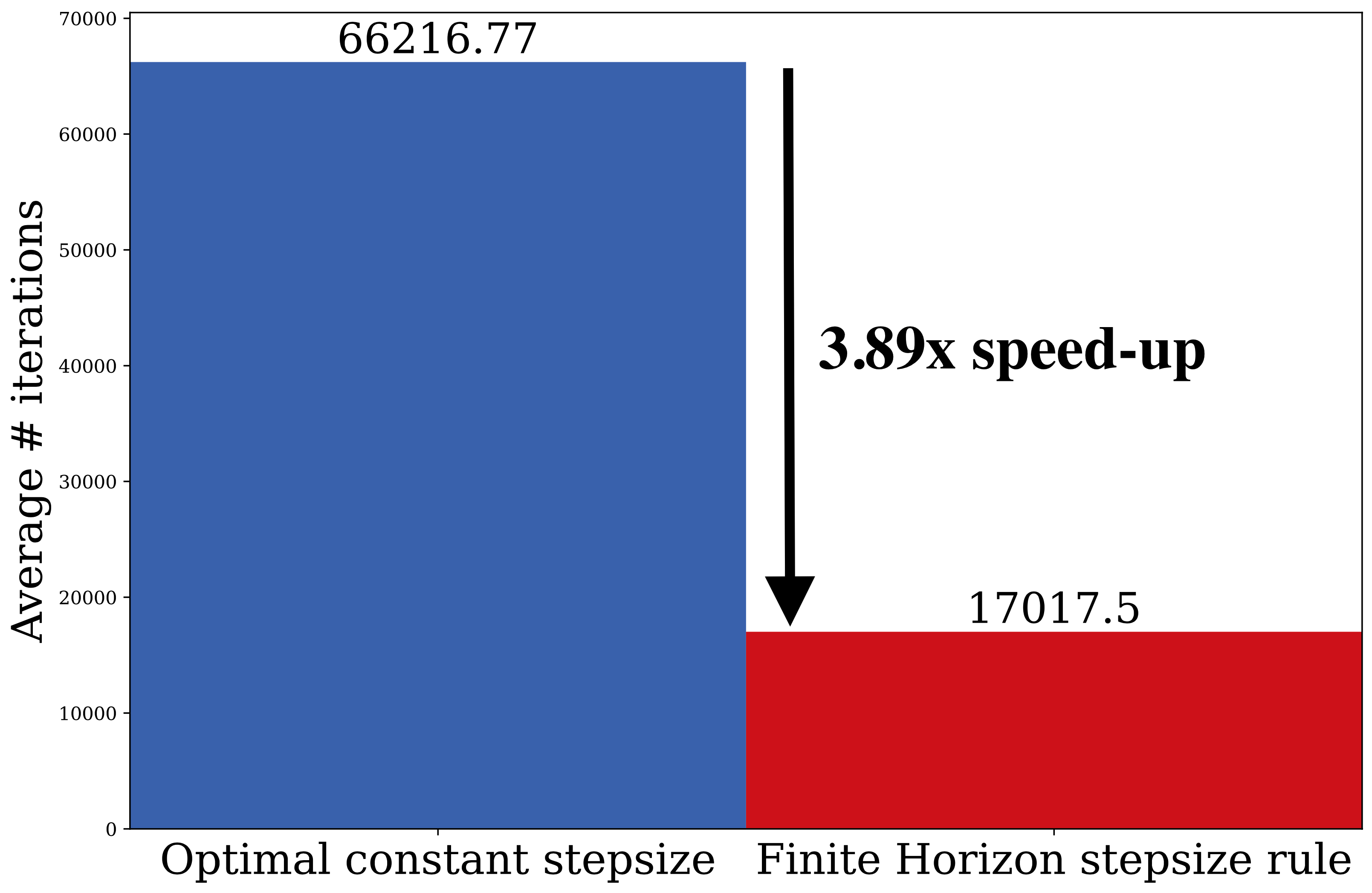}}
  \subfigure[Average running time]{\includegraphics[width=0.30\textwidth]{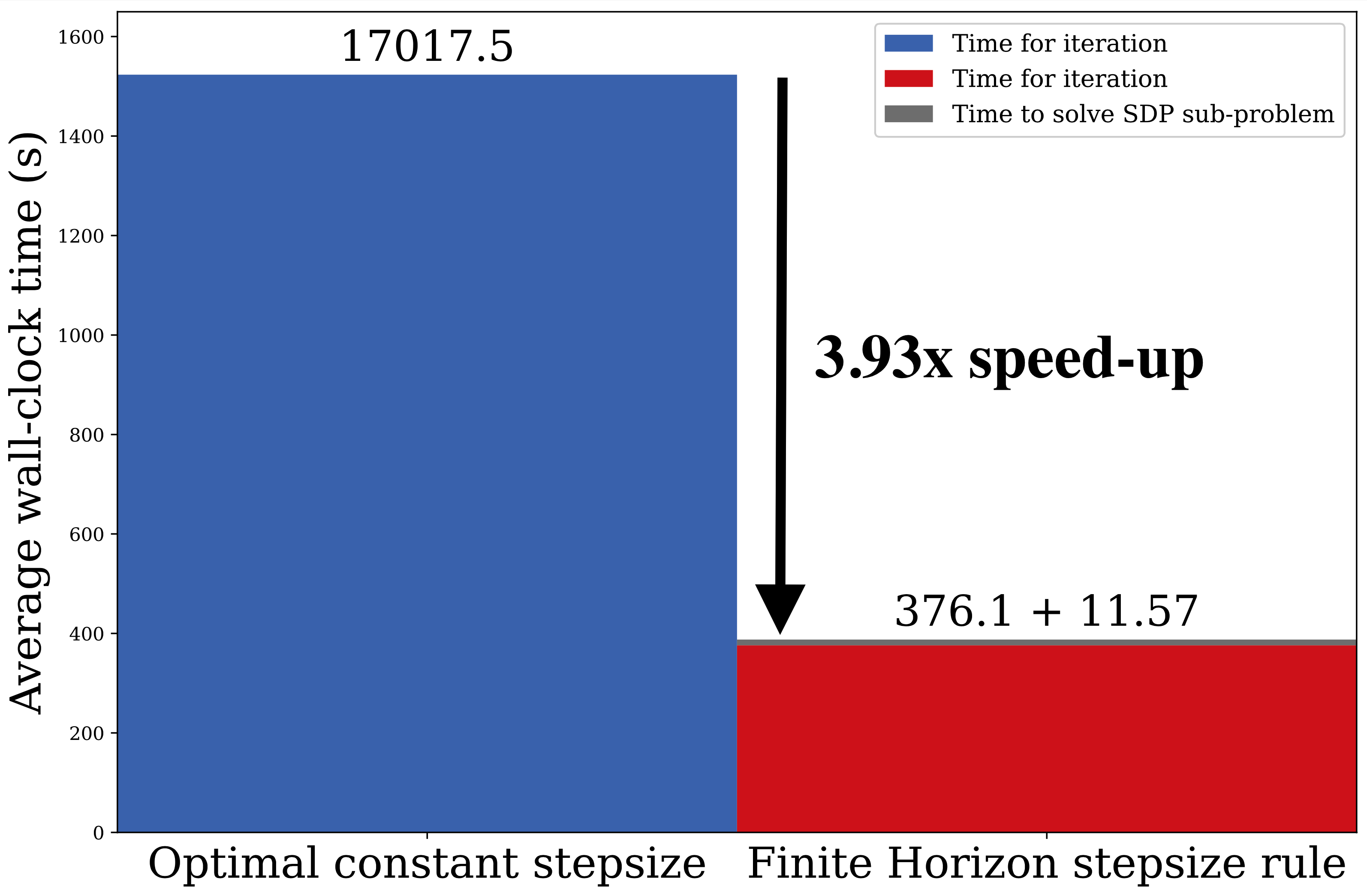}}
  \vspace{-0.2cm}
  \caption{  {\bf (a):} Finite Horizon stepsize rule for $T = 10$ on all the intances in \texttt{Netlib}. We observe that the stepsize at the final iteration ($T=10$) always increases sharply. {\bf (b,c):} Average iteration complexity and running time to  on \texttt{Netlib}. Our stepsize rule achieves on average $3.89 \times$ and $3.93 \times$  speed-up in terms of the number of iterations and wall-clock time, respectively. Further,  the SDP sub-problem takes negligible time to solve (11.57 seconds, $< 3.1\%$ of total time).}  
\label{fig_netlib_complete}
\end{figure}

\paragraph{More Results on \texttt{Netlib} Benchmark.}
We now demonstrate more results on \texttt{Netlib} benchmark. Figure~\ref{fig_netlib_complete} (a) demonstrates Finite Horizon stepsize rule for $T = 10$ on all the instances in \texttt{Netlib} Benchmark. We demonstrate the results up to 30 iterations with 3 repeated cycles of the stepsize rules. We observe a universal pattern across all LP instances: the stepsize at the final iteration ($T=10$) increases sharply, and the magnitude of this increase depends on the specific instance. Interestingly, 
this pattern is similar to the recent stepsize rules for GD, where a large stepsize will appear once in a while among several small ones (e.g., \citep{altschuler2023acceleration}). We leave more investigation on this pattern as a future direction.

Figure~\ref{fig_netlib_complete} (b,c) reports the average speed-up of Finite Horizon stepsize rule over the optimal constant stepsize to hit the same stopping criteria. The results are averaged over all the LP instances from \texttt{Netlib} benchmark. 
In Table~\ref{tab_netlib},  we also present a detailed running time analysis for the top-10 largest instances from \texttt{Netlib}. We summarize the key observations here.
In terms of iteration complexity,  our stepsize rule saves {\bf 74.30\%} iterations over the optimal constant stepsize to reach the same precision, which achieves on average {\bf 3.89$\times$} speed-up.  In terms of running time, our stepsize rule saves  {\bf 75.31\%} wall-clock time, which achieves {\bf 3.93$\times$} speed-up. Further, the SDP sub-problem takes an average 11.57 seconds to solve, which is negligible compared to the overall running time (3.1 \%). Note that the time for solving the SDP sub-problem does not increase with the problem size, since the SDP sub-problem only involves matrices sized $4\times 4$. 

Nevertheless, we acknowledge that the 11.57 seconds might still be too expensive for scenarios with millisecond-level time budgets. We believe there is significant room for improvement in both the design and implementation of this approach, which we leave for future work.

\begin{table}[h]
\vspace{-0.1cm}
\centering
\caption{The wall-clock running time (s) for different methods on the 10 largest datasets in the Netlib Benchmark. Note that the time for solving the SDP sub-problem does not increase with the problem size, since the SDP sub-problem only involves matrices sized $4\times 4$. }
\label{tab_netlib}
\resizebox{0.95\linewidth}{!}{%
\begin{tabular}{lccccccccc}
\toprule[1pt]
\multirow{2}{*}{Dataset} & \multirow{2}{*}{\# constraints} &  \multirow{2}{*}{\# variables} 
 &  \multirow{2}{*}{Time for constant stepsize (s)  }  
 && \multicolumn{2}{c}{Time for Finite Horizon stepsize rule} && \multirow{2}{*}{Total time (speed-up ratio)} 
\\  
 \cmidrule{6-8}
&&& / \# iterations& & Time to solve SDP (s) & Time for iteration (s) / \# iterations  \\
\midrule
DFL001& 6084& 12243 &983.91 /2903 && 7.28 & 231.12  / 661 && 238.40 (\textbf{ 75.77 \% $\downarrow$ })\\
QAP12& 3192& 8856 &4527.73 /12986 && 15.60 & 864.66  / 2199 && 880.26 (\textbf{ 80.56 \% $\downarrow$ })\\
PILOT& 2875& 6294 &17751.92 /100001 && 8.37 & 6314.68  / 33931 && 6323.05 (\textbf{ 64.38 \% $\downarrow$ })\\
SIERRA& 3263& 4771 &21426.60 /22274 && 7.25 & 1860.69  / 33931 && 1867.94 (\textbf{ 91.28 \% $\downarrow$ })\\
D2Q06C& 2171& 5831 &11174.35 /22630 && 10.67 & 5026.00  / 33931 && 5036.66 (\textbf{ 54.93 \% $\downarrow$ })\\
BNL2& 2324& 4486 &7571.56 /54761 && 9.68 & 1354.14  / 12391 && 1363.82 (\textbf{ 81.99 \% $\downarrow$ })\\
WOODW& 1098& 8418 &6853.10 /21540 && 9.36 & 1538.90  / 12391 && 1548.26 (\textbf{ 77.41 \% $\downarrow$ })\\
CYCLE& 1980& 3448 &6265.64 /22637 && 9.23 & 1396.78  / 12391 && 1406.02 (\textbf{ 77.56 \% $\downarrow$ })\\
STOCFOR2& 2157& 3045 &3381.69 /22630 && 10.18 & 712.24  / 12391 && 722.42 (\textbf{ 78.64 \% $\downarrow$ })\\
SHIP12L& 1151& 5533 &3567.07 /22590 && 8.14 & 1153.50  / 12391 && 1161.64 (\textbf{ 67.43 \% $\downarrow$ })\\
\bottomrule
\end{tabular}
 }
\end{table}

\section{Related Works}
\label{sec_related_work}
In addition to all the related works described above, we now mention more literature on the stepsize design of GD, performance estimation of first-order methods, and algorithms for solving LP.

\paragraph{Stepsize design of GD.} 
 There is extensive research on designing stepsize rules for vanilla GD for unconstrained minimization problems. The most simple and classical textbook stepsize choice of GD  is the constant stepsize, which has guarantees on generic smooth functions \citep{nesterov2018lectures}. There are also some classical adaptive strategies including exact line search \citep{boyd2004convex}, Armijo-Goldstein rule \citep{nesterov2018lectures},  Barzilai-Borwein-type stepsize \citep{barzilai1988two}, Polyak-type stepsize \citep{pedregal2004introduction}, etc.. However, these adaptive stepsize strategies either have limited theoretical justification for their benefits, or are already known to have no improvement over constant stepsize in the worst case  \citep{de2017worst, li2021faster}.

 The very first theoretical justification for the time-varying stepsize is from the pioneer work \citep{young1953richardson}. 
 For unconstrained quadratic minimization, \citet{young1953richardson} showed in 1953 that GD can be accelerated using non-constant stepsizes, if the total iteration $T$ is pre-determined (as in the finite horizon setup). The stepsize rule is given explicitly as the inverses of the roots of Chebyshev polynomials, and we will refer to it as Young's stepsize. The results in \citep{young1953richardson}, for the very first time,  revealed that time-varying stepsizes can perform better than constant ones. However, Young's stepsize highly relies on the special connection between quadratic minimization and minimax polynomial theory, and no analysis has yet succeeded in extending it to non-quadratic cases over the past 70 years. Young's stepsize was recently confirmed to be ``not generalizable beyond quadratic" by \citep[Chapter 8]{altschuler2018greed}, where the authors showed that it is provably bad for generic convex settings. 
 For decades, researchers have shifted their attention to other designs, such as 
 adding extra building blocks beyond choosing stepsizes, e.g., via momentum \citep{nesterov1983method}. There has been limited progress on the stepsize design of GD since \citep{young1953richardson}.

A recent influential work \citep{altschuler2018greed} has once again brought people's attention back to stepsize design of GD.
 Perhaps a bit surprisingly, the authors proposed a new stepsize rule of vanilla GD that can bring acceleration beyond quadratic problems. 
 Since then, a series of new stepsize design has emerged.
 These designs have theoretical benefits under different settings. To better put them into our context, we here only highlight their dependency on the iteration horizon $T$. 
 \citep{altschuler2018greed} provided an optimal stepsize for  $T= 2,3$ for strongly convex functions. By cycling through these stepsizes, a constant factor improvement is achieved over the standard unaccelerated rate for constant stepsizes. For the non-strongly-convex setting, 
other stepsize rules of GD are developed under different iteration horizons including $T = 2,3$ \citep{daccache2019performance}, $T = 25$ \citep{das2024branch}, and $T = 127$ \citep{grimmer2024provably}. 
 Recently, 
\citep{altschuler2023acceleration, altschuler2024acceleration} 
proposed the silver stepsize schedules that can achieve acceleration for $T = 2^k -1$ for arbitrary $k\in \mathbb{N}_{+}$. The convergence rate was later refined by a factor of 2 in \citep{wang2024relaxed} under the same requirement of $T$, and has been extended to proximal gradient very recently \citep{bok2024accelerating}.  Concurrent to the silver stepsize schedule, \citep{grimmer2023accelerated,grimmer2024accelerated} also uncovered a new stepsize schedule that brings acceleration at $T = 2^k -1$ for arbitrary $k\in \mathbb{N}_{+}$. Inspired by the aforementioned works, \citet{zhang2024accelerated, grimmer2024composing} further proposed to construct new stepsize schedules by concatenating shorter ones, and \citet{zhang2024accelerated} proved accelerated rate for any pre-determined $T$ (as in the finite horizon setup). \citet{zhang2024anytime} further carefully removed the dependency of $T$ from the stepsize construction and established ``anytime accelerated convergence" over the constant stepsize (albeit at a slower rate than that of the silver stepsize schedule  \citep{altschuler2023acceleration,altschuler2024acceleration}). The results in \citep{zhang2024anytime} answered the open question posted in \citep{kornowski2024open}. There are also new conjectures and insights into the behaviors of constant-stepsize GD \citep{eloi2022worst,grimmer2024strengthened,rotaru2024exact}.

The above works focus on theoretical acceleration. In deep learning field, there are also many new stepsize (learning rate) designs targeted for better empirical performance, e.g., ReduceLROnPlateau \citep{sun2020optimization} \footnote{\url{https://pytorch.org/docs/stable/generated/torch.optim.lr_scheduler.ReduceLROnPlateau.html}.}, warm-up \citep{he2016deep, goyal2017accurate}, cyclic learning rate \citep{smith2017cyclical,fu2019cyclical}, Cosine Annealing \citep{loshchilov2016sgdr}, and many other recent proposals \citep{defazio2023and,ibrahim2024simple,hu2024minicpm,defazio2024road}. Some of these stepsize schedules are theoretically studied in some special function classes such as river valley landscape \citep{wen2024understanding}, and objectives with structured Hessian spectrum (e.g., two-cluster  \citep{oymak2021provable,goujaud2022super} and multiscale \citep{kelner2022big}). The interplay between stepsize choice and logistic loss is also studied in \citep{wu2024large,cai2024large}.

{\bf We now summarize the difference between our work and the works above.} These works primarily focus on GD for unconstrained minimization \footnote{A recent concurrent work \citep{bok2024accelerating} extends the silver stepsize schedule to projected GD. Note that projected GD is not a practical method for LP since the projections onto polyhedral constraints is as hard as solving the original LP \citep{lu2024first}. }, while we study the primal-dual method for LP. One fundamental difference between these two settings is that: for the primal-dual method, its update matrix is non-symmetric (i.e., $M$ is non-symmetric in \eqref{eq_lp_no_inequality}), whereas for GD, its update matrix  (Hessian) is always symmetric. We highlight that non-symmetric update matrix will raise discrepancies between its singular values and eigenvalues, which do {\it not} exist in the symmetric case. 
We find that such discrepancies  will substantially changes the principle for the time-varying stepsize design, especially in the {\it non-asymptotic} case when $T$ is bounded. More detailed discussions on the technical challenges can be seen in Section \ref{sec_main_methods}.

\paragraph{Performance estimation problem (PEP) for first-order methods.} Parallel to the aforementioned works on GD, there is an active line of works analyzing the worst-case performance of a given first-order method (including but not limited to GD) for unconstrained minimization. 
\citet{drori2014performance} proposed to analyze the worst-case performance of a given first-order algorithm by numerically solving an additional optimization problem, which they called the Performance Estimation Problem (PEP). The authors formulated PEP as a nonconvex quadratic matrix program \citep{beck2007quadratic}, which was then relaxed and solved via its dual SDP problem. \citet{taylor2017smooth}  later showed that the SDP reformulation is actually {\it exact} by introducing the notion of convex interpolation. To further improve the efficiency of solving PEP, researchers also developed various frameworks such as the  Integral quadratic constraints (IQC) \citep{lessard2016analysis} and nonconvex QCQP (especially for analyzing nonlinear conjugate gradient methods) \citep{gupta2023nonlinear}.  The development of PEP has led to many exciting advances in the theory of first-order methods, including tightening the classical bounds (e.g., \citep{de2017worst,fazlyab2018analysis, de2020worst,barre2020complexity, drori2022oracle,teboulle2023elementary} and all the works above), the construction of new Lyapunov functions \citep{taylor2018lyapunov,taylor2019stochastic},   performance estimation of a {\it given} monotone operator splitting method with constant ``stepsize" \citep{ryu2020operator,gu2020tight,wang2019robust}, and new analysis of  ADMM \citep{nishihara2015general,seidman2019control}.  PEP can also be extended to tighten the complexity bound of projected or proximal gradient methods \citep{lessard2016analysis, taylor2017exact, taylor2018exact, kim2018another, d2021acceleration, teboulle2023elementary}, mirror descent method \citep{dragomir2022optimal}, and very recently, Frank-Wolfe method \citep{luner2024performance}. It is worth mentioning that PEP requires a pre-determined $T$, so it falls into the finite horizon setup. {\bf One key difference between these works and ours is that:} they focus on estimating the performance of a {\it given} algorithm, while we focus on designing new stepsize rules.

\paragraph{New algorithm design via PEP.} 
In addition to analyzing a given algorithm, PEP also inspired many new algorithm designs. There are two mainstream approaches.

The first approach is to combine PEP with human expertise to inspire new proof ideas.
This process often requires professional researchers to conduct SDP-based analysis, investigate the structure or the solution of PEP,
and then manually generate a rigorous convergence proof. Following this approach, researchers have designed new stepsize design of GD \citep{grimmer2023accelerated,grimmer2024accelerated}, new momentum mechanism \citep{kim2016optimized,van2017fastest, cyrus2018robust,yoon2021accelerated,taylor2023optimal},  new splitting method \citep{ryu2020finding}, and new method for composite nonsmooth saddle-point problem  \citep{drori2014contributions}. All these methods have asymptotic guarantees based on non-trivial design and proof. 
Their development requires human experts to conduct case-by-case exploration.

The second approach is to automatically search for a better design by introducing an ``outer minimization problem" to PEP. While the PEP itself can be modeled as a convex problem (usually a SDP) \citep{drori2014performance, taylor2017smooth},  the overall minimax search problem is usually highly nonconvex \citep{drori2014performance,das2024branch}. \citep[Section 8]{altschuler2018greed} used algebraic techniques to solve the search problem analytically for $T \leq 3$. However, as commented by the authors, it seems difficult to extend these computations for a larger $T$. 
\citet[Section 5]{drori2014performance} proposed to relax the nonconvex search problem to a convex SDP, and then they discovered an optimal stepsize rule for $T \leq 5$ (their Figure 4). Using similar convex relaxation or heuristics, researchers have discovered many new designs such as new momentum mechanism \citep{kim2017convergence,drori2017exact,drori2020efficient,kim2018generalizing,kim2021optimizing,zhou2022practical,park2024optimal}, new proximal point methods  \citep{kim2021accelerated,lieder2021convergence,park2022exact,barre2023principled},
and new cutting-plane methods \citep{drori2016optimal}.
Recently,  
\citet{das2024branch} reformulated the overall minimax problem as a nonconvex QCQP  proposed to solve it via customized spatial branch-and-bound algorithms. They then discovered the optimal stepsize schedule of GD up to $T = 25$.  A similar stepsize schedule of GD for $T \leq  127$  are also discovered by  ``brute force searching" \citep{grimmer2023accelerated}.  The methodologies in \citep{das2024branch} have also been extended to unconstrained nonsmooth problems to design 
new ISTA-type methods \citep{jang2023computer}.

{\bf Our work is substantially different from the aforementioned works in three aspects.}  {\bf First,} PEP framework is primarily designed for unconstrained minimization (and some closely related variants), while we focus on LP. It is widely recognized that linear programming (LP) exhibits many distinct structures and its algorithms are often developed in a tailored or specialized manner \citep{luenberger1984linear}. 
We study LP because: in engineering scenarios with strict time constraints (e.g., wireless communications, power systems, and autonomous vehicles),  LP-related problems frequently emerge (see Section \ref{sec_intro}). This observation motivates us to re-design LP algorithms under a finite-horizon framework.
{\bf Second}, our stepsize design does not stem from PEP. Rather, it is an independently developed method. 
As a result, we have significantly different methodologies from PEP:  PEP requires solving a {\it nonconvex} large-scale QCQP to design the optimal stepsize of GD, whereas we only require a {\it convex} small-scaled SDP with $4\times 4$ matrices. The convexity of our search problem stems from LP. {\bf Third}, PEP frameworks primarilly focus on theoretical benefits, while we focus on gaining speed-up on real-world LP benchmarks (in addition to theoretical guarantee).

\paragraph{Algorithms for LP.} Since the 1940s, extensive research has been devoted to developing algorithms for LP. Today, the de-facto and mostly recognized methods include Dantzig’s simplex method \citep{dantzig2016linear,dantzig1990origins} and interior-point methods (IPMs) \citep{karmarkar1984new,renegar1988polynomial,monteiro1989interior,nesterov1994interior,wright1997primal,ye2011interior}. 
Recently,  first-order methods is becoming increasingly popular for solving large-scale LPs (e.g.  \citep{chambolle2011first,goldstein2015adaptive,applegate2021practical,lin2021admm,deng2022new,basu2020eclipse,o2016conic,o2021operator,applegate2023faster,lu2023cupdlp,lu2023unified,xiong2023computational,xiong2024role,xiong2024accessible,deng2024enhanced,fercoq2024monitoring}. A more comprehensive review can be seen in \citep{lu2024first}). Unlike simplex and IPMs, which require matrix factorization, first-order methods rely solely on matrix-vector multiplications. This makes them highly scalable on GPU and distributed computing platforms \citep{lu2024first}. To our knowledge, we are the first work that re-design LP algorithms under the notion of ``finite horizon". We find that:  if we abandon the asymptotic performance and instead focus on how the algorithm performs within a finite number of iterations, the primal-dual method can reach 3$\times$ to  4$\times$
acceleration by merely changing its stepsize design. We believe similar potential is inherently embedded in a wide range of LP algorithms, and our results might provide new perspectives to further accelerate the current methods.

\section{Conclusion}
In this work, we propose the framework of finite horizon optimization, which focus on algorithm design for the scenario where the total iteration budget is fixed and finite. We revisit the primal-dual method for LP and propose a new stepsize rule called Finite Horizon stepsize rule. We prove that this stepsize rule can reach an accelerated convergence rate at the $T$-th step, where $T$ is the pre-determined total iteration budget. 
 On a wide range of real-world LP instances, Finite Horizon stepsize rule reaches substantial speed-up over the optimal constant stepsize.
 One important future direction is to design Finite Horizon stepsize rules for more advanced algorithms, such as the momentum or the preconditioned variants of primal-dual method \citep{applegate2021practical}.  It is also intriguing to apply the finite horizon framework to broader practical scenarios beyond LP. 

\section*{Acknowledgement}
The authors would like to thank Prof. Ruoyu Sun, Jiawei Zhang, Tian Ding, Prof. Tongxin Li, Wenzhi Gao, Xin Jiang, Ziniu Li, Xueyao Zhang, and group members of Prof. Mingyi Hong for the valuable discussions.

\bibliographystyle{abbrvnat}
\bibliography{reference}

\begin{thebibliography}{158}
\providecommand{\natexlab}[1]{#1}
\providecommand{\url}[1]{\texttt{#1}}
\expandafter\ifx\csname urlstyle\endcsname\relax
  \providecommand{\doi}[1]{doi: #1}\else
  \providecommand{\doi}{doi: \begingroup \urlstyle{rm}\Url}\fi

\bibitem[Adler and {\"O}ktem(2018)]{adler2018learned}
J.~Adler and O.~{\"O}ktem.
\newblock Learned primal-dual reconstruction.
\newblock \emph{IEEE transactions on medical imaging}, 37\penalty0 (6):\penalty0 1322--1332, 2018.

\bibitem[Altschuler(2018)]{altschuler2018greed}
J.~Altschuler.
\newblock \emph{Greed, hedging, and acceleration in convex optimization}.
\newblock PhD thesis, Massachusetts Institute of Technology, 2018.

\bibitem[Altschuler and Parrilo(2023)]{altschuler2023acceleration}
J.~M. Altschuler and P.~A. Parrilo.
\newblock Acceleration by stepsize hedging i: Multi-step descent and the silver stepsize schedule.
\newblock \emph{arXiv preprint arXiv:2309.07879}, 2023.

\bibitem[Altschuler and Parrilo(2024)]{altschuler2024acceleration}
J.~M. Altschuler and P.~A. Parrilo.
\newblock Acceleration by stepsize hedging: Silver stepsize schedule for smooth convex optimization.
\newblock \emph{Mathematical Programming}, pages 1--14, 2024.

\bibitem[Applegate et~al.(2021)Applegate, D{\'\i}az, Hinder, Lu, Lubin, O'Donoghue, and Schudy]{applegate2021practical}
D.~Applegate, M.~D{\'\i}az, O.~Hinder, H.~Lu, M.~Lubin, B.~O'Donoghue, and W.~Schudy.
\newblock Practical large-scale linear programming using primal-dual hybrid gradient.
\newblock \emph{Advances in Neural Information Processing Systems}, 34:\penalty0 20243--20257, 2021.

\bibitem[Applegate et~al.(2023)Applegate, Hinder, Lu, and Lubin]{applegate2023faster}
D.~Applegate, O.~Hinder, H.~Lu, and M.~Lubin.
\newblock Faster first-order primal-dual methods for linear programming using restarts and sharpness.
\newblock \emph{Mathematical Programming}, 201\penalty0 (1):\penalty0 133--184, 2023.

\bibitem[Babaeinejadsarookolaee et~al.(2019)Babaeinejadsarookolaee, Birchfield, Christie, Coffrin, DeMarco, Diao, Ferris, Fliscounakis, Greene, Huang, et~al.]{babaeinejadsarookolaee2019power}
S.~Babaeinejadsarookolaee, A.~Birchfield, R.~D. Christie, C.~Coffrin, C.~DeMarco, R.~Diao, M.~Ferris, S.~Fliscounakis, S.~Greene, R.~Huang, et~al.
\newblock The power grid library for benchmarking ac optimal power flow algorithms.
\newblock \emph{arXiv preprint arXiv:1908.02788}, 2019.

\bibitem[Barr{\'e} et~al.(2020)Barr{\'e}, Taylor, and d’Aspremont]{barre2020complexity}
M.~Barr{\'e}, A.~Taylor, and A.~d’Aspremont.
\newblock Complexity guarantees for polyak steps with momentum.
\newblock In \emph{Conference on learning theory}, pages 452--478. PMLR, 2020.

\bibitem[Barr{\'e} et~al.(2023)Barr{\'e}, Taylor, and Bach]{barre2023principled}
M.~Barr{\'e}, A.~B. Taylor, and F.~Bach.
\newblock Principled analyses and design of first-order methods with inexact proximal operators.
\newblock \emph{Mathematical Programming}, 201\penalty0 (1):\penalty0 185--230, 2023.

\bibitem[Barzilai and Borwein(1988)]{barzilai1988two}
J.~Barzilai and J.~M. Borwein.
\newblock Two-point step size gradient methods.
\newblock \emph{IMA journal of numerical analysis}, 8\penalty0 (1):\penalty0 141--148, 1988.

\bibitem[Basu et~al.(2020)Basu, Ghoting, Mazumder, and Pan]{basu2020eclipse}
K.~Basu, A.~Ghoting, R.~Mazumder, and Y.~Pan.
\newblock Eclipse: An extreme-scale linear program solver for web-applications.
\newblock In \emph{International Conference on Machine Learning}, pages 704--714. PMLR, 2020.

\bibitem[Beck(2007)]{beck2007quadratic}
A.~Beck.
\newblock Quadratic matrix programming.
\newblock \emph{SIAM Journal on Optimization}, 17\penalty0 (4):\penalty0 1224--1238, 2007.

\bibitem[Boher et~al.(2008)Boher, Rabineau, and Helard]{boher2008fpga}
L.~Boher, R.~Rabineau, and M.~Helard.
\newblock Fpga implementation of an iterative receiver for mimo-ofdm systems.
\newblock \emph{IEEE Journal on Selected Areas in Communications}, 26\penalty0 (6):\penalty0 857--866, 2008.

\bibitem[Bok and Altschuler(2024)]{bok2024accelerating}
J.~Bok and J.~M. Altschuler.
\newblock Accelerating proximal gradient descent via silver stepsizes.
\newblock \emph{arXiv preprint arXiv:2412.05497}, 2024.

\bibitem[Boyd and Vandenberghe(2004)]{boyd2004convex}
S.~Boyd and L.~Vandenberghe.
\newblock \emph{Convex optimization}.
\newblock Cambridge university press, 2004.

\bibitem[Cai et~al.(2024)Cai, Wu, Mei, Lindsey, and Bartlett]{cai2024large}
Y.~Cai, J.~Wu, S.~Mei, M.~Lindsey, and P.~L. Bartlett.
\newblock Large stepsize gradient descent for non-homogeneous two-layer networks: Margin improvement and fast optimization.
\newblock \emph{arXiv preprint arXiv:2406.08654}, 2024.

\bibitem[Cauchy et~al.(1847)]{cauchy1847methode}
A.~Cauchy et~al.
\newblock M{\'e}thode g{\'e}n{\'e}rale pour la r{\'e}solution des systemes d’{\'e}quations simultan{\'e}es.
\newblock \emph{Comp. Rend. Sci. Paris}, 25\penalty0 (1847):\penalty0 536--538, 1847.

\bibitem[Chambolle and Pock(2011)]{chambolle2011first}
A.~Chambolle and T.~Pock.
\newblock A first-order primal-dual algorithm for convex problems with applications to imaging.
\newblock \emph{Journal of mathematical imaging and vision}, 40:\penalty0 120--145, 2011.

\bibitem[Chasparis and Shamma(2005)]{chasparis2005linear}
G.~C. Chasparis and J.~S. Shamma.
\newblock Linear-programming-based multi-vehicle path planning with adversaries.
\newblock In \emph{Proceedings of the 2005, American Control Conference, 2005.}, pages 1072--1077. IEEE, 2005.

\bibitem[Chatzivasileiadis(2018)]{chatzivasileiadis2018optimal}
S.~Chatzivasileiadis.
\newblock Optimal power flow (dc-opf and ac-opf).
\newblock \emph{DTU Summer School}, 2018.

\bibitem[Chebyshev(1853)]{chebyshev1853theorie}
P.~L. Chebyshev.
\newblock \emph{Th{\'e}orie des m{\'e}canismes connus sous le nom de parall{\'e}logrammes}.
\newblock Imprimerie de l'Acad{\'e}mie imp{\'e}riale des sciences, 1853.

\bibitem[Chen et~al.(2022)Chen, Chen, Chen, Heaton, Liu, Wang, and Yin]{chen2022learning}
T.~Chen, X.~Chen, W.~Chen, H.~Heaton, J.~Liu, Z.~Wang, and W.~Yin.
\newblock Learning to optimize: A primer and a benchmark.
\newblock \emph{Journal of Machine Learning Research}, 23\penalty0 (189):\penalty0 1--59, 2022.

\bibitem[Chen et~al.(2012)Chen, Han, and Zhao]{chen2012three}
Y.~Chen, J.~Han, and X.~Zhao.
\newblock Three-dimensional path planning for unmanned aerial vehicle based on linear programming.
\newblock \emph{Robotica}, 30\penalty0 (5):\penalty0 773--781, 2012.

\bibitem[Christie et~al.(2000)Christie, Wollenberg, and Wangensteen]{christie2000transmission}
R.~D. Christie, B.~F. Wollenberg, and I.~Wangensteen.
\newblock Transmission management in the deregulated environment.
\newblock \emph{Proceedings of the IEEE}, 88\penalty0 (2):\penalty0 170--195, 2000.

\bibitem[Culligan(2006)]{culligan2006online}
K.~F. Culligan.
\newblock \emph{Online trajectory planning for UAVs using mixed integer linear programming}.
\newblock PhD thesis, Massachusetts Institute of Technology, 2006.

\bibitem[Cyrus et~al.(2018)Cyrus, Hu, Van~Scoy, and Lessard]{cyrus2018robust}
S.~Cyrus, B.~Hu, B.~Van~Scoy, and L.~Lessard.
\newblock A robust accelerated optimization algorithm for strongly convex functions.
\newblock In \emph{2018 Annual American Control Conference (ACC)}, pages 1376--1381. IEEE, 2018.

\bibitem[Daccache et~al.(2019)Daccache, Glineur, and Hendrickx]{daccache2019performance}
A.~Daccache, F.~Glineur, and J.~Hendrickx.
\newblock \emph{Performance estimation of the gradient method with fixed arbitrary step sizes}.
\newblock PhD thesis, Master’s thesis, Universit{\'e} Catholique de Louvain, 2019.

\bibitem[Dantzig(1990)]{dantzig1990origins}
G.~B. Dantzig.
\newblock Origins of the simplex method.
\newblock In \emph{A history of scientific computing}, pages 141--151. 1990.

\bibitem[Dantzig(2002)]{dantzig2002linear}
G.~B. Dantzig.
\newblock Linear programming.
\newblock \emph{Operations research}, 50\penalty0 (1):\penalty0 42--47, 2002.

\bibitem[Dantzig(2016)]{dantzig2016linear}
G.~B. Dantzig.
\newblock Linear programming and extensions.
\newblock In \emph{Linear programming and extensions}. Princeton university press, 2016.

\bibitem[Das~Gupta et~al.(2024)Das~Gupta, Van~Parys, and Ryu]{das2024branch}
S.~Das~Gupta, B.~P. Van~Parys, and E.~K. Ryu.
\newblock Branch-and-bound performance estimation programming: A unified methodology for constructing optimal optimization methods.
\newblock \emph{Mathematical Programming}, 204\penalty0 (1):\penalty0 567--639, 2024.

\bibitem[De~Klerk et~al.(2017)De~Klerk, Glineur, and Taylor]{de2017worst}
E.~De~Klerk, F.~Glineur, and A.~B. Taylor.
\newblock On the worst-case complexity of the gradient method with exact line search for smooth strongly convex functions.
\newblock \emph{Optimization Letters}, 11:\penalty0 1185--1199, 2017.

\bibitem[De~Klerk et~al.(2020)De~Klerk, Glineur, and Taylor]{de2020worst}
E.~De~Klerk, F.~Glineur, and A.~B. Taylor.
\newblock Worst-case convergence analysis of inexact gradient and newton methods through semidefinite programming performance estimation.
\newblock \emph{SIAM Journal on Optimization}, 30\penalty0 (3):\penalty0 2053--2082, 2020.

\bibitem[Defazio et~al.(2023)Defazio, Cutkosky, Mehta, and Mishchenko]{defazio2023and}
A.~Defazio, A.~Cutkosky, H.~Mehta, and K.~Mishchenko.
\newblock When, why and how much? adaptive learning rate scheduling by refinement.
\newblock \emph{arXiv preprint arXiv:2310.07831}, 2023.

\bibitem[Defazio et~al.(2024)Defazio, Yang, Mehta, Mishchenko, Khaled, and Cutkosky]{defazio2024road}
A.~Defazio, X.~A. Yang, H.~Mehta, K.~Mishchenko, A.~Khaled, and A.~Cutkosky.
\newblock The road less scheduled.
\newblock \emph{arXiv preprint arXiv:2405.15682}, 2024.

\bibitem[Deng et~al.(2022)Deng, Feng, Gao, Ge, Jiang, Jiang, Liu, Liu, Xue, Ye, et~al.]{deng2022new}
Q.~Deng, Q.~Feng, W.~Gao, D.~Ge, B.~Jiang, Y.~Jiang, J.~Liu, T.~Liu, C.~Xue, Y.~Ye, et~al.
\newblock New developments of admm-based interior point methods for linear programming and conic programming.
\newblock \emph{arXiv preprint arXiv:2209.01793}, 2022.

\bibitem[Deng et~al.(2024)Deng, Feng, Gao, Ge, Jiang, Jiang, Liu, Liu, Xue, Ye, et~al.]{deng2024enhanced}
Q.~Deng, Q.~Feng, W.~Gao, D.~Ge, B.~Jiang, Y.~Jiang, J.~Liu, T.~Liu, C.~Xue, Y.~Ye, et~al.
\newblock An enhanced alternating direction method of multipliers-based interior point method for linear and conic optimization.
\newblock \emph{INFORMS Journal on Computing}, 2024.

\bibitem[Dolgov et~al.(2008)Dolgov, Thrun, Montemerlo, and Diebel]{dolgov2008practical}
D.~Dolgov, S.~Thrun, M.~Montemerlo, and J.~Diebel.
\newblock Practical search techniques in path planning for autonomous driving.
\newblock \emph{Ann Arbor}, 1001\penalty0 (48105):\penalty0 18--80, 2008.

\bibitem[Dragomir et~al.(2022)Dragomir, Taylor, d’Aspremont, and Bolte]{dragomir2022optimal}
R.-A. Dragomir, A.~B. Taylor, A.~d’Aspremont, and J.~Bolte.
\newblock Optimal complexity and certification of bregman first-order methods.
\newblock \emph{Mathematical Programming}, pages 1--43, 2022.

\bibitem[Drori(2017)]{drori2017exact}
Y.~Drori.
\newblock The exact information-based complexity of smooth convex minimization.
\newblock \emph{Journal of Complexity}, 39:\penalty0 1--16, 2017.

\bibitem[Drori and Drori(2014)]{drori2014contributions}
Y.~Drori and Y.~Drori.
\newblock \emph{Contributions to the complexity analysis of optimization algorithms}.
\newblock Universitat Tel-Aviv Tel Aviv, Israel, 2014.

\bibitem[Drori and Taylor(2022)]{drori2022oracle}
Y.~Drori and A.~Taylor.
\newblock On the oracle complexity of smooth strongly convex minimization.
\newblock \emph{Journal of Complexity}, 68:\penalty0 101590, 2022.

\bibitem[Drori and Taylor(2020)]{drori2020efficient}
Y.~Drori and A.~B. Taylor.
\newblock Efficient first-order methods for convex minimization: a constructive approach.
\newblock \emph{Mathematical Programming}, 184\penalty0 (1):\penalty0 183--220, 2020.

\bibitem[Drori and Teboulle(2014)]{drori2014performance}
Y.~Drori and M.~Teboulle.
\newblock Performance of first-order methods for smooth convex minimization: a novel approach.
\newblock \emph{Mathematical Programming}, 145\penalty0 (1):\penalty0 451--482, 2014.

\bibitem[Drori and Teboulle(2016)]{drori2016optimal}
Y.~Drori and M.~Teboulle.
\newblock An optimal variant of kelley’s cutting-plane method.
\newblock \emph{Mathematical Programming}, 160:\penalty0 321--351, 2016.

\bibitem[d’Aspremont et~al.(2021)d’Aspremont, Scieur, Taylor, et~al.]{d2021acceleration}
A.~d’Aspremont, D.~Scieur, A.~Taylor, et~al.
\newblock Acceleration methods.
\newblock \emph{Foundations and Trends{\textregistered} in Optimization}, 5\penalty0 (1-2):\penalty0 1--245, 2021.

\bibitem[Eloi and Glineur(2022)]{eloi2022worst}
D.~Eloi and F.~Glineur.
\newblock \emph{Worst-case functions for the gradient method with fixed variable step sizes}.
\newblock PhD thesis, Master’s thesis, Universit{\'e} Catholique de Louvain, 2022.

\bibitem[Fazlyab et~al.(2018)Fazlyab, Ribeiro, Morari, and Preciado]{fazlyab2018analysis}
M.~Fazlyab, A.~Ribeiro, M.~Morari, and V.~M. Preciado.
\newblock Analysis of optimization algorithms via integral quadratic constraints: Nonstrongly convex problems.
\newblock \emph{SIAM Journal on Optimization}, 28\penalty0 (3):\penalty0 2654--2689, 2018.

\bibitem[Fercoq(2024)]{fercoq2024monitoring}
O.~Fercoq.
\newblock Monitoring the convergence speed of pdhg to find better primal and dual step sizes.
\newblock \emph{arXiv preprint arXiv:2403.19202}, 2024.

\bibitem[Ferraz et~al.(2021)Ferraz, Subramaniyan, Chinthala, Andrade, Cavallaro, Nandy, Silva, Zhang, Purnaprajna, and Falcao]{ferraz2021survey}
O.~Ferraz, S.~Subramaniyan, R.~Chinthala, J.~Andrade, J.~R. Cavallaro, S.~K. Nandy, V.~Silva, X.~Zhang, M.~Purnaprajna, and G.~Falcao.
\newblock A survey on high-throughput non-binary ldpc decoders: Asic, fpga, and gpu architectures.
\newblock \emph{IEEE Communications Surveys \& Tutorials}, 24\penalty0 (1):\penalty0 524--556, 2021.

\bibitem[Frank et~al.(2012)Frank, Steponavice, and Rebennack]{frank2012optimal}
S.~Frank, I.~Steponavice, and S.~Rebennack.
\newblock Optimal power flow: A bibliographic survey i: Formulations and deterministic methods.
\newblock \emph{Energy systems}, 3:\penalty0 221--258, 2012.

\bibitem[Fu et~al.(2019)Fu, Li, Liu, Gao, Celikyilmaz, and Carin]{fu2019cyclical}
H.~Fu, C.~Li, X.~Liu, J.~Gao, A.~Celikyilmaz, and L.~Carin.
\newblock Cyclical annealing schedule: A simple approach to mitigating kl vanishing.
\newblock \emph{arXiv preprint arXiv:1903.10145}, 2019.

\bibitem[Gay(1985)]{gay1985electronic}
D.~M. Gay.
\newblock Electronic mail distribution of linear programming test problems.
\newblock \emph{Mathematical Programming Society COAL Newsletter}, 13:\penalty0 10--12, 1985.

\bibitem[Goldstein et~al.(2015)Goldstein, Li, and Yuan]{goldstein2015adaptive}
T.~Goldstein, M.~Li, and X.~Yuan.
\newblock Adaptive primal-dual splitting methods for statistical learning and image processing.
\newblock \emph{Advances in neural information processing systems}, 28, 2015.

\bibitem[Goujaud et~al.(2022)Goujaud, Scieur, Dieuleveut, Taylor, and Pedregosa]{goujaud2022super}
B.~Goujaud, D.~Scieur, A.~Dieuleveut, A.~B. Taylor, and F.~Pedregosa.
\newblock Super-acceleration with cyclical step-sizes.
\newblock In \emph{International Conference on Artificial Intelligence and Statistics}, pages 3028--3065. PMLR, 2022.

\bibitem[Goyal et~al.(2017)Goyal, Doll{\'a}r, Girshick, Noordhuis, Wesolowski, Kyrola, Tulloch, Jia, and He]{goyal2017accurate}
P.~Goyal, P.~Doll{\'a}r, R.~Girshick, P.~Noordhuis, L.~Wesolowski, A.~Kyrola, A.~Tulloch, Y.~Jia, and K.~He.
\newblock Accurate, large minibatch sgd: Training imagenet in 1 hour. arxiv 2017.
\newblock \emph{arXiv preprint arXiv:1706.02677}, 2017.

\bibitem[Gregor and LeCun(2010)]{gregor2010learning}
K.~Gregor and Y.~LeCun.
\newblock Learning fast approximations of sparse coding.
\newblock In \emph{Proceedings of the 27th international conference on international conference on machine learning}, pages 399--406, 2010.

\bibitem[Grimmer(2024)]{grimmer2024provably}
B.~Grimmer.
\newblock Provably faster gradient descent via long steps.
\newblock \emph{SIAM Journal on Optimization}, 34\penalty0 (3):\penalty0 2588--2608, 2024.

\bibitem[Grimmer et~al.(2023)Grimmer, Shu, and Wang]{grimmer2023accelerated}
B.~Grimmer, K.~Shu, and A.~L. Wang.
\newblock Accelerated gradient descent via long steps.
\newblock \emph{arXiv preprint arXiv:2309.09961}, 2023.

\bibitem[Grimmer et~al.(2024{\natexlab{a}})Grimmer, Shu, and Wang]{grimmer2024accelerated}
B.~Grimmer, K.~Shu, and A.~L. Wang.
\newblock Accelerated objective gap and gradient norm convergence for gradient descent via long steps.
\newblock \emph{arXiv preprint arXiv:2403.14045}, 2024{\natexlab{a}}.

\bibitem[Grimmer et~al.(2024{\natexlab{b}})Grimmer, Shu, and Wang]{grimmer2024composing}
B.~Grimmer, K.~Shu, and A.~L. Wang.
\newblock Composing optimized stepsize schedules for gradient descent.
\newblock \emph{arXiv preprint arXiv:2410.16249}, 2024{\natexlab{b}}.

\bibitem[Grimmer et~al.(2024{\natexlab{c}})Grimmer, Shu, and Wang]{grimmer2024strengthened}
B.~Grimmer, K.~Shu, and A.~L. Wang.
\newblock A strengthened conjecture on the minimax optimal constant stepsize for gradient descent.
\newblock \emph{arXiv preprint arXiv:2407.11739}, 2024{\natexlab{c}}.

\bibitem[Gu and Yang(2020)]{gu2020tight}
G.~Gu and J.~Yang.
\newblock Tight sublinear convergence rate of the proximal point algorithm for maximal monotone inclusion problems.
\newblock \emph{SIAM Journal on Optimization}, 30\penalty0 (3):\penalty0 1905--1921, 2020.

\bibitem[Gupta et~al.(2023)Gupta, Freund, Sun, and Taylor]{gupta2023nonlinear}
S.~D. Gupta, R.~M. Freund, X.~A. Sun, and A.~Taylor.
\newblock Nonlinear conjugate gradient methods: worst-case convergence rates via computer-assisted analyses.
\newblock \emph{arXiv preprint arXiv:2301.01530}, 2023.

\bibitem[He et~al.(2016)He, Zhang, Ren, and Sun]{he2016deep}
K.~He, X.~Zhang, S.~Ren, and J.~Sun.
\newblock Deep residual learning for image recognition.
\newblock In \emph{Proceedings of the IEEE conference on computer vision and pattern recognition}, pages 770--778, 2016.

\bibitem[Horn and Johnson(2012)]{horn2012matrix}
R.~A. Horn and C.~R. Johnson.
\newblock \emph{Matrix analysis}.
\newblock Cambridge university press, 2012.

\bibitem[Hu et~al.(2024)Hu, Tu, Han, He, Cui, Long, Zheng, Fang, Huang, Zhao, et~al.]{hu2024minicpm}
S.~Hu, Y.~Tu, X.~Han, C.~He, G.~Cui, X.~Long, Z.~Zheng, Y.~Fang, Y.~Huang, W.~Zhao, et~al.
\newblock Minicpm: Unveiling the potential of small language models with scalable training strategies.
\newblock \emph{arXiv preprint arXiv:2404.06395}, 2024.

\bibitem[Ibrahim et~al.(2024)Ibrahim, Th{\'e}rien, Gupta, Richter, Anthony, Lesort, Belilovsky, and Rish]{ibrahim2024simple}
A.~Ibrahim, B.~Th{\'e}rien, K.~Gupta, M.~L. Richter, Q.~Anthony, T.~Lesort, E.~Belilovsky, and I.~Rish.
\newblock Simple and scalable strategies to continually pre-train large language models.
\newblock \emph{arXiv preprint arXiv:2403.08763}, 2024.

\bibitem[Jang et~al.(2023)Jang, Gupta, and Ryu]{jang2023computer}
U.~Jang, S.~D. Gupta, and E.~K. Ryu.
\newblock Computer-assisted design of accelerated composite optimization methods: Optista.
\newblock \emph{arXiv preprint arXiv:2305.15704}, 2023.

\bibitem[Kantorovich(1960)]{kantorovich1960mathematical}
L.~V. Kantorovich.
\newblock Mathematical methods of organizing and planning production.
\newblock \emph{Management science}, 6\penalty0 (4):\penalty0 366--422, 1960.

\bibitem[Karmarkar(1984)]{karmarkar1984new}
N.~Karmarkar.
\newblock A new polynomial-time algorithm for linear programming.
\newblock In \emph{Proceedings of the sixteenth annual ACM symposium on Theory of computing}, pages 302--311, 1984.

\bibitem[Kelner et~al.(2022)Kelner, Marsden, Sharan, Sidford, Valiant, and Yuan]{kelner2022big}
J.~Kelner, A.~Marsden, V.~Sharan, A.~Sidford, G.~Valiant, and H.~Yuan.
\newblock Big-step-little-step: Efficient gradient methods for objectives with multiple scales.
\newblock In \emph{Conference on Learning Theory}, pages 2431--2540. PMLR, 2022.

\bibitem[Kiessling et~al.(2022)Kiessling, Zanelli, Nurkanovi{\'c}, Gillis, Diehl, Zeilinger, Pipeleers, and Swevers]{kiessling2022feasible}
D.~Kiessling, A.~Zanelli, A.~Nurkanovi{\'c}, J.~Gillis, M.~Diehl, M.~Zeilinger, G.~Pipeleers, and J.~Swevers.
\newblock A feasible sequential linear programming algorithm with application to time-optimal path planning problems.
\newblock In \emph{2022 IEEE 61st Conference on Decision and Control (CDC)}, pages 1196--1203. IEEE, 2022.

\bibitem[Kim(2021)]{kim2021accelerated}
D.~Kim.
\newblock Accelerated proximal point method for maximally monotone operators.
\newblock \emph{Mathematical Programming}, 190\penalty0 (1):\penalty0 57--87, 2021.

\bibitem[Kim and Fessler(2016)]{kim2016optimized}
D.~Kim and J.~A. Fessler.
\newblock Optimized first-order methods for smooth convex minimization.
\newblock \emph{Mathematical programming}, 159:\penalty0 81--107, 2016.

\bibitem[Kim and Fessler(2017)]{kim2017convergence}
D.~Kim and J.~A. Fessler.
\newblock On the convergence analysis of the optimized gradient method.
\newblock \emph{Journal of optimization theory and applications}, 172\penalty0 (1):\penalty0 187--205, 2017.

\bibitem[Kim and Fessler(2018{\natexlab{a}})]{kim2018another}
D.~Kim and J.~A. Fessler.
\newblock Another look at the fast iterative shrinkage/thresholding algorithm (fista).
\newblock \emph{SIAM Journal on Optimization}, 28\penalty0 (1):\penalty0 223--250, 2018{\natexlab{a}}.

\bibitem[Kim and Fessler(2018{\natexlab{b}})]{kim2018generalizing}
D.~Kim and J.~A. Fessler.
\newblock Generalizing the optimized gradient method for smooth convex minimization.
\newblock \emph{SIAM Journal on Optimization}, 28\penalty0 (2):\penalty0 1920--1950, 2018{\natexlab{b}}.

\bibitem[Kim and Fessler(2021)]{kim2021optimizing}
D.~Kim and J.~A. Fessler.
\newblock Optimizing the efficiency of first-order methods for decreasing the gradient of smooth convex functions.
\newblock \emph{Journal of optimization theory and applications}, 188\penalty0 (1):\penalty0 192--219, 2021.

\bibitem[Kittaneh(2006)]{kittaneh2006spectral}
F.~Kittaneh.
\newblock Spectral radius inequalities for hilbert space operators.
\newblock \emph{Proceedings of the American Mathematical Society}, pages 385--390, 2006.

\bibitem[Kornowski and Shamir(2024)]{kornowski2024open}
G.~Kornowski and O.~Shamir.
\newblock Open problem: Anytime convergence rate of gradient descent.
\newblock \emph{arXiv preprint arXiv:2406.13888}, 2024.

\bibitem[Latva-Aho et~al.(2019)Latva-Aho, Lepp{\"a}nen, et~al.]{latva2019key}
M.~Latva-Aho, K.~Lepp{\"a}nen, et~al.
\newblock Key drivers and research challenges for 6g ubiquitous wireless intelligence.
\newblock 2019.

\bibitem[Lessard et~al.(2016)Lessard, Recht, and Packard]{lessard2016analysis}
L.~Lessard, B.~Recht, and A.~Packard.
\newblock Analysis and design of optimization algorithms via integral quadratic constraints.
\newblock \emph{SIAM Journal on Optimization}, 26\penalty0 (1):\penalty0 57--95, 2016.

\bibitem[Li et~al.(2024)Li, Yang, Chen, Wang, Chen, Mao, Ma, Wang, Ding, Tang, et~al.]{li2024pdhg}
B.~Li, L.~Yang, Y.~Chen, S.~Wang, Q.~Chen, H.~Mao, Y.~Ma, A.~Wang, T.~Ding, J.~Tang, et~al.
\newblock Pdhg-unrolled learning-to-optimize method for large-scale linear programming.
\newblock \emph{arXiv preprint arXiv:2406.01908}, 2024.

\bibitem[Li and Sun(2021)]{li2021faster}
D.~Li and R.~Sun.
\newblock On a faster $ r $-linear convergence rate of the barzilai-borwein method.
\newblock \emph{arXiv preprint arXiv:2101.00205}, 2021.

\bibitem[Li et~al.(2023)Li, Zhang, Guo, Lenzo, and Guo]{li2023real}
G.~Li, X.~Zhang, H.~Guo, B.~Lenzo, and N.~Guo.
\newblock Real-time optimal trajectory planning for autonomous driving with collision avoidance using convex optimization.
\newblock \emph{Automotive Innovation}, 6\penalty0 (3):\penalty0 481--491, 2023.

\bibitem[Lieder(2021)]{lieder2021convergence}
F.~Lieder.
\newblock On the convergence rate of the halpern-iteration.
\newblock \emph{Optimization letters}, 15\penalty0 (2):\penalty0 405--418, 2021.

\bibitem[Lin et~al.(2021)Lin, Ma, Ye, and Zhang]{lin2021admm}
T.~Lin, S.~Ma, Y.~Ye, and S.~Zhang.
\newblock An admm-based interior-point method for large-scale linear programming.
\newblock \emph{Optimization Methods and Software}, 36\penalty0 (2-3):\penalty0 389--424, 2021.

\bibitem[Loshchilov and Hutter(2016)]{loshchilov2016sgdr}
I.~Loshchilov and F.~Hutter.
\newblock Sgdr: Stochastic gradient descent with warm restarts.
\newblock \emph{arXiv preprint arXiv:1608.03983}, 2016.

\bibitem[Lu(2024)]{lu2024first}
H.~Lu.
\newblock First-order methods for linear programming.
\newblock \emph{arXiv preprint arXiv:2403.14535}, 2024.

\bibitem[Lu and Yang(2023{\natexlab{a}})]{lu2023cupdlp}
H.~Lu and J.~Yang.
\newblock cupdlp. jl: A gpu implementation of restarted primal-dual hybrid gradient for linear programming in julia.
\newblock \emph{arXiv preprint arXiv:2311.12180}, 2023{\natexlab{a}}.

\bibitem[Lu and Yang(2023{\natexlab{b}})]{lu2023unified}
H.~Lu and J.~Yang.
\newblock On a unified and simplified proof for the ergodic convergence rates of ppm, pdhg and admm.
\newblock \emph{arXiv preprint arXiv:2305.02165}, 2023{\natexlab{b}}.

\bibitem[Luenberger and Ye(1984)]{luenberger1984linear}
D.~G. Luenberger and Y.~Ye.
\newblock \emph{Linear and nonlinear programming}, volume~2.
\newblock Springer, 1984.

\bibitem[Luner and Grimmer(2024)]{luner2024performance}
A.~Luner and B.~Grimmer.
\newblock Performance estimation for smooth and strongly convex sets.
\newblock \emph{arXiv preprint arXiv:2410.14811}, 2024.

\bibitem[Markov and Grossmann(1916)]{markoff1916polynome}
W.~Markov and J.~Grossmann.
\newblock {\"U}ber polynome, die in einem gegebenen intervalle m{\"o}glichst wenig von null abweichen.
\newblock \emph{Mathematische Annalen}, 77\penalty0 (2):\penalty0 213--258, 1916.

\bibitem[Mones(2021)]{mones2021gentle}
L.~Mones.
\newblock A gentle introduction to optimal power flow.
\newblock \url{https://invenia.github.io/blog/2021/06/18/opf-intro/}, 2021.

\bibitem[Monga et~al.(2021)Monga, Li, and Eldar]{monga2021algorithm}
V.~Monga, Y.~Li, and Y.~C. Eldar.
\newblock Algorithm unrolling: Interpretable, efficient deep learning for signal and image processing.
\newblock \emph{IEEE Signal Processing Magazine}, 38\penalty0 (2):\penalty0 18--44, 2021.

\bibitem[Monteiro and Adler(1989)]{monteiro1989interior}
R.~D. Monteiro and I.~Adler.
\newblock Interior path following primal-dual algorithms. part i: Linear programming.
\newblock \emph{Mathematical programming}, 44\penalty0 (1):\penalty0 27--41, 1989.

\bibitem[Nesterov(1983)]{nesterov1983method}
Y.~Nesterov.
\newblock A method for solving the convex programming problem with convergence rate o (1/k2).
\newblock In \emph{Dokl akad nauk Sssr}, volume 269, page 543, 1983.

\bibitem[Nesterov and Nemirovskii(1994)]{nesterov1994interior}
Y.~Nesterov and A.~Nemirovskii.
\newblock \emph{Interior-point polynomial algorithms in convex programming}.
\newblock SIAM, 1994.

\bibitem[Nesterov et~al.(2018)]{nesterov2018lectures}
Y.~Nesterov et~al.
\newblock \emph{Lectures on convex optimization}, volume 137.
\newblock Springer, 2018.

\bibitem[Nishihara et~al.(2015)Nishihara, Lessard, Recht, Packard, and Jordan]{nishihara2015general}
R.~Nishihara, L.~Lessard, B.~Recht, A.~Packard, and M.~Jordan.
\newblock A general analysis of the convergence of admm.
\newblock In \emph{International conference on machine learning}, pages 343--352. PMLR, 2015.

\bibitem[O'Donoghue(2021)]{o2021operator}
B.~O'Donoghue.
\newblock Operator splitting for a homogeneous embedding of the linear complementarity problem.
\newblock \emph{SIAM Journal on Optimization}, 31\penalty0 (3):\penalty0 1999--2023, 2021.

\bibitem[Oymak(2021)]{oymak2021provable}
S.~Oymak.
\newblock Provable super-convergence with a large cyclical learning rate.
\newblock \emph{IEEE Signal Processing Letters}, 28:\penalty0 1645--1649, 2021.

\bibitem[O’donoghue et~al.(2016)O’donoghue, Chu, Parikh, and Boyd]{o2016conic}
B.~O’donoghue, E.~Chu, N.~Parikh, and S.~Boyd.
\newblock Conic optimization via operator splitting and homogeneous self-dual embedding.
\newblock \emph{Journal of Optimization Theory and Applications}, 169:\penalty0 1042--1068, 2016.

\bibitem[Pan(2021)]{pan2021deepopf}
X.~Pan.
\newblock Deepopf: deep neural networks for optimal power flow.
\newblock In \emph{Proceedings of the 8th ACM International Conference on Systems for Energy-Efficient Buildings, Cities, and Transportation}, pages 250--251, 2021.

\bibitem[Park and Ryu(2024)]{park2024optimal}
C.~Park and E.~K. Ryu.
\newblock Optimal first-order algorithms as a function of inequalities.
\newblock \emph{Journal of Machine Learning Research}, 25, 2024.

\bibitem[Park and Ryu(2022)]{park2022exact}
J.~Park and E.~K. Ryu.
\newblock Exact optimal accelerated complexity for fixed-point iterations.
\newblock In \emph{International Conference on Machine Learning}, pages 17420--17457. PMLR, 2022.

\bibitem[Pedregal(2004)]{pedregal2004introduction}
P.~Pedregal.
\newblock \emph{Introduction to optimization}, volume~46.
\newblock Springer, 2004.

\bibitem[Pedregosa(2020)]{pedregosa2021residual}
F.~Pedregosa.
\newblock Residual polynomials and the chebyshev method.
\newblock \url{http://fa.bianp.net/blog/2020/polyopt/}, 2020.

\bibitem[Polyak(1964)]{polyak1964some}
B.~T. Polyak.
\newblock Some methods of speeding up the convergence of iteration methods.
\newblock \emph{Ussr computational mathematics and mathematical physics}, 4\penalty0 (5):\penalty0 1--17, 1964.

\bibitem[Qian et~al.(2016)Qian, Altch{\'e}, Bender, Stiller, and de~La~Fortelle]{qian2016optimal}
X.~Qian, F.~Altch{\'e}, P.~Bender, C.~Stiller, and A.~de~La~Fortelle.
\newblock Optimal trajectory planning for autonomous driving integrating logical constraints: An miqp perspective.
\newblock In \emph{2016 IEEE 19th international conference on intelligent transportation systems (ITSC)}, pages 205--210. IEEE, 2016.

\bibitem[Renegar(1988)]{renegar1988polynomial}
J.~Renegar.
\newblock A polynomial-time algorithm, based on newton's method, for linear programming.
\newblock \emph{Mathematical programming}, 40\penalty0 (1):\penalty0 59--93, 1988.

\bibitem[Rotaru et~al.(2024)Rotaru, Glineur, and Patrinos]{rotaru2024exact}
T.~Rotaru, F.~Glineur, and P.~Patrinos.
\newblock Exact worst-case convergence rates of gradient descent: a complete analysis for all constant stepsizes over nonconvex and convex functions.
\newblock \emph{arXiv preprint arXiv:2406.17506}, 2024.

\bibitem[Ryu and V{\~u}(2020)]{ryu2020finding}
E.~K. Ryu and B.~C. V{\~u}.
\newblock Finding the forward-douglas--rachford-forward method.
\newblock \emph{Journal of Optimization Theory and Applications}, 184\penalty0 (3):\penalty0 858--876, 2020.

\bibitem[Ryu et~al.(2020)Ryu, Taylor, Bergeling, and Giselsson]{ryu2020operator}
E.~K. Ryu, A.~B. Taylor, C.~Bergeling, and P.~Giselsson.
\newblock Operator splitting performance estimation: Tight contraction factors and optimal parameter selection.
\newblock \emph{SIAM Journal on Optimization}, 30\penalty0 (3):\penalty0 2251--2271, 2020.

\bibitem[Saad(2003)]{saad2003iterative}
Y.~Saad.
\newblock \emph{Iterative methods for sparse linear systems}.
\newblock SIAM, 2003.

\bibitem[Schrijver(1998)]{schrijver1998theory}
A.~Schrijver.
\newblock \emph{Theory of linear and integer programming}.
\newblock John Wiley \& Sons, 1998.

\bibitem[Seidman et~al.(2019)Seidman, Fazlyab, Preciado, and Pappas]{seidman2019control}
J.~H. Seidman, M.~Fazlyab, V.~M. Preciado, and G.~J. Pappas.
\newblock A control-theoretic approach to analysis and parameter selection of douglas--rachford splitting.
\newblock \emph{IEEE Control Systems Letters}, 4\penalty0 (1):\penalty0 199--204, 2019.

\bibitem[Shi et~al.(2011)Shi, Razaviyayn, Luo, and He]{shi2011iteratively}
Q.~Shi, M.~Razaviyayn, Z.-Q. Luo, and C.~He.
\newblock An iteratively weighted mmse approach to distributed sum-utility maximization for a mimo interfering broadcast channel.
\newblock \emph{IEEE Transactions on Signal Processing}, 59\penalty0 (9):\penalty0 4331--4340, 2011.

\bibitem[Smith(2017)]{smith2017cyclical}
L.~N. Smith.
\newblock Cyclical learning rates for training neural networks.
\newblock In \emph{2017 IEEE winter conference on applications of computer vision (WACV)}, pages 464--472. IEEE, 2017.

\bibitem[Subosits and Gerdes(2019)]{subosits2019racetrack}
J.~K. Subosits and J.~C. Gerdes.
\newblock From the racetrack to the road: Real-time trajectory replanning for autonomous driving.
\newblock \emph{IEEE Transactions on Intelligent Vehicles}, 4\penalty0 (2):\penalty0 309--320, 2019.

\bibitem[Sun et~al.(2018)Sun, Chen, Shi, Hong, Fu, and Sidiropoulos]{sun2018learning}
H.~Sun, X.~Chen, Q.~Shi, M.~Hong, X.~Fu, and N.~D. Sidiropoulos.
\newblock Learning to optimize: Training deep neural networks for interference management.
\newblock \emph{IEEE Transactions on Signal Processing}, 66\penalty0 (20):\penalty0 5438--5453, 2018.

\bibitem[Sun et~al.(2016)Sun, Li, Xu, et~al.]{sun2016deep}
J.~Sun, H.~Li, Z.~Xu, et~al.
\newblock Deep admm-net for compressive sensing mri.
\newblock \emph{Advances in neural information processing systems}, 29, 2016.

\bibitem[Sun and Ye(2021)]{sun2021worst}
R.~Sun and Y.~Ye.
\newblock Worst-case complexity of cyclic coordinate descent: O (n\^{} 2) o (n 2) gap with randomized version.
\newblock \emph{Mathematical Programming}, 185:\penalty0 487--520, 2021.

\bibitem[Sun et~al.(2020)Sun, Luo, and Ye]{sun2020efficiency}
R.~Sun, Z.-Q. Luo, and Y.~Ye.
\newblock On the efficiency of random permutation for admm and coordinate descent.
\newblock \emph{Mathematics of Operations Research}, 45\penalty0 (1):\penalty0 233--271, 2020.

\bibitem[Sun(2020)]{sun2020optimization}
R.-Y. Sun.
\newblock Optimization for deep learning: An overview.
\newblock \emph{Journal of the Operations Research Society of China}, 8\penalty0 (2):\penalty0 249--294, 2020.

\bibitem[Sy(2023)]{sy2023optimization}
M.~Sy.
\newblock Optimization strategies for low-latency 5g nr ldpc decoding on general purpose processor.
\newblock In \emph{2023 International Conference on Control, Communication and Computing (ICCC)}, pages 1--6. IEEE, 2023.

\bibitem[Tang et~al.(2017)Tang, Dvijotham, and Low]{tang2017real}
Y.~Tang, K.~Dvijotham, and S.~Low.
\newblock Real-time optimal power flow.
\newblock \emph{IEEE Transactions on Smart Grid}, 8\penalty0 (6):\penalty0 2963--2973, 2017.

\bibitem[Taylor and Bach(2019)]{taylor2019stochastic}
A.~Taylor and F.~Bach.
\newblock Stochastic first-order methods: non-asymptotic and computer-aided analyses via potential functions.
\newblock In \emph{Conference on Learning Theory}, pages 2934--2992. PMLR, 2019.

\bibitem[Taylor and Drori(2023)]{taylor2023optimal}
A.~Taylor and Y.~Drori.
\newblock An optimal gradient method for smooth strongly convex minimization.
\newblock \emph{Mathematical Programming}, 199\penalty0 (1):\penalty0 557--594, 2023.

\bibitem[Taylor et~al.(2018{\natexlab{a}})Taylor, Van~Scoy, and Lessard]{taylor2018lyapunov}
A.~Taylor, B.~Van~Scoy, and L.~Lessard.
\newblock Lyapunov functions for first-order methods: Tight automated convergence guarantees.
\newblock In \emph{International Conference on Machine Learning}, pages 4897--4906. PMLR, 2018{\natexlab{a}}.

\bibitem[Taylor et~al.(2017{\natexlab{a}})Taylor, Hendrickx, and Glineur]{taylor2017exact}
A.~B. Taylor, J.~M. Hendrickx, and F.~Glineur.
\newblock Exact worst-case performance of first-order methods for composite convex optimization.
\newblock \emph{SIAM Journal on Optimization}, 27\penalty0 (3):\penalty0 1283--1313, 2017{\natexlab{a}}.

\bibitem[Taylor et~al.(2017{\natexlab{b}})Taylor, Hendrickx, and Glineur]{taylor2017smooth}
A.~B. Taylor, J.~M. Hendrickx, and F.~Glineur.
\newblock Smooth strongly convex interpolation and exact worst-case performance of first-order methods.
\newblock \emph{Mathematical Programming}, 161:\penalty0 307--345, 2017{\natexlab{b}}.

\bibitem[Taylor et~al.(2018{\natexlab{b}})Taylor, Hendrickx, and Glineur]{taylor2018exact}
A.~B. Taylor, J.~M. Hendrickx, and F.~Glineur.
\newblock Exact worst-case convergence rates of the proximal gradient method for composite convex minimization.
\newblock \emph{Journal of Optimization Theory and Applications}, 178:\penalty0 455--476, 2018{\natexlab{b}}.

\bibitem[Teboulle and Vaisbourd(2023)]{teboulle2023elementary}
M.~Teboulle and Y.~Vaisbourd.
\newblock An elementary approach to tight worst case complexity analysis of gradient based methods.
\newblock \emph{Mathematical Programming}, 201\penalty0 (1):\penalty0 63--96, 2023.

\bibitem[Toupet(2006)]{toupet2006real}
O.~Toupet.
\newblock \emph{Real-time path-planning using mixed-integer linear programming and global cost-to-go maps}.
\newblock PhD thesis, Massachusetts Institute of Technology, 2006.

\bibitem[Tsallis(1988)]{tsallis1988possible}
C.~Tsallis.
\newblock Possible generalization of boltzmann-gibbs statistics.
\newblock \emph{Journal of statistical physics}, 52:\penalty0 479--487, 1988.

\bibitem[Van~Scoy et~al.(2017)Van~Scoy, Freeman, and Lynch]{van2017fastest}
B.~Van~Scoy, R.~A. Freeman, and K.~M. Lynch.
\newblock The fastest known globally convergent first-order method for minimizing strongly convex functions.
\newblock \emph{IEEE Control Systems Letters}, 2\penalty0 (1):\penalty0 49--54, 2017.

\bibitem[Wang et~al.(2024)Wang, Ma, Yang, and Zhou]{wang2024relaxed}
B.~Wang, S.~Ma, J.~Yang, and D.~Zhou.
\newblock Relaxed proximal point algorithm: Tight complexity bounds and acceleration without momentum.
\newblock \emph{arXiv preprint arXiv:2410.08890}, 2024.

\bibitem[Wang et~al.(2019)Wang, Fazlyab, Chen, and Preciado]{wang2019robust}
H.~Wang, M.~Fazlyab, S.~Chen, and V.~M. Preciado.
\newblock Robust convergence analysis of three-operator splitting.
\newblock In \emph{2019 57th Annual Allerton Conference on Communication, Control, and Computing (Allerton)}, pages 391--398. IEEE, 2019.

\bibitem[Wen et~al.(2024)Wen, Li, Wang, Hall, Liang, and Ma]{wen2024understanding}
K.~Wen, Z.~Li, J.~Wang, D.~Hall, P.~Liang, and T.~Ma.
\newblock Understanding warmup-stable-decay learning rates: A river valley loss landscape perspective.
\newblock \emph{arXiv preprint arXiv:2410.05192}, 2024.

\bibitem[Wright(1997)]{wright1997primal}
S.~J. Wright.
\newblock \emph{Primal-dual interior-point methods}.
\newblock SIAM, 1997.

\bibitem[Wu et~al.(2024)Wu, Bartlett, Telgarsky, and Yu]{wu2024large}
J.~Wu, P.~L. Bartlett, M.~Telgarsky, and B.~Yu.
\newblock Large stepsize gradient descent for logistic loss: Non-monotonicity of the loss improves optimization efficiency.
\newblock \emph{arXiv preprint arXiv:2402.15926}, 2024.

\bibitem[Xiang et~al.(2013)Xiang, Gubian, Suomela, and Hoeng]{xiang2013generalized}
Y.~Xiang, S.~Gubian, B.~Suomela, and J.~Hoeng.
\newblock Generalized simulated annealing for global optimization: the gensa package.
\newblock \emph{R J.}, 5\penalty0 (1):\penalty0 13, 2013.

\bibitem[Xiong(2024)]{xiong2024accessible}
Z.~Xiong.
\newblock Accessible theoretical complexity of the restarted primal-dual hybrid gradient method for linear programs with unique optima.
\newblock \emph{arXiv preprint arXiv:2410.04043}, 2024.

\bibitem[Xiong and Freund(2023)]{xiong2023computational}
Z.~Xiong and R.~M. Freund.
\newblock Computational guarantees for restarted pdhg for lp based on" limiting error ratios" and lp sharpness.
\newblock \emph{arXiv preprint arXiv:2312.14774}, 2023.

\bibitem[Xiong and Freund(2024)]{xiong2024role}
Z.~Xiong and R.~M. Freund.
\newblock The role of level-set geometry on the performance of pdhg for conic linear optimization.
\newblock \emph{arXiv preprint arXiv:2406.01942}, 2024.

\bibitem[Yang et~al.(2018)Yang, Sun, Li, and Xu]{yang2018admm}
Y.~Yang, J.~Sun, H.~Li, and Z.~Xu.
\newblock Admm-csnet: A deep learning approach for image compressive sensing.
\newblock \emph{IEEE transactions on pattern analysis and machine intelligence}, 42\penalty0 (3):\penalty0 521--538, 2018.

\bibitem[Ye(2011)]{ye2011interior}
Y.~Ye.
\newblock \emph{Interior point algorithms: theory and analysis}.
\newblock John Wiley \& Sons, 2011.

\bibitem[Yoon and Ryu(2021)]{yoon2021accelerated}
T.~Yoon and E.~K. Ryu.
\newblock Accelerated algorithms for smooth convex-concave minimax problems with o (1/k\^{} 2) rate on squared gradient norm.
\newblock In \emph{International Conference on Machine Learning}, pages 12098--12109. PMLR, 2021.

\bibitem[Young(1953)]{young1953richardson}
D.~Young.
\newblock On richardson's method for solving linear systems with positive definite matrices.
\newblock \emph{Journal of Mathematics and Physics}, 32\penalty0 (1-4):\penalty0 243--255, 1953.

\bibitem[Young(2014)]{young2014iterative}
D.~M. Young.
\newblock \emph{Iterative solution of large linear systems}.
\newblock Elsevier, 2014.

\bibitem[Yu and Fan(2024)]{yu2024rigid}
M.~Yu and C.~Fan.
\newblock Rigid body path planning using mixed-integer linear programming.
\newblock \emph{IEEE Robotics and Automation Letters}, 2024.

\bibitem[Zhang et~al.(2021)Zhang, Chen, and Zhang]{zhang2021convex}
L.~Zhang, Y.~Chen, and B.~Zhang.
\newblock A convex neural network solver for dcopf with generalization guarantees.
\newblock \emph{IEEE Transactions on Control of Network Systems}, 9\penalty0 (2):\penalty0 719--730, 2021.

\bibitem[Zhang and Jiang(2024)]{zhang2024accelerated}
Z.~Zhang and R.~Jiang.
\newblock Accelerated gradient descent by concatenation of stepsize schedules.
\newblock \emph{arXiv preprint arXiv:2410.12395}, 2024.

\bibitem[Zhang et~al.(2024)Zhang, Lee, Du, and Chen]{zhang2024anytime}
Z.~Zhang, J.~D. Lee, S.~S. Du, and Y.~Chen.
\newblock Anytime acceleration of gradient descent.
\newblock \emph{arXiv preprint arXiv:2411.17668}, 2024.

\bibitem[Zhou et~al.(2022)Zhou, Tian, So, and Cheng]{zhou2022practical}
K.~Zhou, L.~Tian, A.~M.-C. So, and J.~Cheng.
\newblock Practical schemes for finding near-stationary points of convex finite-sums.
\newblock In \emph{International Conference on Artificial Intelligence and Statistics}, pages 3684--3708. PMLR, 2022.

\end{thebibliography}

\appendix
\onecolumn

\section{More Results}
\label{appendix_more_results}

\paragraph{More results related to Figure~\ref{fig_intro} (b).} Figure~\ref{fig_intro} (b)presents the curves of optimality gap and the detailed stepsize rules. We find that the optimal constant stepsize cannot reach optimal solutions until $>200$ steps, while the Finite Horizon stepsize rule solves the problem in 2 steps. The experimental setup is the same as in Section~\ref{sec_toy_example} except that we change the initialization to $N(80,1)$ for better visualization.

\begin{figure}[h]
    \centering
\subfigure[Optimality gap]{\includegraphics[width=0.30\textwidth]{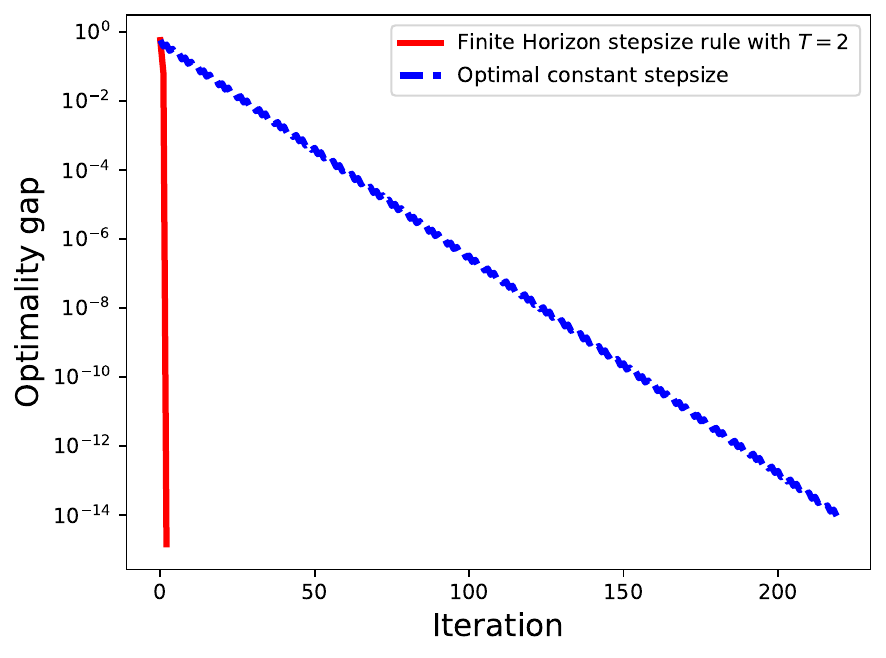}}
\subfigure[Stepsize rule]
{\includegraphics[width=0.30\textwidth]{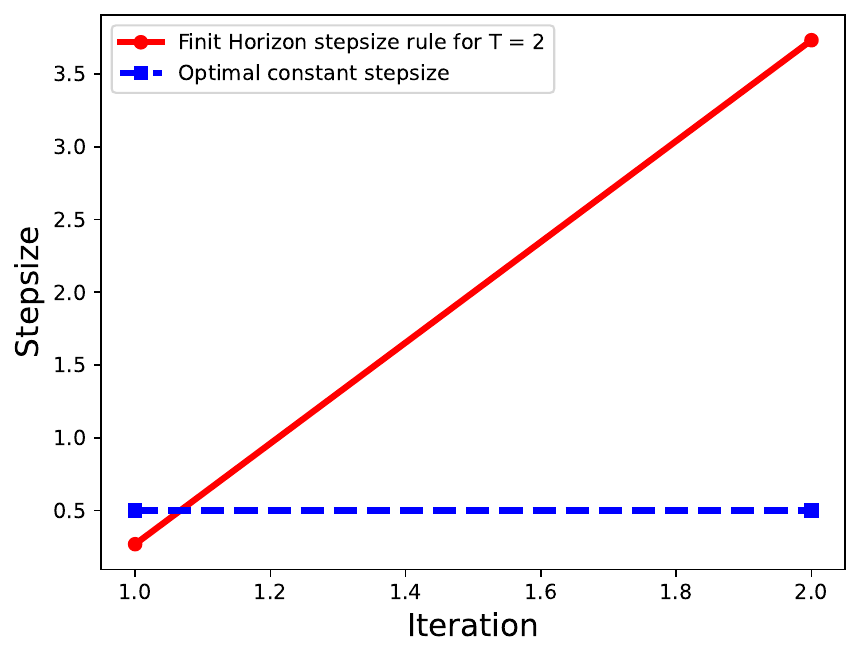}}
    \caption{ The curves of optimality gap and the detailed stepsize rules (for the first 2 steps) of the results in Figure~\ref{fig_intro}.  }
  \label{fig_landscape_more}
  \vspace{-0.3cm}
\end{figure}

\end{document}